\documentclass[aos]{imsart}

\usepackage[utf8]{inputenc}
\usepackage{xr}

\makeatletter
\newcommand*{\addFileDependency}[1]{
  \typeout{(#1)}
  \@addtofilelist{#1}
  \IfFileExists{#1}{}{\typeout{No file #1.}}
}
\makeatother

\newcommand*{\suppl}[1]{%
    \externaldocument{#1}%
    \addFileDependency{#1.tex}%
    \addFileDependency{#1.aux}%
}

\suppl{rockova_rousseau_supp}

\RequirePackage[OT1]{fontenc}
\RequirePackage{amsthm,amsmath}
\RequirePackage[numbers]{natbib}
\RequirePackage[colorlinks,citecolor=blue,urlcolor=blue]{hyperref}

\usepackage{xr-hyper}
\usepackage{amsmath}
\usepackage{amsthm}
\usepackage{color}
\usepackage{amsfonts}
\usepackage{amssymb}
\usepackage{graphicx}
\usepackage{subfigure}
\usepackage{multirow}
\usepackage{natbib}
\usepackage{color}
\usepackage{natbib}
\usepackage{rotating}
\usepackage{mathtools}
\usepackage[all,cmtip]{xy}
\usepackage{amsmath}

\usepackage{accents}

\newcommand\munderbar[1]{
  \underaccent{\bar}{#1}}

\usepackage[all,cmtip]{xy}

\newcommand{\iid}{\stackrel{iid}{\sim}}

\newtheorem{definition}{Definition}
\newtheorem{lemma}{Lemma}
\newtheorem{theorem}{Theorem}
\newtheorem{example}{Example}
\newtheorem{remark}{Remark}

\usepackage[dvipsnames]{xcolor}

\definecolor{pink1}{rgb}{0.858, 0.188, 0.478}

\usepackage{graphicx}
\usepackage{subfigure}
\usepackage{multirow}
\usepackage{color}
\usepackage{natbib}
\usepackage{stmaryrd}
\usepackage{rotating}

\newcommand{\bm}[1]{\boldsymbol{#1}}

\newcommand{\mX}{\mathcal{X}}

\newcommand{\wtLmax}{\wt L_{max}}

\newcommand{\cc}{\mathcal{C}}

\newcommand{\Lmax}{L_{max}}

\newcommand{\La}{\Lambda}

\newcommand{\R}{\mathbb{R}}
\newcommand{\N}{\mathbb{N}}
\newcommand{\bmT}{\backslash \mT}
\newcommand{\I}{\mathbb{I}}

\newcommand{\Y}{\bm{Y}}
\newcommand{\mE}{\mathcal{E}}
\newcommand{\given}{ \,|\,}
\newcommand{\ix}{\mathbb{X}}
\newcommand{\mL}{\mathcal{L}}
\newcommand{\cA}{\mathcal{A}}
\newcommand{\mA}{\mathcal{A}}
\newcommand{\mC}{\mathcal{C}}

\newcommand{\M}{\mathcal{M}}

\renewcommand{\1}{\mathbb{I}}

\newcommand{\wh}[1]{\widehat{#1}}

\newtheorem{assumption}{Assumption}[section]

\def\C {\,|\:}

\def\e{\mathrm{e}}
\def\d{\mathrm{d\,}}

\def\bT{\mathbb{T}}
\def\cT{\mathcal{T}}
\def\mT{\mathcal{T}}
\def\b{\bm{\beta}}
\def\wt{\widetilde}

\def\be {\bm{\beta}}

\usepackage{amssymb,tikz}

\numberwithin{equation}{section}

\begin{document}

\begin{frontmatter}
\title{Ideal Bayesian  Spatial Adaptation}
\runtitle{Bayesian Trees are Spatially Adaptive}
\begin{aug}

\author{\fnms{Veronika} \snm{Ro\v{c}kov\'{a} }\thanksref{t2,m2}\ead[label=e2]{veronika.rockova@chicagobooth.edu}}
\author{and \fnms{Judith} \snm{Rousseau}\thanksref{m1}\ead[label=e1]{judith.rousseau@stats.ox.ac.uk}}

\thankstext{t2}{
The author gratefully acknowledges  support from the James S. Kemper Foundation Faculty Research Fund at  the University of Chicago Booth School of Business and the National Science Foundation (DMS:1944740).}


\affiliation{University of Chicago\thanksmark{m2}}
\affiliation{University of Oxford \thanksmark{m1}}

\address{University of Chicago\\ Booth School of Business\\ 5807 S. Woodlawn Avenue\\ Chicago, IL, 60637 \\ USA\\
 \printead{e2}}

\address{University of Oxford \\
Department of Statistics\\
 24-29 St Giles'\\ Oxford OX1 3LB \\ United Kingdom\\
  \printead{e1}}

\end{aug}

\begin{abstract}
Many real-life applications involve estimation of curves  that exhibit complicated shapes including jumps or varying-frequency oscillations. 
Practical methods have been devised that can adapt to a locally varying complexity of an unknown function (e.g. variable-knot splines, sparse wavelet reconstructions, kernel methods or trees/forests). However, the overwhelming majority of  existing asymptotic minimaxity  theory  is predicated on homogeneous smoothness  assumptions. Focusing on locally H\"{o}lderian functions, we provide new {\em locally adaptive} posterior concentration rate results under the supremum loss for widely used Bayesian machine learning techniques in white noise and non-parametric regression. In particular, we show that popular spike-and-slab priors and Bayesian CART  are uniformly locally adaptive. 
In addition, we propose a new class of repulsive partitioning priors which relate to variable knot splines and which are exact-rate adaptive.
 For uncertainty quantification, we construct locally adaptive confidence bands whose width depends on the local smoothness and which achieve uniform asymptotic coverage under local self-similarity. 
To illustrate that spatial adaptation is not at all automatic, we provide  lower-bound results showing  that popular hierarchical Gaussian process priors fall short of spatial adaptation.

\end{abstract}



\end{frontmatter}


\section{Spatial Adaptation}\label{sec:intro}
The key to practically successful curve estimation  is the ability to adapt to subtle qualitative structures  of the analyzed curve.
Very often, interesting aspects of the estimated curve are related to  spatial inhomogeneities, e.g. discontinuities or oscillations with varying frequency and/or amplitude. There is a wealth of techniques for   function estimation   (e.g. kernel methods, local polynomial fitting, nearest neighbor techniques or splines) which exert various degrees of global and local adaptation. For functions with a {\em locally} varying complexity, however,  global  procedures can be woefully inaccurate,  leading to overfitting   smooth domains and underfitting   wiggly domains.

Such weaknesses have long been recognized. Numerous methodological developments have spawned that are capable of local adaptation, e.g.,  
local polynomial regression   \citep{gijbels} or  kernel estimation \citep{muller_stadtmuller,lepski_etal,breiman_local} with a local bandwidth selection. In the context of spline smoothing, \cite{luo_wahba} suggested  adaptively selecting subsets of basis functions which
  pertains  to selective wavelet reconstructions \citep{donoho94,chipman_wavelet} and variable knot  spline techniques  \cite{friedman_silverman,stone,mars,variable_knot}. Notably, \cite{mammen_geer,tibshirani_trend} proposed  (total-variation) penalized least square estimates which correspond to regression splines with data-adaptive knot points.
  An alternative approach is to allow the smoothing parameter to vary locally (see \cite{pintore} for piecewise constant smoothing parameters). 
For example, \cite{ruppert_carroll} suggest   spline fitting with a roughness penalty whose logarithm   is itself a linear spline   with knot values   chosen  by  cross-validation.  
 Variants of such spatially adaptive penalty parameters have  been  widely used in practice \citep{lang_bretzger, krivobokova, craina,bala}.
 Besides splines and wavelets, tree based methods (CART \citep{breiman_book,cart1,cart2}, random forests \citep{breiman2} or BART \citep{bart}) are particularly appealing for recovering spatially inhomogeneous functions by adapting the placement of splits to the function itself via recursive partitioning \citep{ideal_spatial}.  Deep learning methods are  also expected to perform well  for  structured curves \citep{suzuki,hayakawa_suzuki}.

From a methodological perspective,   spatially adaptive curve estimation has been tackled rather broadly. From a theoretical perspective, however, there are still gaps in understanding whether these techniques are indeed optimal and (uniformly) spatially adaptive. The majority of existing asymptotic minimaxity  theory (for density or regression function estimation)  is  predicated on homogeneous smoothness  assumptions. 
For example, existing results for random forests \citep{wager_athey,wager}, deep learning \citep{schmidt,polson_rockova}, 
Bayesian forests \citep{rockova_vdp,seonghyun} and other non-parametric methods such as Gaussian processes \citep{vdv_zanten,gp_pati} have been  concerned with  convergence rates for spatially homogeneous H\"{o}lderian functions  under the $L^2$ global estimation loss.  Here, we extend the scope of such theoretical results in two important directions. First, we focus  on {\em both} global and local (supremum) loss providing results for uniform local adaptation. Second, we provide a frequentist framework for  uncertainty quantification   in the form of locally adaptive bands. Our goal is to investigate the extent to which widely used Bayesian priors (spike-and-slab priors \cite{chipman_wavelet,ray,hoffmann,yoo}, Gaussian process priors 
\cite{arbel,shen_ghosal,belitzer_ghosal,knapik,szabo_vdv_vz}  and Bayesian CART priors \cite{cart1,cart2,castillo_rockova}) can  optimally adapt to local smoothness. 
Before listing our contributions, we  review existing theoretical results for  spatially inhomogeneous functions.

The first natural question is  how well an estimator performs {\em globally}. For the stereotypical Besov  classes $B_{p,q}^\alpha$, one way to assess the global quality of an estimator is in terms of a  $L^{r}$ loss  
that is sharper than the norm of the Besov functional class (i.e. $p<r<\infty$), see e.g. \cite{donoho98}  and \cite{lepski_etal}. For $p<2$, linear estimators are known to be incapable of achieving the optimal rate \citep{donoho98}. 
For a discussion on minimax rates in Besov spaces, we refer to \cite{asymptopia} and \cite{delyon_juditsky}.
Unlike linear estimators, wavelet thresholding  offers a powerful technology for spatially adaptive curve estimation \citep{asymptopia}. In particular,
\cite{donoho94} describe a selective wavelet reconstruction method called RiskShrink based on shrinkage of wavelet coefficients and show that this procedure mimics an oracle `as well as it is possible to do so'.
RiskShrink is an automatic model selection method which picks a subset of wavelet vectors and fits a model consisting only of wavelets in that subset. In this work, we investigate Bayesian variants of such strategies.
Positive findings for global estimation in Besov spaces have also been reported for deep learning  \citep{suzuki,hayakawa_suzuki}, 
penalized least squares  \citep{mammen_geer}, locally variable kernels \citep{muller_stadtmuller}.
Notably, \cite{lepski_etal} propose  a kernel estimator with a variable data-driven bandwidth  that achieves the minimax rate of estimation over a wide scale of Besov classes and hence shares rate optimality properties with wavelet estimators.  

Another, and perhaps more transparent, approach is to assess the quality of an estimator {\em locally}.
For density estimation, \cite{gach_etal}  study adaptivity to heterogeneous smoothness, simultaneously for every point  in a fixed interval, in a supremum-norm loss. 
 The authors consider a certain notion of pointwise H\"{o}lder continuity and study dyadic histogram estimators with a variable bin size  and with  a Lepski-type adaptation.
We adopt a similar estimation setup here, but approach it entirely from a Bayesian perspective. 

Practical deployments of the Lepski-based adaptation require tuning parameters (especially of the threshold used for comparing two estimates from different scales) for which
theoretical justifications may not always be available \cite{gach_etal,kathkovnik}. Bayesian procedures, on the other hand, are known to adapt automatically to the unknown aspects of the estimation problem, even yielding rate-exact adaptation \cite{hoffmann}.  This work  studies  whether (rate-exact) uniform adaptation is attainable for popular Bayesian learning procedures in terms of local (supremum-norm) concentration rates. We are not aware of any other Bayesian investigation of this type in the literature. Our contributions can be grouped into four types of results.
First, we show that spike-and-slab priors achieve uniform exact-rate optimal adaptation   in a supremum-norm sense under the white noise model. We relax the prior assumptions of \cite{hoffmann}, allowing for considerably less sparse priors. Next, building on \cite{castillo_rockova} we show that Bayesian CART is {\em also} uniformly locally adaptive but sacrifices a logarithmic factor. These results are obtained in the white-noise model as well as non-parametric regression with suitably regular (not necessarily equi-distant) fixed design points. Second, we show how to construct locally adaptive credible bands (with asymptotic coverage $1$) whose width depends on local smoothness. Third, we provide negative results   for Gaussian process priors showing that they are incapable of local adaptation. Fourth, in the context of non-parametric regression, we propose a new class of `repulsive' partitioning priors which penalize irregular partitions and which are locally rate-exact. These priors can be viewed as a simplified (zero-degree)  version of data-adaptive knot splines. Our results thus provide a stepping stone towards studying Bayesian variable-knot spline techniques.

The manuscript is organized as follows. Section \ref{sec:setup} describes the estimation setup and reviews some facts about spatially inhomogeneous functions. Section \ref{sec:wn} shows results for spike-and-slab and Bayesian CART priors in the white noise model. Section \ref{sec:npreg} then shows  our results for non-parametric regression. Section \ref{sec:discussion} wraps up with a discussion and Section \ref{sec:proofs} shows  proofs of two of our main theorems. The rest of the proofs is in the Supplemental Materials.

\section{Statistical Setting}\label{sec:setup}
For our theoretical development, we will consider {\em both}  the non-parametric regression  model as well as its idealized  white noise counterpart. 
The regression model assumes  $n$ noisy samples  $Y=(Y_1,\dots, Y_n)'$ of a function $f_0:[0,1]\rightarrow\R$, where 
  \begin{equation}\label{model}
   Y_i = f_0(x_i) +\epsilon_i, \quad \epsilon_i \stackrel{iid}{\sim} \mathcal N(0, \sigma^2),
   \end{equation}
with $x_i \in[0,1]$ and where $\sigma^2>0$ is a known scalar.  
The white noise model is an elegant continuous version of \eqref{model} defined via a stochastic differential equation
\begin{equation}\label{model2}
dY(t)=f_0(t)dt+\frac{1}{\sqrt n}dW(t),\quad t\in[0,1],\quad n\in\N,
\end{equation}
where $Y(t)$ is an observation process,  $W(t)$ is the standard Wiener process on $[0,1]$ and $f_0\in L^2[0,1]$ is an unknown bounded function on $[0,1]$ to be estimated.  
The model \eqref{model2} 
is observationally equivalent to a Gaussian sequence   model after projecting the observation process onto a wavelet basis $\{\psi_{lk}:l\geq0,0\leq k\leq 2^{l}-1\}$ of $L^2[0,1]$. This sequence model writes 
as 
\begin{equation}\label{eq:model2}
Y_{lk}=\beta_{lk}^0+\frac{\varepsilon_{lk}}{\sqrt n},\quad\varepsilon_{lk}\iid\mathcal{N}(0,1),\quad l\geq 0,\quad k=0,\dots, 2^l-1,
\end{equation}
 where  the wavelet coefficients $\beta_{lk}^0=\langle f_0,\psi_{lk}\rangle= \int_0^1f_0(t)\psi_{lk}(t)dt$ of $f_0$ are indexed by a scale parameter $l$ and a location parameter $k$. 
Throughout this work, we will be using  the standard Haar wavelet basis 
\begin{equation}\label{haar}
\psi_{-10}(x)=\I_{[0,1]}(x)\quad\text{and}\quad\psi_{lk}(x)=2^{l/2}\psi(2^lx-k)
\end{equation}
obtained with orthonormal dilation-translations of  $\psi=\I_{(0,1/2]}-\I_{(1/2,1]}$.
We denote with $I_{lk}$ the {\em dyadic}  intervals  which correspond to the domain of the balanced Haar wavelets  $\psi_{lk}$ in \eqref{haar}, i.e.
\begin{equation}\label{eq:domain}
 I_{00}=(0,1],\quad I_{lk}=(k2^{-l},(k+1)2^{-l}]\quad\text{for $l\ge 0$ and $0\le k<2^l$}.
 \end{equation}

\begin{figure}[!t]
\scalebox{0.37}{\includegraphics{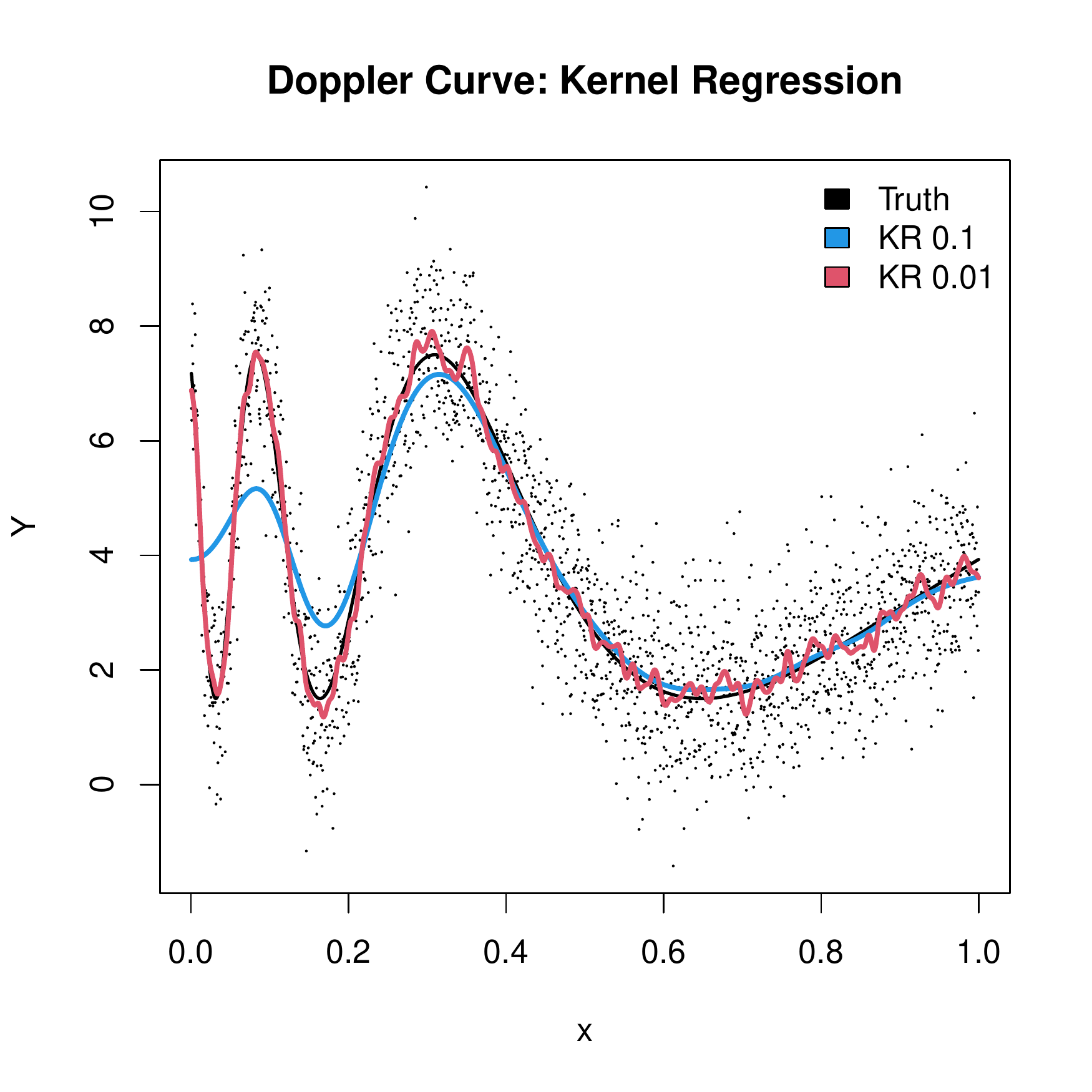}}
\scalebox{0.37}{\includegraphics{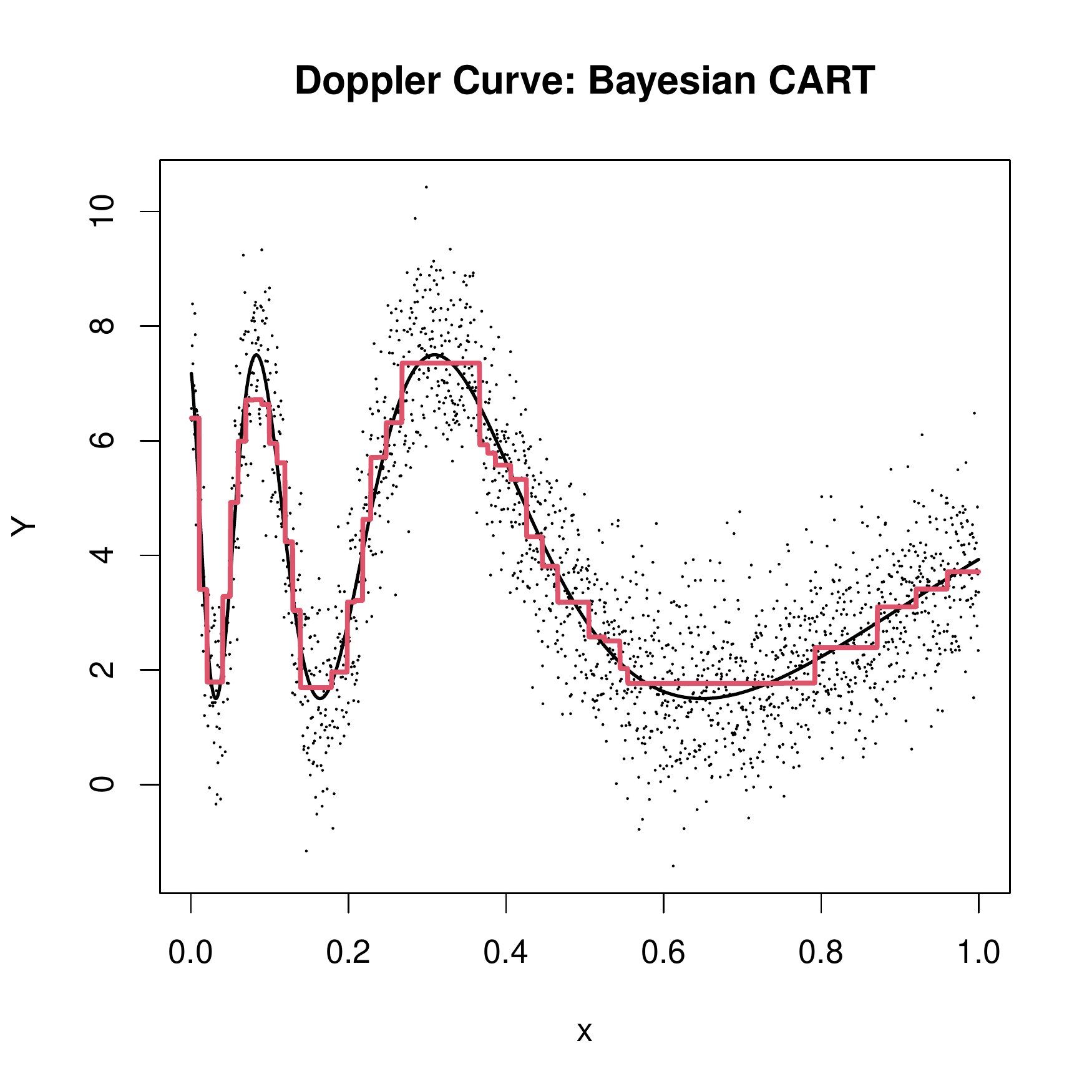}}
\caption{Doppler curve and (left) kernel regression estimates (\texttt{ksmooth} in R) with a  bandwidth $0.1$ (capturing well the flat part)  and $0.01$ (capturing  well the wiggly part), (right) Bayesian CART posterior mean fit (\texttt{wbart} in the R package \texttt{BART}).}\label{doppler1}
\end{figure}

Brown and Low \cite{brown_low} showed asymptotic equivalence between \eqref{model} (with an equispaced fixed design) and \eqref{model2}  under a uniform smoothness assumption which is satisfied by, e.g., $\alpha$-H\"{o}lderian smooth  functions with $\alpha>0.5$. From a sequence of optimal procedures in one problem, they also prescribed a construction of an asymptotically equivalent sequence in the other. This recipe is particularly convenient for linear estimators. For Bayesian methods, however, it is generally {\em not known} whether the knowledge of a (wavelet shrinkage/non-linear) minimax procedure in one problem automatically implies the optimality in the other. 
This is why we study Bayesian procedures in {\em both} models (see Section \ref{sec:wn} for white noise and Section \ref{sec:npreg} for regression). We obtain local rate-optimality in regression under the assumption $\alpha>1/2$ in Section \ref{sec:np1} and, finally,  in Section \ref{sec:repulsive} we show that a new class of adaptive-split priors (related to variable-knot spline techniques) yields exact-rate optimality {\em without assuming}  $\alpha>1/2$.

In both models \eqref{model} and \eqref{model2}, the goal is to estimate a  possibly {\em spatially inhomogeneous} function $f_0$ (see Section \ref{sec:spatial_functions} below). We assess the quality of an estimator using both the global $L^r$ loss as well as the locally re-weighted supremum-norm loss (similarly as in \cite{gach_etal}). In particular, for the (near) minimax rate  $r_n(x)$ of adaptive estimation of $f_0$ at the point $x$, we will show that, with probability $P_{f_0}$ tending to one, the random variable 
\begin{equation}\label{eq:sup_norm}
 \sup_{x\in[0,1]}\frac{1}{r_n(x)}|f(x)-f_0(x)|,\quad\text{with}\quad f\sim \Pi(f\C Y)
\end{equation}
denoting the posterior distribution, is stochastically bounded, thereby implying a uniform local adaptation. Such spatial adaptation is not automatic for many standard estimators.
We illustrate this phenomenon on two  simulated examples below.

\begin{example}(Doppler Curve)
Similarly as in \cite{asymptopia} and  \cite{tibshirani_trend}, we generate $n=2048$ observations from a noisy Doppler curve \eqref{model} with 
$f_0(x)=3\sin[4/(x+0.2)]+1.5$ and $\sigma=1$ with $x_i=i/n$. This function has  heterogeneous smoothness which cannot be captured with prototypical global smoothing methods such as global kernel regression   (Figure  \ref{doppler1} on the left) which leads to over/undersmoothing depending on the choice of a fixed kernel width. Tree-based methods, such as Bayesian CART, are better suited for this task by placing the splits more often in areas where the function is less smooth  (Figure  \ref{doppler1} on the right). 
\end{example}

\begin{figure}[!t]
\scalebox{0.37}{\includegraphics{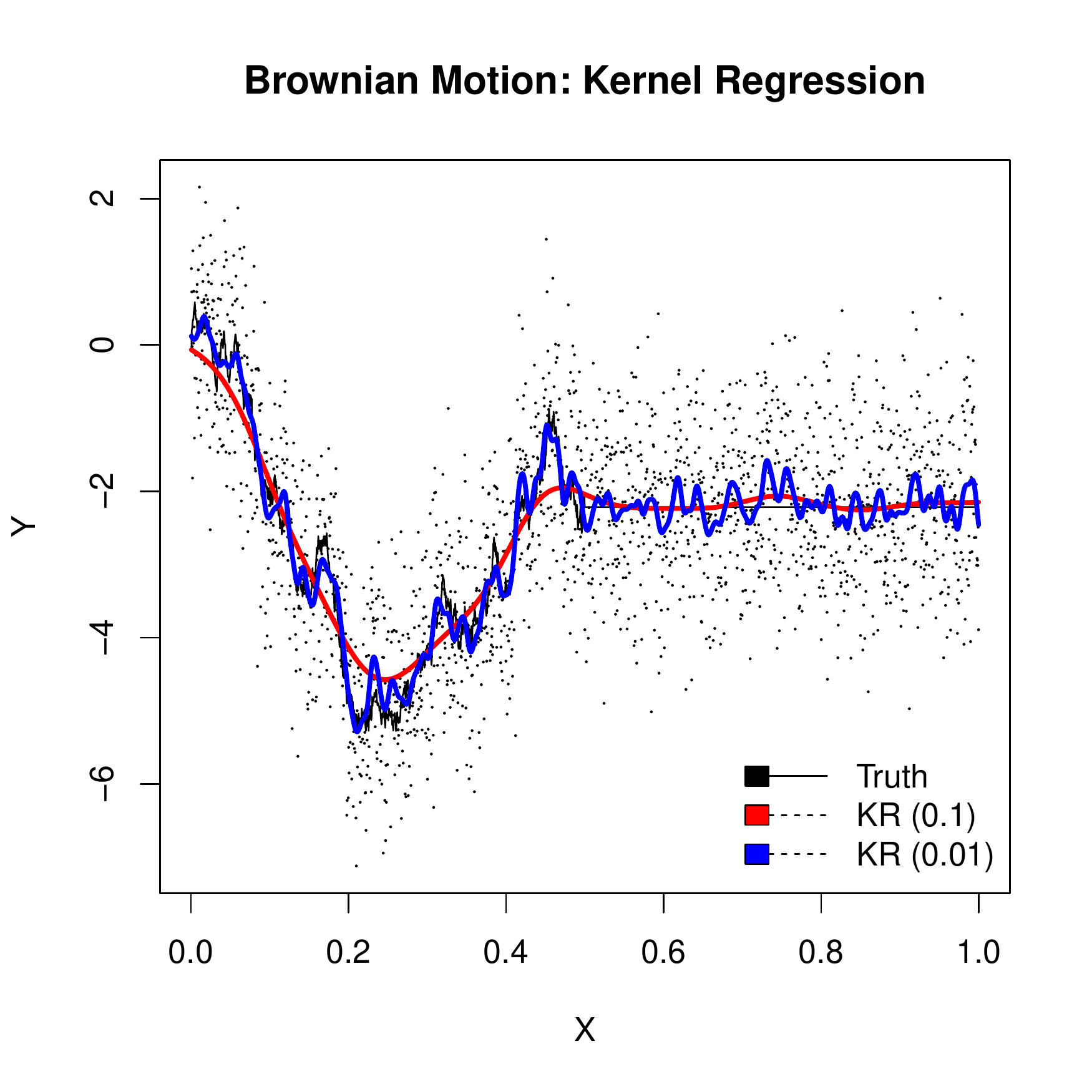}}
\scalebox{0.37}{\includegraphics{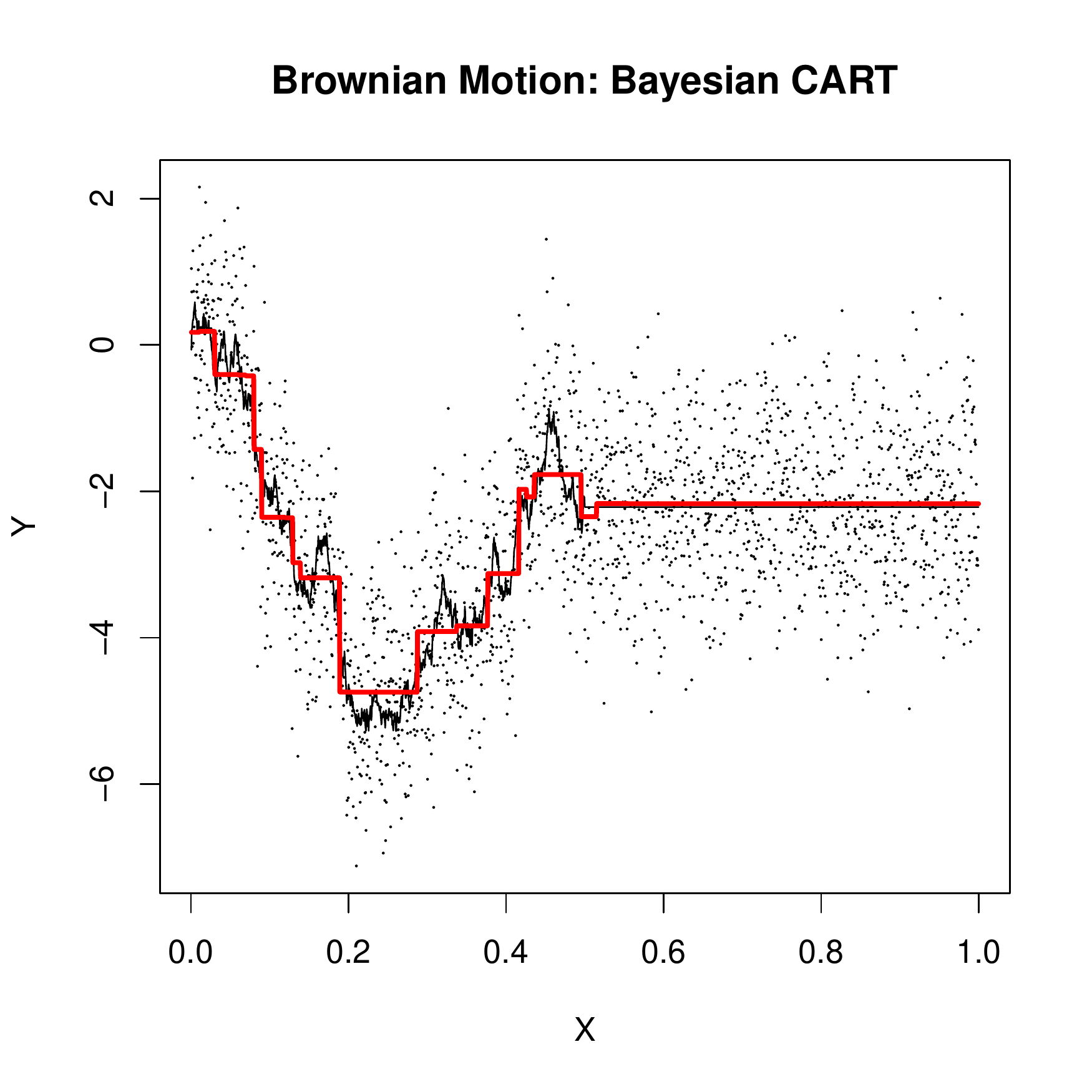}}
\caption{The Brownian motion example  with  (left) kernel regression estimates  (\texttt{ksmooth} in R)  with a  bandwidth $0.1$ (capturing well the flat part of the curve)  and $0.01$ (capturing  well the wigglier part of the curve) and (right) Bayesian CART (posterior mean fit).}
\label{fig:bm}
\end{figure}

\begin{example}(Brownian Motion)\label{ex:bm}
The second example\footnote{We are thankful to Isma\"{e}l Castillo for suggesting this example.} assumes  $n=2^{11}$ where $f_0$ has been generated from a Brownian motion on $[0,1/2)$ (whose almost all trajectories are locally $\alpha$-H\"{o}lder continuous with $\alpha<1/2$) and a constant function on $[1/2,1]$. 
The plots of the kernel regression and Bayesian CART estimates are in Figure \ref{fig:bm}. Bayesian CART wastes  no splits on the flat domain, showcasing its spatial adaptivity. 
We will investigate this example theoretically in Section \ref{sec:LBGaussian} where we show that hierarchical Gaussian processes adapt to the worse regularity (determined by the Brownian motion). Beyond adaptability, Bayesian methods can also quantify uncertainty via the posterior (as seen from Figure \ref{bm} in Section \ref{sec:uq}). The width of the optimal band should be wider when the function is less smooth. 
In Section \ref{sec:uq}, we propose one such construction and show its frequentist validity.
\end{example}

\vspace{-0.5cm}
\subsection{Spatially Inhomogeneous Functions}\label{sec:spatial_functions}
Below, we review several known facts about function classes with  inhomogeneous smoothness.
The Besov class $B_{p,q}^\alpha$   (which contains   H\"{o}lder and Sobolev classes by setting $p=q=\infty$ and $p=q=2$, respectively)   permits spatial inhomogeneity when $p<2$. For example, the  Bump algebra (consisting of infinite mixtures of Gaussian bumps) coincides with  $B_{1,1}^1$  \citep{donoho98} and  constitutes an interesting caricature of  smoothness inhomogeneity which would {\em not} be allowed within the  H\"{o}lder class. Another example is the total variation (TV) class  (contained inside $B_{1,\infty}^1$ and containing $B_{1,1}^1$) which includes functions that may have jumps localized in one part of the domain and be very flat elsewhere  \citep{donoho98}. 
For a discussion on global minimax rates in Besov spaces we refer to \cite{donoho_johnstone95,delyon_juditsky}.

Function spaces where the smoothness can vary from point to point have  quite a rich history.  Besov spaces with variable smoothness were defined by \cite{leopold} and later developed by  many others
(see \cite{schneider} and references therein). We focus on H\"{o}lderian  functions which are more intuitive for  a supremum-norm analysis.
Indeed, the perhaps more widely accepted notion of pointwise regularity has been formalized for H\"{o}lderian functions where the  exponent itself is a function taking its values in $[0,\infty)$  \citep{jaffard,anderson}. 
The question of which functions may be H\"{o}lder exponents was   raised and partially answered in \cite{daoudi}. 
Andersson \cite{anderson} later showed that  a non-negative function is an exponent of a pointwise H\"{o}lder function if and only if it can be written as a limit inferior of a sequence of continuous functions.
For example, 'typical' functions in the Besov space exhibit a multifractal behavior where the H\"{o}lder exponent is a continuous function \citep{jaffard,gach_etal}. 

Following \cite{gach_etal}, we define a set of bounded functions that are locally $t$-H\"{o}lder at $x\in \mathcal [0,1]$
\begin{equation}\label{eq:local_holder}
\mC(t,x,M,\eta)=\left\{f:[0,1] \rightarrow\R;\max\left(\|f\|_{\infty}, 
 \sup_{0<|m|\leq\eta} \frac{|f(x+m)-f(x)|}{|m|^{t}}\right)\leq M\right\},
\end{equation}
for $t \leq 1$.  We denote with $\mC(t,M,\eta)$ the set of functions that are locally H\"{o}lder  for each $x\in [0,1]$, i.e. $\mC(t,M,\eta)=\{f: f(x)\in\mC(t,x,M,\eta)\,\,\forall x\in [0,1] \}$ .
 Throughout this work, we will make the following assumption on $f_0$.
\begin{assumption}\label{ass:local_holder}
Assume  $f_0\in\mathcal C(t,M,\eta)$ where  $M(\cdot)$ and $\eta(\cdot)$ are bounded and uniformly bounded away from zero and where the smoothness function $t(\cdot)$ satisfies $0< t_1\leq \inf_{x\in[0,1]}t(x)$. 
\end{assumption}
It is well known that  regularity, both local and global, of a function is reflected in the speed at which its wavelet coefficients decay.
The following lemma formally characterizes the magnitude of  multiscale coefficients in terms of the {\em local} H\"{o}lder smoothness.

\begin{lemma} \label{lemma:coef2}
Denote with $\beta_{lk}^0=\langle f_0,\psi_{lk}\rangle$ the  multiscale coefficient of a function $f_0$ that satisfies Assumption \ref{ass:local_holder}. 
 Let $x \in [0,1]$ and for all $l >0$ define $k_l(x)\in\{0,1,\dots, 2^l-1\}$ such that $x\in I_{l \, k_l(x) }$. For $l >0$ and $k\in\{0,1,\dots,2^l-1\}$ let 
\begin{equation}\label{eq:def_t_eta}
\eta_{lk}=\min_{x\in I_{lk}}\eta(x)\quad\text{ and}\quad M_{lk} =\max_{x\in I_{lk}} M(x).
\end{equation}
When $l \geq \log_2 [1/(2\eta_{l\,k_l(x)})]$,   we have
\begin{equation}\label{eq:beta_ineq}
|\beta_{l\,k_l(x) }^0| \leq 2\, M_{l\,k_l(x)} 2^{-l[t(x)+1/2]}.
\end{equation}
\end{lemma}
\begin{proof}
For $k=k_l(x)$, we have 
$|\beta_{lk}^0|=2^{l/2} \left| \int_{k/2^l}^{(k+1/2)/2^l}[ f_0(y)  - f_0(y+2^{-l-1})]dy\right|$.
Then
 \begin{equation*}
 \begin{split}
 |\beta_{lk}^0| &\leq 2^{l/2} \int_{\frac{k}{2^l}}^{\frac{(k+1/2)}{2^l}}\left[ |f_0(y)-f_0(x)|+|f_0(x)  - f_0(y+2^{-l-1})|  \right]d y \leq 2\, M_{lk} 2^{-l[t(x)+1/2]}.
 \end{split}\qedhere
 \end{equation*}
\end{proof}

\begin{remark} \label{rem:coef2}
One can obtain an alternative bound $|\beta_{lk}^0| \leq  M_{lk} 2^{-l(t_{lk}+1/2)}$  in terms of $t_{lk}=\min_{x\in I_{lk}}t(x)$  which follows directly from the definition of $\mC(t,x,M,\eta)$.
\end{remark}

\begin{remark} \label{rem:smooth:extension}


While in this paper we study the case $t(\cdot)\leq 1$,  our results (Theorems \ref{thm:one}, \ref{thm:uq}, \ref{th:conc:spikeslab} and \ref{thm:np})  can be generalized  to higher-order H\"older functions for which Lemma \ref{lemma:coef2} applies with $S$-regular wavelet basis (see e.g. Extension 1 in \cite{castillo_rockova} for more discussion and references). 
\end{remark}
 
 \section{Spatial Adaptation in White Noise}\label{sec:wn}
 Donoho et al. \cite{asymptopia} characterized pointwise (as well as global)  properties of selective wavelet reconstructions showing their near-optimality  for estimating H\"{o}lderian functions at a given point $x\in[0,1]$. Here, we establish  {\em uniform}  (supremum-norm) local adaptation for {\em all} $x\in[0,1]$  focusing on \eqref{eq:sup_norm} under the white noise model \eqref{eq:model2} and various priors $\Pi(f)$.
 Adaptive supremum-norm concentration rate results (in white noise and regression) are still few and far between with pathbreaking progress made by multiple authors including \cite{hoffmann,yoo}. To date, results exist  {\em only} for homogeneous H\"{o}lderian functions under the spike-and-slab prior \citep{hoffmann,yoo} and, more recently, the Bayesian CART prior \citep{castillo_rockova}. Both of these priors leverage certain
 sparsity structure on the wavelet coefficients $\{\beta_{lk}\}$. We will show that {\em both} of these priors achieve uniform spatial adaptation.

\subsection{Bayesian CART}\label{sec:bcart1}
CART methods \citep{donoho,cart1,cart2} and other successful software developments including MARS \citep{mars} capture local aspects of the function being estimated by recursively subdividing the predictor space. Donoho et al. \cite{asymptopia} pointed out that {\em `the spatial adaptivity camp is, to date, a-theoretical and largely motivated by heuristic plausibility of the methods'}. While it has been more than 20 years since this seminal paper, there is a shortage of theoretical justifications focusing on spatial adaptation with  practically used machine learning methods. Here, we resurrect this question by focusing on Bayesian CART.

Bayesian CART corresponds to a wavelet prior  that prescribes a particular sparsity structure in the wavelet reconstruction according to a {\em binary tree $\mT$} (see \cite{castillo_rockova} for a more thorough exposition). A tree $\mT$  is defined as  a collection of hierarchically organized nodes $(l,k)$ where 
\begin{equation}\label{eq:tree_condition}
(l,k)\in\mT\Rightarrow (j,\lfloor k/2^{l-j})\in\mT\quad \text{for}\quad j=0,\dots, l-1.
\end{equation}
It will be useful to distinguish between two types of nodes: {\em internal} ones $\mT_{int}=\{(l,k)\in\mT: \{(l+1,2k),(l+1,2k+1)\}\in\mT\}\cup (-1,0)$  and {\em external} ones $\mT_{ext}=\mT\backslash\mT_{int}$ which are at the bottom of the tree. 
We then denote with  $\b_\mT=(\beta_{lk}:(l,k)\in\mT_{int})$ the `active' wavelet coefficients. 
Similarly as with the selective wavelet reconstruction ({\em RiskShrink} of \cite{donoho94}), Bayesian CART weeds out wavelet coefficients that are outside the tree, i.e. $\beta_{lk}=0$ when $(l,k)\notin\mT_{int}$.
Namely, for  $\Lmax\equiv\lfloor \log_2 n\rfloor$ we assume the tree-shaped wavelet shrinkage prior (Section 3 in \cite{castillo_rockova})
\begin{align}
\mT \qquad & \sim \qquad \Pi_{\mT} \label{eq:prior_beta1} \\
\{\beta_{lk}\}_{l\le \Lmax,k} \C \mT\ & \sim\  \pi(\b_\mT) \, \otimes 
\bigotimes_{(l,k)\notin\mT_{int}} \delta_0(\beta_{lk}),\label{eq:prior_beta2} 
\end{align}
where 
$$
\pi(\b_\mT)=\prod_{(l,k)\in\mT_{int}}\phi(\beta_{lk};0,1)
$$ 
is an independent product of standard Gaussians\footnote{\cite{castillo_rockova} also consider correlated wavelet coefficients.} and where $\Pi_\mT$ is the Bayesian CART prior \cite{cart1}.
This prior is essentially a heterogeneous Galton-Watson process with a node split probability 
\begin{equation}\label{eq:split_prob}
p_{lk}={P}[(l,k)\in\mT_{int}]=p_l=(1/\Gamma)^{l}
\end{equation}
for some $\Gamma>2$ (see  Section 2.1 in \cite{castillo_rockova} and  \cite{rockova_saha}). 
\subsubsection{Uniformly Adaptive Rate}
The following theorem establishes  uniform spatial adaptation of Bayesian CART in the supremum-norm sense. In other words,  the posterior is shown to contract at a locally minimax rate, up to a log factor, uniformly for all $x\in[0,1]$. While very intuitive,  such a result has not yet been formalized in the Bayesian literature.

\begin{theorem}\label{thm:one}
Assume the Bayesian CART  prior  \eqref{eq:prior_beta1} and \eqref{eq:prior_beta2}
 with a split probability \eqref{eq:split_prob} for some sufficiently large $\Gamma>0$.
Under the model \eqref{eq:model2} and with $t,M$ and $\eta$ satisfying Assumption \ref{ass:local_holder},  we have
\begin{equation}\label{eq:supnorm_rate}
\sup\limits_{f_0\in\mathcal C(t,M,\eta)}E_{f_0}\Pi\left[f: \sup_{x\in[0,1]}\zeta_n(x)|f(x)-f_0(x)|>M_n\,\big|\,  Y\right]\rightarrow 0
\end{equation}
for $\zeta_n(x)= \left(\frac{n}{\log n}\right)^{\frac{t(x)}{2t(x)+1}}$ and for any $M_n\rightarrow\infty$ that is faster than $\sqrt{\log n}$.
\end{theorem}
The proof is provided in Section \ref{sec:proof_thm1}.

\smallskip

 The first step in the proof of Theorem \ref{thm:one} is showing that trees, a-posteriori, grow deeper in domains where $f_0$ is less smooth.
 This property is summarized in Lemma \ref{lemma:locally_small} in the Supplemental Materials.
Supremum-norm convergence rate results are valuable for constructing  confidence bands.
For example,  
Theorem \ref{thm:one} implies the non-parametric Bernstein-von Mises phenomenon in the multiscale space   which can be   used to construct credible bands with {\em exact} asymptotic coverage (see, e.g., Theorem 4.1 in \cite{castillo_rockova}). This set, however, is not guaranteed to have the optimal size (i.e. its diameter shrinking at the minimax rate).
Here, we will focus on  constructing  valid {\em adaptive} confidence bands. With spatially varying functions (such as the local H\"{o}lder functions from Section \ref{sec:spatial_functions}), one would expect the width of the confidence band to vary with the smoothness $t(\cdot)$ and be wider where $t$ is smaller. Keeping the diameter constant throughout may yield bands that are more conservative in certain areas of the sample space.

\subsubsection{Locally Adaptive Bands}\label{sec:uq}
A reasonable requirement for band construction is that their diameter shrinks at the minimax rate of estimation, up to possibly a slow multiplication factor. 
When the degree of smoothness is known, multiscale\footnote{They resemble the $L^\infty$ balls \citep{ray}.} credible balls  can be constructed (see  (5) in \citep{castillo_nickl1}) and  intersected  with qualitative restrictions on $f_0$  to obtain `optimal'  frequentist confidence sets (which shrink at the optimal rate).

We construct optimal confidence sets when the smoothness $t(\cdot)$ is {\em unknown} and {\em varying} over $[0,1]$.  Confidence bands that are simultaneously adaptive and honest, of course, do not exist in full generality  \citep{low}. Gine and Nickl \cite{gine_nickl} point out, however, that such confidence sets exist for certain generic subsets of H\"{o}lderian functions, the so-called self-similar functions  \citep{picard,bull, gine_nickl2,nickl_szabo,ray}, whose complement was shown to be negligible \citep{bull}.
Under self-similarity, \cite{ray,castillo_rockova} constructed  adaptive credible bands for  homogeneous H\"{o}lderian functions under the spike-and-slab prior and the Bayesian CART, respectively.

Here, we extend the notion of  self-similarity  to {\em inhomogeneous} 
H\"{o}lder classes for which it is possible to construct a {\em locally adaptive} confidence set $\mathcal C_n$ in the sense  that
\begin{equation}\label{eq:diameter}
\sup_{f,g\in\mathcal C_n}\left[\sup_{x\in[0,1]} \frac{\zeta_n(x)}{v_n} |f(x)-g(x)|\right]=\mathcal{O}_{P_{f_0}}(1)
\end{equation}
for some suitable sequence $v_n\rightarrow\infty$ and where $\zeta_n(x)=(n/\log n)^{t(x)/(2t(x)+1)}$. Note that the diameter of $\mathcal C_n$ depends on $x$ and equals the minimax rate of estimation (inflated by $v_n$) {\em at every point} $x\in[0,1]$. Below, we  formally introduce the notion of locally self-similar functions.
\begin{definition}(Local Self-Similarity)   We say that $f\in\mathcal C(t,M,\eta)$ is {\em locally self-similar} at $x\in[0,1]$ if, for some   $c_1>0$ and an integer $j_0$, 
\begin{equation} \label{sscond}
| K_j(f)(x) - f(x) | \ge  2^{-j \, t(x)} c_1 \quad  \text{for all } j \ge j_0,
\end{equation}
where $K_j(f)=\sum_{l\le j-1}\sum_k \, \langle\psi_{lk},f\rangle \, \psi_{lk}$ is the wavelet projection  at level $j$.
The class of all   self-similar functions at $x$ will be denoted by $\mathcal{C}_{SS}(t(x),x,M(x),\eta(x))$.
Moreover, we denote with $\mathcal{C}_{SS}(t,M,\eta)$ a set of functions  that are self-similar for all $x\in[0,1]$.
\end{definition}

\begin{figure}[!t]
\scalebox{0.35}{\includegraphics{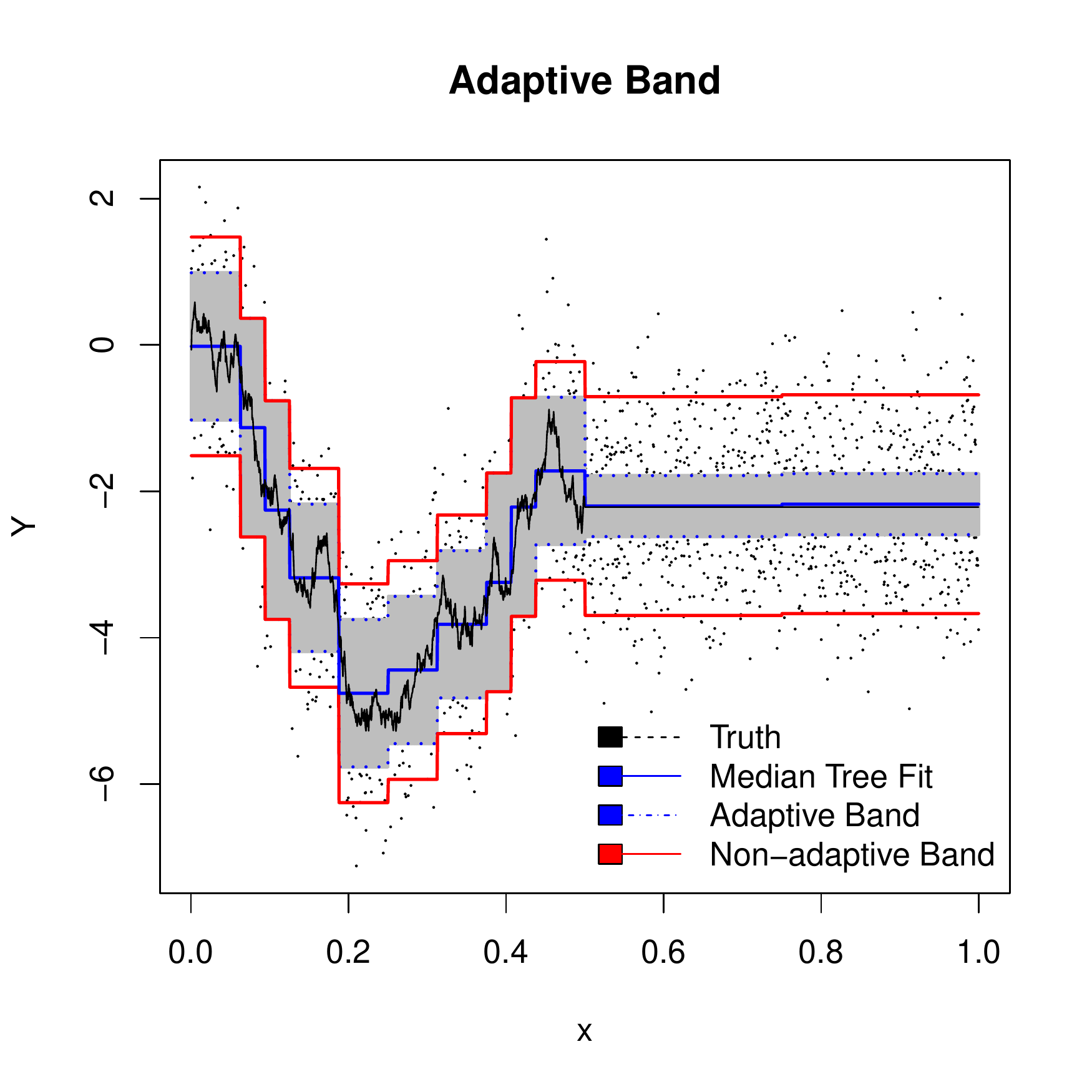}}
\scalebox{0.35}{\includegraphics{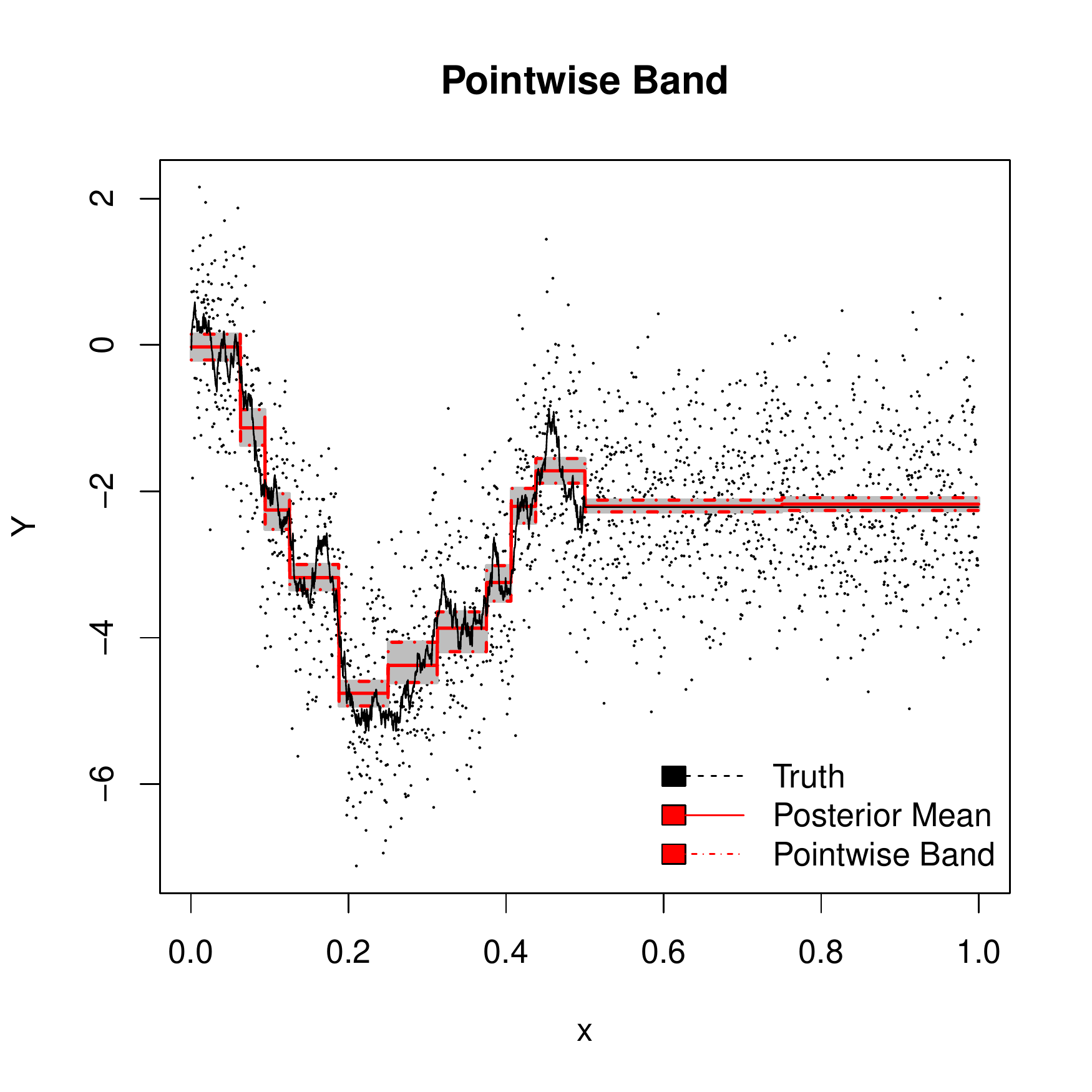}}
\caption{The Brownian motion example (Example \ref{ex:bm}) with dyadic Bayesian CART. (left) Adaptive confidence band \eqref{credible} (blue) and non-adaptive band (red). (right) Point-wise $95\%$ credible intervals.}
\label{bm}
\end{figure}

For spatially heterogeneous H\"{o}lderian functions, we construct locally adaptive confidence bands whose  width is varying and  reflects smoothness at each given $x$
While related to previous constructions (see e.g. \cite{castillo_rockova} for the homogeneous case), its simplicity and ease of computability 
make our band    particularly appealing in practice (see Figure \ref{bm}). In addition, we are not aware of any other related frequentist band for the case of heterogeneous smoothness.
 We center our confidence bands around  a pivot estimator, the median tree estimator  \cite{castillo_rockova}.

\begin{definition}(The Median Tree)
Given a posterior distribution $\Pi_\mT[\cdot\given Y]$ over tree-shaped coefficient subsets, we define the {\em median tree}  $\cT^*_Y$ as the following set of nodes  
\begin{equation} \label{bulktree}
\cT^*_Y=\left\{ (l,k),\ l\le L_{max},\ \ \Pi[(l,k)\in \cT_{int} \given Y] \ge 1/2 \right\}.
\end{equation}
\end{definition}
We  define the resulting median tree estimator as 
\begin{equation} \label{bulkest}
\wh f_T (x) = \sum_{(l,k)\in \cT^*_Y} Y_{lk} \psi_{lk}(x)
\end{equation}
which is shown to attain the near-minimax rate of estimation at each point (see the proof of Theorem \ref{thm:uq}).
Next, we define the {\em local radius} (which varies with $x$) as
\begin{equation}  \label{radiusprox}
\sigma_n(x) = v_n \sqrt{\frac{\log{n}}{n}} \sum_{l=0}^{L_{max}}
\mathbb{I}\{(l,k_{l}(x))\in \cT_Y^*\} |\psi_{lk_l(x)}(x)|
\end{equation}
for some $v_n\rightarrow\infty$ to be chosen.
Finally, we construct the confidence band according to the following prescription
\begin{equation} \label{credible}
\cc_n=\left\{f :\sup_{x\in[0,1]}\left[\frac{1}{\sigma_n(x)} |f(x) -\wh f_T(x)|\right] \le  1 \right\},
\end{equation}
where $\wh f_T$ is as in \eqref{bulkest} and $\sigma_n(x)$ is as in \eqref{radiusprox}.

\begin{theorem} \label{thm:uq}
Let $\Pi$ be the prior as in the statement of Theorem \ref{thm:one}. 
Then for $\cc_n$ defined in \eqref{credible}, uniformly over $t, M$ and $\eta$  that satisfy the Assumption \ref{ass:local_holder}, as $n\to\infty$,
\[  \inf_{f_0\in \mathcal{C}_{SS}(t,M,\eta)} P_{f_0}(f_0\in \cc_n)  \to 1.\]
Uniformly over $f_0\in \mathcal{C}_{SS}(t,M,\eta)$, the  diameter verifies \eqref{eq:diameter}, as $n\to \infty$, and  the credibility of the band  satisfies
\begin{align}
 \Pi[\cc_n\given X] & = 1 + o_{P_{f_0}}(1). 
\end{align} 
\end{theorem}

The proof is provided in Section \ref{proof:thm_uq}.

\smallskip
According to Theorem \ref{thm:uq}, the band \eqref{credible} has asymptotic coverage $1$. It is possible to intersect \eqref{credible} with a multi-scale ball (similarly to \cite{ray,castillo_rockova}) to obtain asymptotic coverage $1-\gamma$ for some small $\gamma>0$ as a consequence of the non-parametric Bernstein-von Mises (BvM) theorem.   In order to illustrate the practical usefulness of Theorem \ref{thm:uq}, we revisit the Brownian motion example (Example \ref{ex:bm}) from Figure \ref{bm}. We implement a dyadic version of the Bayesian CART algorithm \citep{cart1} which splits only at dyadic rationals (as opposed to Figure \ref{fig:bm} (on the right) where Bayesian CART splits at observed values). We plot the adaptive band \eqref{credible} choosing $v_n=2$ together with a non-adaptive band obtained by taking the maximal diameter $\sigma(x)$ over the domain $[0,1]$ (as in Theorem 4 of \cite{castillo_rockova}). We can see the benefits of our adaptive construction, where the width is larger in the first half of the domain where the function meanders according to the Brownian motion (expected since the smoothness is smaller than $1/2$). Compared with Figure \ref{bm} on the right, the point-wise $95\%$ credible band does not yield satisfactory coverage on this example.

\subsection{Spike-and-Slab Priors} 
Spike-and-slab priors are arguably one of the most ubiquitous priors in statistics (see references \cite{chipman_wavelet,hoffmann,ray} for wavelet shrinkage contexts).
Compared with the Bayesian CART prior from Section \ref{sec:bcart1},
spike-and-slab priors allocate positive prior mass to {\em any} subset  $\mT$ of $\{ (l,k): l \leq L_{max}\}$, not just tree-shaped subsets.  We  define the spike-and-slab prior  through the following hierarchical model.
\begin{assumption}(Spike-and-Slab Prior)\label{ass:prior}
\begin{itemize}  
\item Prior on $\mT\subseteq\{(l,k): 0\leq k < 2^l, l\leq L_{max} \}$: There exist constants $c_T, C_T>0$ such that 
\begin{equation}\label{cond:priorMT:SS}
c_T\, \omega_{l} \leq \frac{ \Pi( \mT \cup \{ (l,k)\} )}{ \Pi( \mT  )} \leq C_T\, \omega_{l} \quad\forall \mT\,\,\text{such that}\,\, (l,k)\notin\mT
\end{equation}
 for some   positive sequence $\omega_l$ such that, for some  $B_\omega>0$ and $\delta>0$,
\begin{equation}\label{cond:omegal}
n^{-B_\omega} \leq \omega_{l} \leq  n^{(1-\delta)/2} 2^{-l} \quad \text{for}\quad l \leq L_{max}.
\end{equation} 
\item There exist probability densities $\pi_{lk}(\cdot)$ on $\mathbb R$ such that, conditionally on $\mT$,
$$ 
\beta_{lk} \stackrel{ind}{\sim} \pi_{lk}  \quad \forall (l,k)\in \mT \quad \text{and}\quad \beta_{lk}=0, \quad \forall (l,k)\notin \mT
$$
and there exist $R,c_R, C_R>0$ such that 
\begin{equation}\label{pibeta}
 c_R \leq \inf_{|\beta|\leq  R } \pi_{lk} (\beta) \leq \sup_{\beta \in \R} \pi_{lk} (\beta) \leq C_R.
\end{equation}
\end{itemize}
\end{assumption}

While seemingly similar to the prior considered in \cite{hoffmann}, our  Assumption \ref{ass:prior} is much weaker. Indeed,  our prior construction  subsumes the spike-and-slab prior of \cite{hoffmann} by imposing weaker constraints on the decay of inclusion probabilities. Note that  $\omega_l$'s in \eqref{cond:omegal} are allowed to be much larger than in \cite{hoffmann} which assume $n^{-B}\leq \omega_l\leq 2^{-j(1+\tau)}$ for some $\tau>1/2$. This perhaps subtle difference is of great practical importance and indicates that optimal sup-norm adaptation occurs in far less sparse situations than originally perceived. Another important difference is that  we {\em do not require} the binary indicators $\1(\beta_{lk}\neq 0)$ for $0\leq k <2^l$ and $l\leq L_{max}$ to be  {\em iid} Bernoulli random variables. This extension allows us to consider, for example, Ising prior constructions \cite{besag} which allow the inclusion indicators to be related through a Markovian model.

\begin{theorem}\label{th:conc:spikeslab}
Consider the model \eqref{eq:model2}  with a prior $\Pi$ on $\{\beta_{lk}\}_{lk}$ following the Assumption \ref{ass:prior}. Let $t,M$ and $\eta$ satisfy the Assumption \ref{ass:local_holder}. Then
$$ \sup_{f_0\in \mathcal C(t,M,\eta)} 
E_{f_0}\Pi\left[f: \sup_{x\in[0,1]}\zeta_n(x)|f(x)-f_0(x)|>\wt M\,\big|\,  Y\right]\rightarrow 0,
$$
 for all sufficiently large $\wt M>0$ where  $\zeta_n(x)=(n/\log n)^{t(x)/(2t(x)+1)}$.
\end{theorem}

The proof is provided in  Section \ref{sec:proof:th:conc:spikeslab}.
\smallskip

Theorem \ref{th:conc:spikeslab} shows that, unlike Bayesian CART, the spike-and-slab priors achieve the {\em exact} rate uniformly over the entire domain $[0,1]$
without any additional logarithmic penalty (\cite{castillo_rockova} showed that the log-factor in Bayesian CART is non-negotiable). In Lemma \ref{lem:smallcoef} (an analog of Lemma 1 in \cite{hoffmann}) we show that the posterior concentrates on a subset of large enough coefficients. This fact can be used to show that  the median probability model (MPM) \citep{barbieri,barbieri2, ray} consisting of all coefficients with at least $50\%$-posterior probability of being active is an (exact) rate-optimal estimator. Following the strategy of Proposition 4.5 in \cite{ray} one can then show that Theorem \ref{thm:uq} remains true for the spike-and-slab prior when replacing the median tree estimator with MPM and with $v_n$ that can grow slower at a rate at least $\sqrt{\log n}$ \citep{ray}.

\section{Spatial Adaptation in Non-parametric Regression}\label{sec:npreg}
Throughout this section, we  assume the canonical non-parametric regression setup \eqref{model} with $\sigma^2=1$.
While  nonparametric regression  with a {\em regular} design and the white noise model are 
asymptotically equivalent (e.g. under the usual smoothness assumption $t(x)>1/2$ \citep{brown_low}),
 optimality of a procedure in one setup {\em does not} automatically imply optimality in the other. In Section \ref{sec:np1}, we show rate-optimality of Bayesian CART in  non-parametric regression when $t(x)>1/2$ {\em without} assuming regular designs. Later in Section \ref{sec:repulsive}, we relax the restriction $t(x)>1/2$ and propose  new `repulsive' partitioning priors (related to adaptive-knot splines) and show that they are exact-rate adaptive. First, we describe some not so optimistic findings for hierarchical Gaussian priors.

\subsection{Gaussian Processes} \label{sec:LBGaussian}
 In Figure \ref{doppler1} and \ref{fig:bm} we have seen  that spatial adaptation is {\em not attainable} by methods which are not sufficiently localized.
In this section, we formally show that several practically used Gaussian process priors
  {\em do not} lead to spatially adaptive concentration rates.  While this phenomenon is not entirely surprising, it is nevertheless worthwhile to document it formally. In particular, we provide lower bound results showing sub-optimality of Gaussian processes in terms of a {\em global} estimation loss.
To this end, we consider the following heterogeneous-smoothness assumption which aligns with Figure \ref{fig:bm} where the function $f_0$ has smoothness $\alpha_1<1/2$ on $[0,1/2]$ and $\alpha_2=1$ on $(1/2, 1]$. 
\begin{assumption} 
Assume that the Haar wavelet decomposition of a function $f_0$  satisfies  (with $I_l=\{0,1,\dots,2^{l}-1\}$)
 \begin{equation}\label{eq:f0}
 f_0(x) =\psi_{-10}(x)\beta_{-10}^0+ \sum_{l=0}^\infty \sum_{k \in I_l} \beta^0_{lk} \psi_{l k}(x)
 \end{equation} 
 where (for  $I_{l1}\cup  I_{l2} = I_l $ with $|I_{l1}| \asymp |I_{l2}|$ and some $M_1>0$)
 \begin{equation}\label{spatialBeta}
 \max_{k \in I_{l1}} |\beta^0_{lk}| \leq M_1 2^{-l(\alpha_1 +1/2)}\quad\text{and}\quad \max_{k \in I_{l2}} |\beta^0_{lk}| \leq M_1 2^{-l(\alpha_2 +1/2)} \quad\text{with}\quad \alpha_1 < \alpha_2.
 \end{equation}
 \end{assumption}
 Methods that are globally, but not locally,  adaptive are expected to adapt to the worse-case scenario and attain the slower rate  determined by the smaller smoothness  $\alpha_1$.
 We will formalize this intuition  by assessing the quality of the reconstruction with an $L^2$ loss over the {\em entire} domain as well as the {\em smoother} domain determined by $\alpha_2$, i.e. we define
 $$
 \|f-f_0\|_2^2\equiv\int_0^1|f(x)-f_0(x)|^2d x\quad\text{and}\quad  \|f-f_0\|_{1/2,1}^2\equiv\int_{1/2}^1|f(x)-f_0(x)|^2d x.
 $$
 

We now consider three hierarchical Gaussian processes on 
\begin{equation}\label{eq:expansion}
f(x) =\psi_{-10}(x)\beta_{-10}+ \sum_{l=0}^\infty \sum_{k \in I_l} \beta_{lk} \psi_{l k}(x)
\end{equation}
induced through a prior on the sequence
$\{\beta_{lk}\}_{lk}$. These priors were studied in \cite{rousseau:szabo:17}. This section assumes a regular design $x_i=i/n$ for  $n=2^{\Lmax+1}$  with some $\Lmax>0$.
 \begin{itemize}
 \item (T1) ({\em Sieve Prior}) Given a truncation level $1\leq L\leq L_{max}$ and some  fixed $\tau>0$ we have
  $$ 
   \beta_{lk}  \stackrel{iid}{\sim} \mathcal N(0,\tau^2)\quad \text{for}\quad l \leq L\quad\text{and}\quad \beta_{lk}=0\quad\text{for}\quad l>L.
  $$
  The truncation level $L$ is assigned a prior $\pi(L)$ which satisfies Lemma 3.4 in \cite{rousseau:szabo:17} (e.g. a hypergeometric or a Poisson prior distribution).
  \smallskip
  
 \item (T2) ({\em Scale Parameter}) Given $\alpha>0$ assume
  $$ 
  \beta_{lk}= \begin{cases}
  2^{-l(\alpha+1/2)} Z_{lk} \quad \text{where}\quad Z_{lk}  \stackrel{iid}{\sim} \mathcal N(0,\tau^2) \quad&\text{for}\quad l \leq  L_{max}\\
  0\quad&\text{for}\quad l>L_{max}
  \end{cases}
 $$
where the scale parameter arrives from $ \tau \sim \pi_\tau$, with $\pi_\tau$ satisfying the assumptions of Lemma 3.5 of \cite{rousseau:szabo:17} (e.g. an inverse gamma or a Weibull distribution). 

\smallskip

  \item (T3)  ({\em Rate Parameter}) Given  $\tau>0$  assume
$$
 \beta_{lk}=\begin{cases}
  2^{-l(\alpha+1/2)} Z_{lk} \quad \text{where}\quad Z_{lk}  \stackrel{iid}{\sim} \mathcal N(0,\tau^2) \quad&\text{for}\quad l \leq L_{max} \\
  0\quad&\text{for}\quad l>L_{max}
  \end{cases}
$$
where the `smoothness' parameter $\alpha$ arrives from
 $
 \alpha \sim \pi_\alpha
 $ satisfying the assumptions of Lemma 3.6 of \cite{rousseau:szabo:17} (e.g. a gamma distribution). 

 \end{itemize} 
These priors have been studied in a multitude of works, see e.g. \cite{arbel,shen_ghosal} for the setup (T1) \cite{belitzer_ghosal,knapik,szabo_vdv_vz} for the Gaussian process priors (T2) and (T3).
More recently,  this framework has been studied in \cite{rousseau:szabo:17} in the case of Fourier-series priors where both  lower bounds and an upper bound have been obtained in the case of non-linear regression.  We adapt their proof to the wavelet basis case with functions satisfying \eqref{spatialBeta}. 
 In this case, we note that for any $L $ such that $2^L \leq n$ and for $f_{\beta,L}$ denoting the  Haar wavelet expansions \eqref{eq:expansion} truncated at $L\geq 1$ we have
   for $x_i=i/n$ 
   $$ 
\|f_{\beta,L} -f_{\beta_0,L} \|_n^2=\frac{1}{n}   \sum_i [f_{\beta,L} (x_i) -f_{\beta_0,L} (x_i) ]^2  = \| \b_L - \b_{0,L} \|_2^2=\|f_{\beta,L} -f_{\beta_0,L}\|_2^2.
   $$
   where $\b_L=(\beta_{lk}:l\leq L)'$ (resp. $\b_{0,L} $ ) is  the truncated version of $\be$ (resp $\b_{0}$).
  
   \begin{theorem}\label{th:GaussianP}
   Let $f_0$ satisfy \eqref{spatialBeta} and consider either of the priors (T1)-(T3). Then for $Y=(Y_1,\dots, Y_n)'$ arising from \eqref{model} with $x_i=i/n$ we have
   \begin{equation}\label{eq:statement1}
   \Pi\left[ \|f-f_0\|_2 \leq M_n \epsilon_n(\alpha_1) \C Y\right] = 1 + o_{P_{f_0}}(1) \quad\text{as $n\rightarrow\infty$}
   \end{equation}
   with 
   $$\epsilon_n(\alpha_1) =
   \begin{cases} 
   &(n/\log n)^{-\alpha_1/(2\alpha_1+1)} \quad\text{under (T1), (T2) if  $\alpha_1 < 2\alpha+1$, and (T3)} \\
   & (n/\log n)^{-(2\alpha+1)/(4\alpha+4)} \quad\text{under (T2) if $\alpha_1 \geq 2\alpha+1$}.
   \end{cases}
   $$
Moreover, for all functions satisfying \eqref{spatialBeta}  and, for some $c>0$,
   \begin{equation}\label{spatialBeta2}
 \min_{k \in I_{l1}} |\beta^0_{lk}| \geq c\, 2^{-l(\alpha_1 +1/2)}
 \end{equation}
 we have for all $\delta >0$ 
  \begin{equation}\label{eq:statement2}
  \Pi\left[ \|f-f_0\|_{1/2,1} \leq n^{-\delta} \epsilon_n(\alpha_1) \C Y\right] =  o_{P_{f_0}}(1)\quad\text{as $n\rightarrow\infty$}. 
  \end{equation}
   \end{theorem}
 The proof is provided in Section \ref{sec:proof_lb}.

\smallskip

The first statement \eqref{eq:statement1} shows that the posterior under the Gaussian sequence priors adapts to the {\em worse} smoothness $\alpha_1$. 
Moreover,  the second statement \eqref{eq:statement2} implies that, under a suitable identifiability condition, the posterior is {\em incapable} of achieving a faster rate on the smoother domain (determined by $\alpha_2>\alpha_1$), rendering adaptation to $\alpha_2$ impossible. Note that similarly to \cite{rousseau:szabo:17}, 
the same conclusions holds if one deploys an empirical Bayesian procedure based on the marginal maximum likelihood estimator on $L$ for T1 (resp. $\tau$ for T2 and $\alpha$ for T3). 

\subsection{Bayesian CART}\label{sec:np1}
This section reports positive findings in the context Bayesian CART. In particular, we show a regression analogue of Theorem \ref{thm:one} assuming  $t(x)>1/2$.
\cite{castillo_rockova} also study Bayesian CART in regression but  with a {\em regular} design where the prior is assigned to {\em empirical} wavelet coefficients. This re-parametrization closely resembles the white noise model, enabling a more direct transfer of the results.
Here, we follow an alternative route. A  perhaps more transparent approach is to assign a prior {\em directly} to the actual (not empirical) wavelet coefficients (similarly as in \cite{yoo}). This strategy aligns more closely with what is done in practice. We pursue this direction here and, in addition, consider designs that are not necessarily regular.

With a vector of observations $Y=(Y_1,\dots,Y_n)'$ and $F_0=(f_0(x_1),\dots, f_0(x_n))'$, we can re-write \eqref{model} in a matrix notation as
\begin{equation}\label{eq:regression}
Y=X\b_0+ \bm\nu,\quad\text{where}\quad \bm\nu=F_0-X\b_0+\bm \epsilon\quad\text{with}\quad  \bm \epsilon\sim\mathcal N(0,I_n),
\end{equation}
where $\b_0=(\beta^0_1,\dots,\beta^0_{p})'\in\R^{p}$  is  a sparse vector  of  multiscale coefficients $\langle f_0,\psi_{lk}\rangle$ 
ordered according to $2^l+k$ and where $X=[X_1,\dots, X_p]$ is  an $(n\times p)$ Haar wavelet design matrix constructed as follows.
 The matrix $X=(x_{ij})$  consists of $p$ wavelet bases $\psi_{lk}$ up to a resolution $\Lmax$ which satisfies $p=2^{\Lmax}=\lfloor C^* \sqrt{n/\log n}\rfloor$ for some suitable $C^*>0$. The columns in $X$ are  implicitly ordered according to the index $2^l+k$ (from the smallest to the largest), where
$$
x_{ij}=\psi_{lk}(x_i)\quad\text{with}\quad j=2^l+k.
$$ 
Because we assume $t(x)>1/2$, we do not need resolutions larger than $\Lmax$ to be able to approximate $f_0$ well. We will be denoting with $\bT$ all tree-shaped subsets of nodes $(l,k)$ such that $l\leq \Lmax$.
For a tree $\mT\in\bT$ and a vector $\b\in\R^{p}$, we denote with  $\b_{\mT}$ the subset of coefficients inside the tree and with $\b_{\bmT}$ the complement. 
Similarly, we split the design matrix $X$ into active covariates $X_{\mT}$ (that correspond to $(l,k)\in\mT_{int}$) and the complementary inactive ones $X_{\bmT}$. 
It will be advantageous\footnote{We can take advantage of certain properties of projection matrices. Other priors can be considered as well.} to use the unit-information $g$-prior for $\b_\mT$ 
\begin{equation}\label{eq:g_prior}
\b_\mT\sim\mathcal N(0,g_n(X_\mT'X_\mT)^{-1})\quad\text{with $g_n=n$}
\end{equation}
which yields the following marginal likelihood under each tree $\mT$ 
$$
N_\mT(Y)=\frac{1}{(2\pi)^{n/2}(1+g_n)^{|\mT|/2}}\exp\left\{-\frac{1}{2}Y'[I-X_\mT\Sigma_\mT X_\mT']Y\right\},
$$
where 
$$
\Sigma_\mT= c_n(X_\mT'X_\mT)^{-1}\quad\text{with $c_n=g_n/(g_n+1)$}.
$$

Throughout this section, we will denote with  $n_{lk}^L$ (resp. $n_{lk}^R$)  the number of observations that fall inside the domain of the left (resp. right) wavelet  piece $\psi_{lk}$, i.e. $n_{lk}=\sum_{i=1}^n\mathbb I(x_i\in I_{lk})=n_{lk}^R+n_{lk}^L$ and we define
$$
\bar n_{lk}=\max\{n_{lk}^R,n_{lk}^L\}\quad\text{and}\quad \munderbar n_{lk}=\min\{n_{lk}^R,n_{lk}^L\}.
$$ 
The regular design $x_i=i/n$ satisfies $\bar n_{lk}=\munderbar n_{lk}=n/2^{l+1}$. 
Here, we assume a fixed design $\mX=\{x_i\in[0,1]:1\leq i\leq n\}$ that is {\em not necessarily} regular. Instead, we make the following   design balance assumption.

\begin{assumption}(Balanced Design)\label{ass:design}\,
Let $\wtLmax$ be such that $2^{\wtLmax}=\mathcal O(n/\log n)$. We
say that the design $\mX$ is $\upsilon$-regular for some $\upsilon>0$  if for any $(l,k)$ s.t. $0\leq l\leq \wtLmax$
\begin{equation}\label{ass:balance1}
\frac{c\,n}{2^l}\leq \munderbar n_{lk}\leq  \bar n_{lk}\leq \frac{(C+l)n}{2^l}\quad\text{for some $c,C>0$}
\end{equation}
and, for some  $C_d>0$,
\begin{equation}\label{ass:balance2}
0\leq \bar n_{lk}-\munderbar n_{lk}\leq  C_d \frac{\sqrt{n}\log^\upsilon n}{2^{l/2}}.
\end{equation}
\end{assumption}
 Note that the threshold $\Lmax$ used to construct the design matrix $X$ is smaller than the threshold $\wtLmax$ in  Assumption \ref{ass:design} and all  partitioning cells induced by $\mT\in\bT$ are guaranteed to have  at least one observation.  The  Assumption \ref{ass:design} is not overly restrictive. Indeed,  we show in Lemma \ref{lemma:design} (Supplemental Materials) that  \eqref{ass:balance2}  is satisfied with   probability approaching one when $\mX$ arises from a uniform distribution on $[0,1]$ for $\upsilon\geq 1/2$.
Irregular observations ultimately induce {\em correlated} Haar wavelet designs $X$ where the correlation pattern has a particular hierarchical structure described in the following lemma.

\begin{lemma}\label{lemma:cor}
Let $X_i$ and $X_j$ be columns of $X$ that correspond to nodes $(l_2,k_2)$ and $(l_1,k_1)$, respectively. Then we have
\begin{align*}
|X_j'X_i|=0&\,\,\text{when}\,\, (l_2,k_2)\,\,\text{is {\em not} a descendant of $(l_1,k_1)$},\\
|X_j'X_i| \leq 2^{\frac{l_1+l_2}{2}}|n_{l_2k_2}^L-n_{l_2k_2}^R|&\,\,\text{when}\,\, (l_2,k_2)\,\,\text{is a descendant of $(l_1,k_1)$}.
\end{align*}
\end{lemma}
\proof
When $(l_2,k_2)$ is {\em not} a descendant of $(l_1,k_1)$, the domains of $\psi_{l_1k_1}$ and $\psi_{l_2k_2}$ do not overlap, yielding orthogonality.
When  $(l_2,k_2)$ is a descendant of $(l_1,k_1)$, the wavelet domains satisfy $I_{l_2k_2}\subset I_{l_1k_1}$ and 
$X_j'X_i$ will be (up to a sign) equal  to the size of the amplitude product $2^{(l_1+l_2)/2}$ multiplied by the excess number of observations falling inside the longer wavelet piece $\psi_{l_2,k_2}$.\qed

\smallskip
As with other related consistency results in regression (see e.g. \cite{Narisetty}), we cannot allow for too much correlation in the design $X$.
Fortunately, {\em balanced} designs that satisfy Assumption \ref{ass:design} are not too collinear and yield diagonally dominant covariance matrices 
with well-behaved  eigenvalues (see Lemma \ref{lemma:eigen} below).

\begin{lemma}(Eigenvalue Bounds)\label{lemma:eigen}
Under the Assumption \ref{ass:design} with $0\leq \upsilon\leq 1/2$ and some suitable $c>0$, the eigen-spectrum of $X_\mT'X_\mT$ for each $\mT\in \bT$ satisfies  
\begin{equation}\label{eq:min_eigen1}
\munderbar{\lambda}\, n\leq \lambda_{min}(X_\mT'X_\mT)\leq \lambda_{max}(X_\mT'X_\mT)\leq \bar\lambda\, n\,\log n\quad\text{for some $0<\munderbar\lambda\leq \bar\lambda$}.
\end{equation}
\end{lemma}
The proof is provided in Section \ref{sec:proof_lemma_eigen}.

\medskip

We are now ready to state a non-parametric regression version of Theorem \ref{thm:one}.

\begin{theorem}\label{thm:np}
Assume the regression model \eqref{eq:regression} under  Assumption \ref{ass:design} for some  $0\leq \upsilon\leq 1/2$ and   $c,C,C_d>0$. 
Assume the Bayesian CART tree prior $\Pi_\mT$ with a split probability $p_l=(\Gamma)^{-l^{2}}$  for some sufficiently large $\Gamma>0$
and the $g$-prior \eqref{eq:g_prior}. Assume that $t,M$ and $\eta$ satisfy Assumption \ref{ass:local_holder} with $t_1>1/2$,  then
$$
\sup_{f_0\in\mathcal C(t,M,\eta)}E_{f_0}\Pi\left[f: \sup_{x\in[0,1]}\zeta_n(x)|f(x)-f_0(x)|>M_n\,\big|\,  Y\right]\rightarrow 0
$$
for $\zeta_n(x)= \left(\frac{n}{\log n}\right)^{\frac{t(x)}{2t(x)+1}}$ and for any $M_n\rightarrow\infty$ that is faster than ${\log^{1/2} n}$.
\end{theorem}
The proof is provided in  Section \ref{pr:thm:np}.
\smallskip

Note that the prior split probability decays more rapidly, i.e. $p_l=(\Gamma)^{-l^2}$, to accommodate the irregular design assumption.
An analogous statement can be obtained for the spike-and-slab prior using a similar approach as in the proof of Theorem \ref{th:conc:spikeslab}. 
In addition,  the confidence set construction in \eqref{credible} remains valid also under the non-parametric regression setting.
Indeed,  rate-optimality of the posterior  implies   rate-optimality of the median-tree estimator and the regression variant of Theorem  \ref{thm:uq}  thus holds under the assumption $t_1> 1/2$.

\subsection{Repulsive Partition Prior} \label{sec:repulsive}
In the previous section, we studied a prior on wavelet coefficients that corresponds to recursive partitioning.  In this section, we  propose a different  partitioning prior on piecewise constant functions which  relates to variable knot spline techniques \citep{mammen_geer,tibshirani_trend}. 
We assume that  a partition  $S = (I_1^S, \cdots, I_J^S) $ of $[0,1]$ is not necessarily induced by a tree but arrives from a `determinantal-type' prior  
\begin{equation}\label{prior:repulsive_S}
S \sim \pi_S \propto \prod_{j=1}^J |I_j^S|^B\mathbb I(S\in \mathbb S)\quad\text{for some $B>0$}, 
\end{equation}
where $\mathbb S$ contains all partitions made out of blocks with endpoints belonging to a fixed grid $\mathcal I_n = (z_\ell: \ell \leq N_n)$ such that $z_0=0$ and $z_{N_n}=1$ and
\begin{equation}\label{eq:z}
0\leq  z_\ell< z_{\ell+1}, \,\,\frac{ C_1 \log n}{n }\leq z_{\ell+1}-z_\ell \leq \frac{C_2 \log n}{n },\,\, \ell\leq N_n \quad \text{for some}\quad C_1, C_2>0.
\end{equation}

Note that the size of an  interval in $S$ can be measured either in terms of its length or in terms of its number of units, i.e. number of elements in the grid $\mathcal I_n$ belonging to it. 
 
We refer to \eqref{prior:repulsive_S} as a repulsive partitioning prior because it prevents  the splits from occurring too close to one another. 
The prior \eqref{prior:repulsive_S} thus rewards   partitions that are more regular.
The set $\mathcal I_n$ contains candidate knots for possible split, e.g.   a subset of  observed design points. 
While in variable knot spline techniques (such as MARS \citep{mars}) knot points are added, removed and allocated recursively using cross-vaildiation, here we let the posterior distribution choose the knots in a data-adaptive way.

Given the partition $S\in\mathbb S$ we reconstruct the regression surface with
\begin{equation}\label{prior:repulsive}
\begin{split}
f_\beta^S (x)& = \sum_{j=1}^J \beta_j \I_{I_j^S}(x), \quad \text{where}\quad (\beta_j: j\leq J)  \stackrel{ind}{\sim } g_j. \\ 
\end{split}
\end{equation} 
Regarding the prior density $g_j$, we will assume that there exist $0<c_0\leq c_1$ and  $B_0>0$ such that
\begin{equation}\label{eq:prior_g}
 c_0 \leq g_j(\beta)   \quad \forall |\beta| \leq B_0,\quad\text{and}\quad \| g_j\|_\infty \leq c_1\quad \forall j\leq J.
\end{equation}

While in Section \ref{sec:np1} we obtained the near-minimax rate under the assumption $t_1>1/2$, here we show rate-exactness {\em without} necessarily assuming $t_1>1/2$.
We are interested in bounding $\Pi (A_{\varepsilon_n}^c(\wt M) \C D_n)$ where $D_n=\{(Y_i,x_i)\}_{i=1}^n$ and
\begin{equation}\label{eq:setAM}
A_{\varepsilon_n}(\wt M) = \left\{\sup_{x\in (0,1)}\frac{  |f_\beta^S (x)  -f_0(x) | }{ \varepsilon_n(x) } \leq \wt M \right\} \quad\text{and}\quad  \varepsilon_n(x)  = \left(\frac{n}{\log n}\right)^{-\frac{t(x)}{2t(x)+1}}.
\end{equation}

For a given point $x\in[0,1]$ and a partition $S=\{I_j^S\}_{j=1}^J$ we denote with $I_x^S\in S$ the interval containing $x$ and with $R(I_x^S)$ (resp. $L(I_x^S)$) its right (resp. left) neighbor. We then define $\mathcal I(x)$ as the set of intervals   which  contain $x$  and the two neighboring intervals, i.e.
\begin{equation}\label{eq:Ix}
\mathcal I(x)=\bigcup_{S\in\mathbb S}\{I_x^S, R(I_x^S), L(I_x^S)\},
\end{equation}
 and we define, for a given $x\in[0,1]$ and some $u_1>0$,
 \begin{equation}\label{eq:omegas}
 \Omega_{n,x} (u_1) = \left\{  |n_I - n\times p_I| \leq u_1\sqrt{ \log n  \times  n\times p_I}\quad \forall I \in \mathcal I(x)\right\}
 \end{equation}
where $n_I=\sum_{i=1}^n\mathbb I(x_i\in I)$ and $p_I$ is a function of $I$ which satisfies $p_0 |I| \leq p_I \leq p_1 |I|$ for some $0<p_0\leq p_1$. Our results will be conditional on a large probability event which can be loosely regarded as a design assumption.
Namely, we require that   the cells containing the knot points $z_l$ (and the neighboring cells) are large enough in terms of the number of observations falling inside, i.e. we consider an intersection of events in \eqref{eq:omegas}
$$
 \Omega_{n}(u_1) = \bigcap_{l=1}^{N_n} \Omega_{n,z_l} (u_1).
 $$
 
 Our result below will hold on this event. If the design is regular, then $\Omega_{n}(u_1) $ holds for any $u_1>0$ and $p_0 = p_1 = 1$. If the design is random (with a density bounded away from zero) then \eqref{eq:omegas} holds with large probability if $u_1$ is large enough, as shown in Lemma \ref{lemma:random_design} (Supplemental Materials).
Unlike in Section \ref{sec:np1} where we assumed $\inf_{x\in[0,1]}t(x)> 1/2$, now we assume that $0<t(\cdot)\leq 1$  and that $t(\cdot)$ is  piecewise H\"older.
 \begin{assumption}(Piecewise H\"{o}lder)\label{ass:piecewise_holder}
 Assume that there exists a fixed partition of $[0,1]$ into $k$ intervals  say $[a_j,a_{j+1})$ (resp. $(a_j,a_{j+1}]$)  with $a_0=0$ and $a_{k+1} = 1$ such that $t(\cdot) $ is 
 $\alpha_j$-H\"older  on $(a_j, a_{j+1})$,  i.e. for $L_0>0$
 $$
 |t(x)-t(y)|\leq L_0|x-y|^{\alpha_j}\quad\text{for $x,y\in (a_j, a_{j+1})$},
 $$ 
 and such that $t(\cdot)$ is right (resp. left)  H\"older at $a_j$.
 \end{assumption}
We have the following theorem. 
\begin{theorem}\label{th:repulsive}
Consider the prior defined by \eqref{prior:repulsive}, \eqref{eq:prior_g}  and \eqref{prior:repulsive_S} and with $B>9$ and $C > 32 p_0$ and $\|f_0\|_\infty<B_0$.  Under the Assumption \ref{ass:local_holder} with $t(\cdot)$ piecewise H\"older according to the Assumption \ref{ass:piecewise_holder},  there exists $M>0$ such that 
 $$ 
E_{f_0} \left[ \1_{\Omega_{n}(u_1) } \Pi ( A_{\varepsilon_n}^c(\wt M) \C D_n) \right] = o(1). 
 $$
\end{theorem}
The proof is provided in  Section \ref{sec:proof_repulsive}.

\smallskip

Theorem \ref{th:repulsive} shows that {\em rate-exactness} can be achieved in regression {\em uniformly} over $[0,1]$ for local H\"{o}lderian functions whose exponents are piece-wise H\"{o}lder. The prior construction \eqref{prior:repulsive}, \eqref{eq:prior_g}  and \eqref{prior:repulsive_S} can be regarded as a version of variable-knot zeroth order splines.
Note that the partition in Assumption \ref{ass:piecewise_holder} needs not be known for our procedure to be valid. 

While we have presented our result in the case of univariate densities, extension to the multivariate case are possible but perhaps a bit more tedious. More interestingly,  the proving technique in Section \ref{sec:proof_repulsive} may be extended to  free-knot splines, which have typically been devised to adapt spatially but for which no proofs exist.  Finally, although the repulsive prior used in Theorem \ref{th:repulsive}  on the partition is not proved to be necessary, we believe that some form repulsion is necessary. 

\section{Discussion}\label{sec:discussion}
This work studies spatial adaptivity aspects of popular Bayesian machine learning procedures including Bayesian CART, Gaussian processes, spike-and-slab wavelet reconstructions and   variable-knot splines. We have focused on H\"{o}lderian classes where the smoothness is varying over the function domain. We have shown uniform (near)-minimax local adaptation in the supremum-norm sense in white noise as well as non-parametric regression for Bayesian CART and spike-and-slab priors. We have also provided a valid frequentist framework for uncertainty quantification with confidence set with asymptotic coverage 1 and whose width is optimal and  varies with local smoothness. We proposed a new class of repulsive partitioning priors which relate to variable-knot spline techniques and showed that they are locally rate-exact.
 
\section{Proofs of Theorems \ref{thm:np} and \ref{th:repulsive}} \label{sec:proofs}

\subsection{Proof of Theorem \ref{thm:np}} \label{pr:thm:np}
We write  $L_n=L_{max}=\lfloor\log_2 n\rfloor$ and denote with $\bT$ the set of binary trees whose  deepest internal node depth is smaller than $L_n$.
Recall the notation from Section \ref{sec:bcart1} where we denoted the set of internal tree nodes with $\mT_{int}$ and the set of external tree nodes with $\mT_{ext}$. 
Using the definition of $M_{lk}$ and $\eta_{lk}$ in  \eqref{eq:def_t_eta} and $k_l(x)$ in Lemma \ref{lemma:coef2} we first define,  for some $\bar\gamma>0$,
\begin{equation}\label{eq:def_d}
d_{l}(x)=\Big \lfloor  \log_2 \left[C_{l}(x)\left(\frac{n}{\log n}\right)^{\frac{1}{2t(x)+1}}\right]\Big\rfloor \quad\text{where}\quad  C_{l}(x)=(2M_{lk_l(x)}/\bar\gamma)^{\frac{1}{t(x)+1/2}}.
\end{equation}
It turns out that when\footnote{Note that when $\eta$ is bounded away from zero, we have $\wt d_l(x)=d_l(x)$ when $n$ is large enough.} 
\begin{equation}\label{eq:lowerb}
l\geq \wt d_{l}(x)\equiv\max\{\log_2(1/2\eta_{lk_l(x)}), d_{l}(x)\},
\end{equation}
the multiscale coefficient satisfies (from Lemma \ref{lemma:coef2})
\begin{equation}\label{eq:beta_bound}
|\beta_{lk_l(x)}^0|\leq\bar\gamma \sqrt{\frac{\log n}{n}}.
\end{equation}
Moreover, \eqref{eq:lowerb} implies that  $|\beta_{l'k_{l'}(x)}^0|\leq \bar\gamma \sqrt{\frac{\log n}{n}}$ for {\em all} $(l',k_{l'}(x))$ where $l'>l$. Indeed,
since $I_{l'k_{l'}(x)}\subset I_{lk_l(x)}$ we have $M_{l'k_{l'}(x)}\leq M_{lk_l(x)}$ and   thereby
$$
|\beta_{l'k_{l'}(x)}^0|\leq 2M_{l'k_{l'}(x)}2^{-l'(t(x)+1/2)}\leq 2 M_{lk_l(x)}2^{-l(t(x)+1/2)}\leq \bar\gamma \sqrt{\log n/n}.
$$
For a tree $\mT$, we denote with $\wt\mT_{int}$ a set of {\em pre-terminal nodes} such that both children are external nodes, i.e. 
\begin{equation}\label{eq:preterminal}
\wt\mT_{int}=\left\{(l,k)\in\mT_{int}\quad s.t.\quad \{(l+1,2k),(l+1,2k+1)\}\in\mT_{ext} \right\}.
\end{equation}
Note that for all $x\in[0,1]$ we have 
$
\wt d_{l}(x)\geq  \wt d_{l+1}(x).
$

 The main difference between the regression case and the white noise model is the dependence of parameters in the posterior distribution due to the fact that the design is not necessarily  regular. 
 Let $A_n = \{ \sup_{x\in\mX} \zeta_n(x) |f(x)-f_0(x)|>M_n\} $
and  $T$ denote a set of trees $\mT\in\bT$ that (a) {\em capture signal} and (b) {\em that are suitably small locally}.
Formally, we define the set $T$  as 
\begin{equation*}
T=\left\{\mT\in\bT:   l\leq \min_{x\in I_{lk}}\wt d_{l}(x)\quad\forall (l,k)\in\wt\mT_{int}\quad\text{and}\quad S(f_0,A;\upsilon)\subseteq \mT_{int} \right\}
\end{equation*}
for some $A,\upsilon>0$ where 
$
S(f_0,A;\upsilon)\equiv \{(l,k): |\beta_{lk}^0|>A\log^{1\vee \upsilon} n/\sqrt{n}\}
$
where $a\vee b=\max\{a,b\}$.
Going further, with $\mE(\mT)$ we denote the set of functions $f=\sum_{(l,k)\in\mT_{int}}\psi_{lk}\beta_{lk}$ that live on the tree skeleton $\mT$ and
\begin{equation}\label{eq:E}
\mE=\bigcup\limits_{\mT\in T}\, \mE(\mT)=\{f:\mT\in T\}.
\end{equation}

With $\mE$ introduced in \eqref{eq:E}, we show in Section \ref{sec:small_trees2} and Section \ref{sec:signal_catch}  that $E_{f_0}\Pi(\mE^c\C Y)\rightarrow 0$.
We can write, for $\mA$ defined in \eqref{eq:setA_np} with $P_{f_0}(\mA^c)\leq 2/ p\rightarrow 0$,
\begin{align*}
&E_{f_0}\Pi\left[f \in A_n \,\big|\,  Y\right]\leq P_{f_0}[\mA^c]+E_{f_0}\Pi[\mE^c\C Y]+ E_{f_0} \Pi\left[f \in A_n \cap \mE \,\big|\,  Y\right]\I_\mA
\end{align*}
Using the Markov's inequality, one can bound the last display above with (denoting $\mX=\{x_i:1\leq i\leq n\}$)
\begin{align}
\Pi\left[f \in A_n \cap \mE \,\big|\,  Y\right] &\leq  M_n^{-1} \int_\mE \sup_{x\in\mX} \zeta_n(x)|f(x)-f_0^{d}(x)|d\Pi(f\C Y)+M_n^{-1}B \label{Pi_E},
\end{align}
 where $B$ is the bias term defined in \eqref{eq:B_term} and is shown to be $\mathcal{O}(1)$ in Lemma \ref{lemma:bias} and 
$ 
f_0^{  d}(x)=\sum_{l\leq L_n}\sum_{k=0}^{2^l-1}\mathbb{I}[l \leq  \wt d_{l}(x)]\,\psi_{lk}(x)\beta_{lk}^0.$
 Since trees $\mT\in T$ catch large signals (according to the definition of $T$ above) we have $|\beta_{lk}^0|<A\log^{1\vee\upsilon} n/\sqrt{n}$ for $(l,k)\notin\mT_{int}$ and 
\begin{equation*}
 \sup_{x\in\mX}\left[\zeta_n(x)\sum_{(l,k)\notin \mT_{int}; l\leq \wt d_l(x)}2^{l/2}\I_{x\in I_{lk}}|\beta_{lk}^0|\right]\lesssim \frac{\log^{1\vee\upsilon} n}{\sqrt{n}}\sup_{x\in\mX}\zeta_n(x)2^{\wt d_l(x)/2}\lesssim {\log ^{1\vee \upsilon-1/2} n}.
\end{equation*}
It thereby suffices to focus on the active coordinates inside $\mT_{int}$. We now  show that on the event $\mA$
$$
\int \max_{(l,k)\in \mT_{int}}|\beta_{lk}-\beta_{lk}^0|d\Pi(\b\C\mT,Y)\lesssim{\frac{\log^{1\vee\upsilon} n}{\sqrt n}}.
$$
Set  $\Sigma_{\mT}=c_n(X_\mT'X_\mT)^{-1}$  with $c_n=g_n/(1+g_n)$ and $\mu_\mT=\Sigma_{\mT}X_\mT'[X\b_0+\bm \nu]$  we have 
$
\b_\mT\C Y\sim\mathcal{N}\left(\mu_\mT,\Sigma_{\mT}\right)
$ 
and we use Lemma 8 in \cite{castillo_rockova}  which yields for $\bar\sigma=\max \mathrm {diag} (\Sigma_\mT)$
\begin{equation}\label{eq:bound_expect}
E\|\b_\mT-\b_\mT^0\|_\infty\leq \|\mu_\mT-\b_\mT^0\|_\infty+\sqrt{2\bar\sigma^2\log |\mT_{int}|}+2\sqrt{2\pi\bar\sigma}.
\end{equation}
For the first term, we note (denoting $\|A\|_\infty=\max_i\sum_{j}|a_{ij}|$)
$$
\|\mu_\mT-\b^0_\mT\|_\infty\leq (1-c_n)\|\b^0_\mT\|_\infty+ \|\Sigma_\mT\|_\infty\|X_\mT'(X_{\bmT}\b^0_{\bmT}+ F_0-X\b_0+\bm\varepsilon)\|_\infty 
$$
From Lemma \ref{lemma:bias} and Section \ref{sec:small_trees2} we  have on the event $\mA$
$$ 
\|X_\mT'(X_{\bmT}\b^0_{\bmT}+ F_0-X\b_0+\bm\varepsilon)\|_\infty \lesssim\sqrt{n}\log^{1\vee\upsilon} n.
$$
Denoting with $a(i,\mT)$ (resp. $a(\backslash i,\mT)$) the $i^{th}$ diagonal (resp. off-diagonal) entry in the matrix $X_\mT'X_\mT$, we can write 
using the Gershgorin theorem (see e.g. \citep{matrix_book}) and Lemma \ref{lemma:eigen}
$$
\|\Sigma_\mT\|_\infty\leq \frac{c_n}{\min_{i}[a(i,\mT)-a(\backslash i,\mT)]}\leq \frac{1}{\munderbar\lambda n}.
$$ 
Next, $\bar \sigma\leq \|\Sigma_\mT\|_\infty\leq 1/(\munderbar\lambda n)$ and from \eqref{eq:bound_expect} 
we obtain $E\|\b_\mT-\b_\mT^0\|_\infty\lesssim \log^{1\vee\upsilon} n/\sqrt{n}$.
Therefore, on the event $\mA$
\begin{align*}
 A(\mT)\equiv&\int  \sup_{x\in \mathcal X}\left[\zeta_n(x)\sum_{(l,k)\in\mT_{int}} \I_{x\in I_{lk}}2^{l/2}  |\beta_{lk}-\beta_{lk}^0|\right]d\,\Pi(\b\C \mT, Y)\\
\lesssim& \int  \max_{(l,k)\in\mT_{int}}|\beta_{lk}-\beta_{lk}^0|\sup_{x\in \mathcal X}\left[\zeta_n(x) 2^{\wt d_l(x)/2}\right]d\,\Pi(\b\C \mT, Y)\nonumber\\
\lesssim& \sqrt{\frac{n}{\log n}} \int  \max_{(l,k)\in\mT_{int}}|\beta_{lk}-\beta_{lk}^0|d\,\Pi(\b\C\mT, Y) \leq B_A\times {\log^{1\vee\upsilon-1/2} n},
\end{align*}
uniformly for all $\mT\in T$ for some $B_A>0$.
We now put the pieces together. From the considerations above, we continue the calculations in \eqref{Pi_E} to obtain, on the event $\mA$,
\begin{align*}
\Pi\left[f \in A_n \cap \mE \,\big|\,  Y\right] &\leq M_n^{-1}\sum_{\mT\in T}\Pi[\mT\C Y]\int_{\mE(\mT)}   \sup_{x\in\mX}\zeta_n(x)\left|f(x)-f_0^{{d}}(x)\right|d\Pi(f\C Y,\mT) + o(1)\\
&\leq M_n^{-1}\mathcal{O}({\log^{1\vee\upsilon-1/2} n}) +o(1).
\end{align*}
The upper bound goes to zero as long as $M_n$ is strictly faster than ${\log^{1\vee\upsilon-1/2} n}$.

\subsection{Proof of Theorem \ref{th:repulsive}}\label{sec:proof_repulsive}
First, we show that  $A_{\varepsilon_n}^c(\wt M)$  is  contained in
$$ 
\bigcup_{l=1}^{N_n} \left\{ \frac{  |f_\beta^S (z_l)  -f_0(z_l) | }{ \varepsilon_n(z_l) } > \wt M/2 \right\}
$$ 
so that it suffices to focus on the discretization $\mathcal I_n$ of $[0,1]$.
To show this, we note that for all $x \in [z_l, z_{l+1})$ we have $f_\beta^S (x) =f_\beta^S (z_l)$. Next, from the Assumption\footnote{Since in Assumption \ref{ass:local_holder}, $M(\cdot)$ and $\eta(\cdot)$ are bounded, they could be regarded as constants.}
\ref{ass:local_holder} where $M(\cdot)\leq\bar M$ and $\eta(\cdot)\geq\munderbar\eta>0$, and by using \eqref{eq:z} and the Assumption \ref{ass:piecewise_holder} we obtain for a sufficiently large $n$ and a suitable $\alpha_l>0$
$$
|f_0(x) - f_0(z_l)| \leq \bar M \left(\frac{C_2 \log n}{n}\right)^{t(z_l)} \leq  \bar M \left(\frac{C_2 \log n}{n}\right)^{t(x) + L_0C_2^{\alpha_l} \left(\frac{\log n}{n}\right)^{\alpha_l}}=o(\varepsilon_n(x)).
$$
Hence 
$$
\sup_{x\in (0,1)}\frac{  |f_\beta^S (x)  -f_0(x) | }{ \varepsilon_n(x) }  \leq \max_{l\leq N_n}\frac{  |f_\beta^S (z_l)  -f_0(z_l) | }{ \varepsilon_n(z_l) }+o(1)
$$ 
and thereby
$$
\Pi(A_{\varepsilon_n}^c \C D_n ) \leq \sum_{l\leq N_n}\Pi \left(\frac{  |f_\beta^S (z_l)  -f_0(z_l) | }{ \varepsilon_n(z_l) }>\wt M/2 \right).
$$
We now focus on one particular knot value $x=z_l$ for some $l$. 
For a given partition $S$, recall that   $I_x^S$ denotes the interval in $S$ which contains $x$. We consider two types of partitions $S$ (`small-bias' versus `large-bias'), i.e. for some $M_1>0$ we distinguish between partitions $S$ satisfying
$\{ |\bar y_{I_x^S} - f_0(x) | \leq M_1\varepsilon_n\} \quad
\text{and}\quad \{ |\bar y_{I_x^S} - f_0(x) | > M_1\varepsilon_n\},$
where
$$
\bar y_{I} = \sum_{i:x_i\in I} \frac{Y_i}{n_I}\quad\text{and}\quad  n_I = \sum_{i=1}^n \I(x_i\in I).
$$
We further split the `small-bias' partitions $\{ |\bar y_{I_x^S} - f_0(x) | \leq M_1\varepsilon_n\}$ into two types (a `small cell' $I_x^S$ versus a `large cell' $I_x^S$), i.e.  
for some small  $\delta >0$ we distinguish between
\begin{equation}\label{eq:event_split2}
 \{n_{I_x^S}>s_n(\delta) \} \quad\text{and}\quad \{n_{I_x^S} \leq s_n(\delta)\}, \quad s_n(\delta) = \frac{ \delta \log n }{\varepsilon_n^2}.
 \end{equation}
 
We first prove that if $S$ is a favorable partition, i.e. if it belongs to 
\begin{equation*} 
B_n  = \left\{S:   \left\{ |\bar y_{I_x^S} - f_0(x) | \leq M_1\varepsilon_n\} \cap \{ n_{I_x^S}> s_n(\delta)\right\} \right\}
\end{equation*} 
then the conditional posterior distribution given $S$ concentrates on $\{| f_0(x) - f_\beta^S(x) | \leq  2 M_1 \varepsilon_n\}$. We then prove that the posterior probability of  the set of non-favorable partitions, i.e. $B_n^c$, goes to zero as $n$ goes to infinity.

  Recall the definition of $\mathcal I(x)$ in \eqref{eq:Ix} as the set of intervals which either contain $x$ or are neighboring intervals to the one which contains $x$.
We now define the following events for $u_0,  u_2>0$ and $ \bar \epsilon_{I} = \sum_{i:x_i\in I} {\epsilon_i}/{n_I}$ 
$$
\Omega_{n,y}(u_0) = \left\{ \forall S\in\mathbb S :  |\bar {\epsilon}_{I_x^S}| \leq u_0\sqrt{\frac{\log n}{n_{I_x^S}}}\right\},\,\, \Omega_{n,y,2}(u_2) = \left\{ \forall I \in \mathcal I(x): |\bar{\epsilon}_{I}| \leq u_2  
\sqrt{\frac{\log n}{n_I}}\right\}.
$$
Since for a given $I_x^S$ and $X=(x_1,\dots, x_n)'$  the standard Hoeffding Gaussian tail bound (see e.g. (2.10) in \cite{bernstein}) yields 
\begin{align*}
 P\left( \left. | \bar{\epsilon}_{I_x^S}| > u_0\sqrt{\frac{\log n }{n_{I_x^S}}} \,\right|\, X \right) \leq 2  \exp\left\{ -\frac{  u_0^2 \log n }{ 2 } \right\}  
 \end{align*}
and since the number of possible intervals $I_x^S$ in the definition of $\Omega_{n,y}(u_0)$ is of the order $O((n/\log n)^2)$,
we have 
$$
P(\Omega_{n,y}(u_0)^c) = o(\log n/n)\quad\text{ as soon as } u_0^2 \geq 6.
$$ 
Similarly, note that the number of intervals involved in the definition of $\Omega_{n,y,2}$ is of order $O((n/\log n)^4)$. 
By choosing $u_2^2 > 8$ we thus obtain  that $P(\Omega_{n,y,2}(u_2)^c) =o(1/n)$. 
In the following lemmata,  we will thus condition on the high-probability events $\Omega_{n,y}(u_0)$ and $\Omega_{n,y,2}(u_2)$ and we set $\Omega_n= \Omega_{n,x} (u_1)\cap \Omega_{n,y}(u_0)\cap \Omega_{n,y,2}(u_2)$.

Given the structure of the prior, for a given  partition $S$ the marginal likelihood density has a product form and is proportional to 
$$ 
m(S) = \prod_{j} m(I_j^S), \quad m(I_j^S) = \e^{-\sum_{i \in I_j^S} (Y_i - \bar y_{I_j^S})^2/2} \int_{\mathbb R}\e^{-n_{I_j^S}(\beta - \bar y_{I_j^S})^2/2} g_j(\beta)d\beta.
$$ 
We will use repeatedly the following inequality\footnote{as soon as $ |\bar y_{I_j^S}| < B_0 -\epsilon$ for some arbitrarily small but fixed $\epsilon$} 
\begin{equation}\label{bound:mj}
\frac{ c_0(1+o(1)) \e^{-\sum_{i \in I_j^S} (Y_i - \bar y_{I_j^S})^2/2} \sqrt{2\pi} }{ \sqrt{n_{I_j^S}}} \leq m(I_j^S) \leq \frac{ c_1 \e^{-\sum_{i \in I_j^S} (Y_i - \bar y_{I_j^S})^2/2}\sqrt{2\pi} }{ \sqrt{n_{I_j^S}}}. 
\end{equation}

\begin{lemma}\label{lem:thetax}
Assume   the prior \eqref{prior:repulsive} with \eqref{prior:repulsive_S} and \eqref{eq:prior_g}. For any $a>0$ and if $M_1\geq  {2} a/\sqrt{\delta}$ then 
$$ 
 E \left[ \1_{ {\Omega_n}}  \Pi (   \{  | f_0(x) - f_\beta^S(x) | >  2 M_1 \varepsilon_n\} \cap B_n | D_n ) \right] \lesssim  n^{-a}. 
$$ 
\end{lemma}
\begin{proof}[Proof of Lemma \ref{lem:thetax}]
If  $S \in B_n$ and $| f_0(x) - f_\beta^S(x) | > 2M_1 \varepsilon_n$ then  we have
 $$
 |\bar y_{I_x^S} - f_\beta^S(x)  | \geq | f_0(x) - f_\beta^S(x) |  - |\bar y_{I_x^S} - f_0(x)| \geq M_1\varepsilon_n.
 $$
Using \eqref{bound:mj}, we then have
 \begin{equation*} 
\begin{split}
\Pi \left( | f_0(x) - f_\beta^S(x) | > 2M_1 \varepsilon_n | D_n, S\right) &\leq  \frac{2c_1 \sqrt{n_{I_x^S}}  }{ c_0 \sqrt{2\pi}(1+o(1))}\int\limits_{| \beta - \bar y_{I_x^S}| >M_1	\varepsilon_n}  \e^{ - \frac{n_{I_x^S}}{2}(\beta- \bar y_{I_x^S})^2}g(\beta)d\beta\\
& \leq \frac{ 2c_1\sqrt{n_{I_x^S}}  }{ c_0 \sqrt{2\pi}(1+o(1))} \exp\left\{ -\frac{ n_{I_x^S} M_1^2\varepsilon_n^2}{2}\right\} \\
&\lesssim   \exp\left\{ - \delta M_1^2 \log n/4 \right\}= o(n^{-a})
\end{split}  
\end{equation*} 
if $\delta M_1^2 >4a$.  
\end{proof}

 We now prove that the unfavorable partitions have posterior probability going to 0. Using  Lemma \ref{lem:mIx} (below), with $a>1$ we obtain  on $\Omega_n$ that 
 $$
 \Pi \left( S:  \{ |\bar y_{I_x^S} - f_0(x) | > M_1\varepsilon_n\} \cap \{ n_{I_x^S}> s_n(\delta) \}  \C  D_n \right)  = o_p(n^{-1})
 $$
 and that, using Lemma \ref{lem:small:nx} (below),   
 $ 
\Pi \left( S:  \{ n_{I_x^S}\leq  s_n(\delta) \}  \C  D_n \right)  = o_p(n^{-1}).
$ 
Combining these two results with Lemma \ref{lem:thetax}, we then have on $\Omega_n$
 \begin{equation*} 
\begin{split}
\Pi &\left( | f_0(x) - f_\beta^S(x) | > 2M_1 \varepsilon_n | D_n \right) \leq \sum_{S \in B_n} \Pi \left( | f_0(x) - f_\beta^S(x) | > 2M_1 \varepsilon_n   | D_n, S \right) \Pi(S|D_n) \\
& + 
\Pi \left( S:  \{ |\bar y_{I_x^S} - f_0(x) | > M_1\varepsilon_n\} \cap \{ n_{I_x^S}> s_n(\delta) \}  \C  D_n \right) + \Pi \left( S:  \{ n_{I_x^S}\leq  s_n(\delta) \}  \C  D_n \right) \\
&= o_p(n^{-1}).
\end{split}  
\end{equation*}

\begin{lemma}\label{lem:mIx}
Assume   the prior \eqref{prior:repulsive} with \eqref{prior:repulsive_S} and \eqref{eq:prior_g}. Let $x \in (0,1)$, 
then for all $a,u_0,u_1, u_2>0$ and  $\delta>0$ small enough  we have for  $M_1> 2u_0/\sqrt{\delta}$
$$ 
 E \left[\I_{\Omega_n} \Pi \left( S:  \{ |\bar y_{I_x^S} - f_0(x) | > M_1\varepsilon_n\} \cap \{ n_{I_x^S}> s_n(\delta)\}  \C  D_n \right) \right] = O(n^{-a}).
$$ 
\end{lemma}

\begin{lemma}\label{lem:small:nx}
Assume  the prior \eqref{prior:repulsive} with \eqref{prior:repulsive_S} and \eqref{eq:prior_g}. With $ \delta>0$ as in Lemma \ref{lem:mIx} we have
$$ 
E \left[ \1_{\Omega_{n,x}(u_1)} \Pi \big( \left\{  n_{I_x^S} \leq s_n(\delta)\right\} \C D_n \big) \right]= o(1/n).
$$
\end{lemma}

Lemma \ref{lem:mIx}  is proved by showing that if $ |\bar y_{I_x^S} - f_0(x) |$ and $ n_{I_x^S}> s_n(\delta)$ then the partition has much smaller posterior probability than the one obtained by splitting $I_x$  into smaller intervals. 
The proof of Lemma \ref{lem:mIx} is given below while the proof of Lemma \ref{lem:small:nx} is given in Section \ref{pr:small:nx} of the Supplementary Material \cite{rockova:rousseau:supp}.  The idea of the proof of Lemma \ref{lem:small:nx} is that  partitions verifying $\{ n_{I_x^S}> s_n(\delta)\}$ have either much smaller  probability than the one resulting from merging $I_x^S$ with a neighboring interval, say $I_{x,1}$, or much smaller  probability than the one resulting from splitting $I_{x,1}^S$  into smaller intervals. The latter result comes from the fact that if $I_{x,1}^S$ is too large then there is a point $x_1$ in $I_{x,1}^S$, such that $ |\bar y_{I_{x,1}^S} - f_0(x_1) | > M_0\varepsilon_n(x_1) $  and $ n_{I_{x,1}^S}> s_n(\delta_1)$ for some appropriate values $M_0, \delta_1$ and Lemma \ref{lem:mIx} can then be used. 

\begin{proof}[Proof of Lemma \ref{lem:mIx}]
Throughout the rest of the proof, we suppress the index $S$ when referring to intervals $I_x^S$ or $I_j^S$.  
On the event $\Omega_{n,y}(u_0)$, and if $n_{I_x} > s_n(\delta)$ for a given $\delta$, we have {for $\bar \beta_{0,I_x}=\sum_{x_i\in I_x}f_0(x_i)/n_{I_x}$}
$$ |\bar y_{I_x}  - \bar \beta_{0,I_x} |=| \bar \epsilon_{I_x} | \leq u_0 \frac{ \sqrt{\log n}}{\sqrt{n_{I_x}}} \leq   \frac{u_0 \varepsilon_n }{ \sqrt{\delta} } \leq \frac{ M_1 \varepsilon_n}{ 2} $$
as soon as $M_1 > 2 u_0/\sqrt{\delta} $. In particular if  $ |\bar y_{I_x}  - f_0(x)| > M_1\varepsilon_n$ then 
we have from  Assumption \ref{ass:local_holder} that  as soon as $|I_x| \leq \munderbar\eta$, 
$$
M_1 \varepsilon_n/2\leq |\bar \beta_{0,I_x}-f_0(x) | \leq M |I_x|^{t(x)}
$$
so that in all cases 
 $|I_x|\geq (M_1 \varepsilon_n/2M)^{1/t(x)}$.

Since the cell $I_x$ has a large bias, we compare the partition $S$ with a partition obtained by splitting $I_x$ into 2 or 3 intervals, say $I_1, I_2$, and possibly $ I_3$ if $x$ is too far from the boundary of $I_x$. We do the splitting\footnote{Without loss of generality, we can assume that  cutting an interval of such a size is possible otherwise we would replace it with $|I_1| =  (\tau M_1 \varepsilon_n/2M )^{1/t(x)}(1 + o(1))$ which makes no difference. }  so that $x \in I_1$ and $|I_1| =  (\tau M_1 \varepsilon_n/2M )^{1/t(x)}$ for some $\tau <1$.  We choose also $\tau>0$ small so that both $|I_2|, |I_3| \geq  (\tau M_1 \varepsilon_n/2M )^{1/t(x)}$.
 Then on $\Omega_{n,y}(u_0)$, 
 $$ 
 |\bar \beta_{0,I_1}  - f_0(x)|\leq \tau M_1\varepsilon_n/2 \quad \text{and} \quad |\bar y_{I_1}  - f_0(x)|  \leq \tau M_1\varepsilon_n/2 + \frac{ u_0 \sqrt{ \log n}}{ \sqrt{n_{I_1}}}.  
 $$

 In the following, we write the computations in the case where we have split $I_x$ into 3 intervals. Computations for the case of 2 intervals can be derived similarly. 
Note that, by construction, $|I_2| \geq |I_1|$ and $ |I_3|\geq |I_1|$. In addition, on the event $ \Omega_{n,x} (u_1)$ defined in \eqref{eq:omegas} we have  $n_{I_j} \geq n p_0 |I_j|/2$ for $j=1,2,3$. 
Hence, there exists a constant $C_0>0$ such that   
\begin{equation}\label{eq:ubound}
\frac{ u_0 \sqrt{ \log n} }{ \sqrt{n_{I_1} } } \leq \frac{ u_0 C_0 }{ (\tau M_1
)^{\frac{1}{2t(x)} }} \varepsilon_n\leq M_1 \varepsilon_n/2
\end{equation}
by choosing $M_1$ large enough so that 
$ |\bar y_{I_1}  - f_0(x)|  \leq  M_1\varepsilon_n( 1 + \tau )/2. $
On the event $\Omega_{n,y,2}(u_2)$, for all $u_2>0$,
we have  
$$
|\bar y_{I_1}| \leq  |f_0(x)| + \epsilon\quad\text{and}\quad   |\bar y_{I_2}| \leq  \|f_0\|_\infty + \epsilon\leq B_0
$$ 
for any $\epsilon>0$ small when $n$ is large enough  since $\|f_0\|_\infty<B_0$. 
Hence using \eqref{bound:mj},
\begin{equation*} 
\begin{split}
\frac{ m(I_x) }{ m(I_1) m(I_2)m(I_3) } &\leq \frac{ 2c_1\sqrt{n_{I_1} n_{I_2}n_{I_3}} }{ 2\pi c_0^3\sqrt{ n_{I_x}} } \exp\left(-\frac{ \sum_{i \in I_x} (y_i - \bar y_{I_x})^2}{2} +\frac{\sum_{j=1}^3\sum_{i \in I_j} (y_i - \bar y_{I_j})^2}{2}\right)  \\
& = \frac{ 2c_1\sqrt{n_{I_1} n_{I_2}n_{I_3}} }{ 2\pi c_0^3\sqrt{ n_{I_x}} } \exp\left({-\sum_{j=1}^3\frac{ n_{I_j}(\bar y_{I_x} - \bar y_{I_j})^2 }{2}}\right).
\end{split}
\end{equation*} 
Moreover, we have
$$
|\bar y_{I_x} -   \bar y_{I_1}| > | \bar y_{I_x} -   f_0(x) | - | f_0(x) -\bar y_{I_1}| \geq M_1(1-\tau) \varepsilon_n/2 \geq M_1\varepsilon_n/4
$$
by choosing $\tau \leq 1/2$. 
Finally,  by noting that $n_I\asymp n|I|$ on the event $\Omega_{n,x}(u_1)$ we obtain  
\begin{equation*} 
\begin{split}
\frac{ m(I_x) }{ m(I_1) m(I_2)m(I_3) } &\lesssim  n  \sqrt{|I_1| |I_2|} \exp\left({-n_{I_1} M_1^2 \varepsilon_n^2  /8}\right)\\
&\leq n^{-a}
\end{split}
\end{equation*} 
for all $a>0$ by choosing $n_{I_1} M_1^2 \varepsilon_n^2  > 8 (a+1)$. Since  
$$ n_{I_1} M_1^2 \varepsilon_n^2 \geq \frac{ n p_0 |I_1| M_1^2 \varepsilon_n^2 }{2 } \geq 
(\tau / (2M))^{1/t(x)}M_1^{(2t(x)+1)/t(x) }\log n/2 $$ 
this is true as soon as $M_1 \geq 2u_0/\sqrt{\delta} $ and $\delta$ is small enough. Then we have
\begin{equation*} 
\begin{split}
Z_x(S)\equiv \frac{ m(I_x)|I_x|^B }{ m(I_1) m(I_2)m(I_3)|I_1|^B |I_2|^B|I_3|^B} &\lesssim  n (|I_1| |I_2|)^{-(B-1/2)} \exp\left({-n_{I_1} M_1^2 \varepsilon_n^2  /8}\right) \\
& \lesssim (|I_1| |I_2|)^{-B}n^{-a}  \leq n^{-a/2},
\end{split}
\end{equation*} 
as soon as $a> 4 B / (2t(x) +1) $ since $|I_1| |I_2| \gtrsim \epsilon_n^{1/t(x)}$ . 
This implies that on  $\Omega_n$ we have for
 $$
 \Pi_1\equiv \Pi ( S:  \{ |\bar y_{I_x} - f_0(x) | >M_1 \varepsilon_n\} \cap \{ n_{I_x} >s_n(\delta) \} \C D_n )
 $$
and $\mathbb I_1(S)\equiv \I\{S: |\bar y_{I_x} - f_0(x) | >M_1 \varepsilon_n\}$ and $\mathbb I_2(S)\equiv  \I\{S: n_{I_x}>s_n\} $ the following bound 
\begin{equation*}
\begin{split}
\Pi_1 & = 
\frac{ \sum\limits_{S=S'\cup I_x} \mathbb I_1(S) \times \mathbb I_2(S) \times m(S')\times  m(I_x)  \times \pi_S(S'\cup I_x) }{ \sum\limits_{S=S'\cup I_x}   m(S') \times m(I_x)  \times \pi_S(S'\cup I_x)  } \\
&\leq \frac{ \sum\limits_{S=S'\cup I_x}\mathbb I_1(S) \times \mathbb I_2(S) \times m(S') m(I_1) m(I_2) m(I_3) \times \pi_S( S'\cup I_1\cup I_2\cup I_3)\times Z_x(S)}{ 
\sum\limits_{S=S'\cup I_1\cup I_2\cup I_3} \I_2(S)\times  m(S')\times  m(I_1) m(I_2)m(I_3) \times \pi_S( S'\cup I_1\cup I_2\cup I_3)} \\
&\leq  C n^{-a/2}.\hspace{11.5cm}\qedhere 
\end{split} \end{equation*}


\end{proof}

\section*{Acknowledgements}

The project leading to this work has received funding from the European Research Council
(ERC) under the European Union’s Horizon 2020 research and innovation programme (grant agreement No 834175).
The authors also gratefully acknowledge  support from the James S. Kemper Foundation Faculty Research Fund at  the University of Chicago Booth School of Business and the National Science Foundation (DMS:1944740).

\clearpage

\begin{frontmatter}
\title{Ideal Bayesian  Spatial Adaptation: supplementary material}
\runtitle{Bayesian Trees are Spatially Adaptive}
\begin{aug}

\author{\fnms{Veronika} \snm{Ro\v{c}kov\'{a} }\thanksref{t2,m2}\ead[label=e2]{veronika.rockova@chicagobooth.edu}}
\author{and \fnms{Judith} \snm{Rousseau}\thanksref{m1}\ead[label=e1]{judith.rousseau@stats.ox.ac.uk}}

\thankstext{t2}{
The author gratefully acknowledges  support from the James S. Kemper Foundation Faculty Research Fund at  the University of Chicago Booth School of Business and the National Science Foundation (DMS:1944740).}


\affiliation{University of Chicago\thanksmark{m2}}
\affiliation{University of Oxford \thanksmark{m1}}

\address{University of Chicago\\ Booth School of Business\\ 5807 S. Woodlawn Avenue\\ Chicago, IL, 60637 \\ USA\\
 \printead{e2}}

\address{University of Oxford \\
Department of Statistics\\
 24-29 St Giles'\\ Oxford OX1 3LB \\ United Kingdom\\
  \printead{e1}}

\end{aug}

\begin{abstract}
This supplementary material contains the remaining proofs of the main text. In particular, the proofs  for Theorems \ref{thm:one} and \ref{thm:uq} under the white noise model are presented in Section  \ref{sec:pr:whitenoise}, the proof of the non-spatial adaptation for common classes of hierarchical Gaussian process prior is presented in Section \ref{sec:proof_lb} and some of the technical lemmas used in the proof of Theorem \ref{thm:np} are presented in Section \ref{sec:proof_thm_np}. In Section \ref{pr:small:nx}, we provide the proof of Lemma \ref{lem:small:nx} used in the proof of Theorem \ref{th:repulsive} and Section \ref{sec:aux} contains some auxiliary results. 
\end{abstract}

\begin{keyword}[class=MSC]
\kwd[Primary ]{62G20, 62G15}
\end{keyword}

\begin{keyword}
\kwd{Bayesian CART}
\kwd{Partitioning Priors} 
\kwd{Supremum Norm}
\kwd{Spike-and-Slab}
\kwd{Spatial Adaptation}
\end{keyword}

\end{frontmatter}

\tableofcontents


\section{Proofs for the White Noise Model } \label{sec:pr:whitenoise}
\subsection{Proof of Theorem \ref{thm:one}}\label{sec:proof_thm1}
The proof is similar to the proof in the regression case but is simpler. For the sake of self-sufficiency, we recall some the definitions used in the proof of Theorem \ref{thm:np}, see also Section \ref{pr:thm:np}.  We write  $L_n=L_{max}=\lfloor\log_2 n\rfloor$ and denote with $\bT$ the set of binary trees whose  deepest internal node depth is smaller than $L_n$.
Recall the notation from Section \ref{sec:bcart1} where we denoted the set of internal tree nodes with $\mT_{int}$ and the set of external tree nodes with $\mT_{ext}$. 
Using the definition of $M_{lk}$ and $\eta_{lk}$ in  \eqref{eq:def_t_eta} and $k_l(x)$ in Lemma \ref{lemma:coef2} we first define,  for some $\bar\gamma>0$,
\begin{equation}\label{eq:def_d}
d_{l}(x)=\Big \lfloor  \log_2 \left[C_{l}(x)\left(\frac{n}{\log n}\right)^{\frac{1}{2t(x)+1}}\right]\Big\rfloor \quad\text{where}\quad  C_{l}(x)=(2M_{lk_l(x)}/\bar\gamma)^{\frac{1}{t(x)+1/2}}.
\end{equation}
It turns out that when\footnote{Note that when $\eta$ is bounded away from zero, we have $\wt d_l(x)=d_l(x)$ when $n$ is large enough.} 
\begin{equation}\label{eq:lowerb}
l\geq \wt d_{l}(x)\equiv\max\{\log_2(1/2\eta_{lk_l(x)}), d_{l}(x)\},
\end{equation}
the multiscale coefficient satisfies (from Lemma \ref{lemma:coef2})
\begin{equation}\label{eq:beta_bound}
|\beta_{lk_l(x)}^0|\leq\bar\gamma \sqrt{\frac{\log n}{n}}.
\end{equation}
Moreover, \eqref{eq:lowerb} implies that  $|\beta_{l'k_{l'}(x)}^0|\leq \bar\gamma \sqrt{\frac{\log n}{n}}$ for {\em all} $(l',k_{l'}(x))$ where $l'>l$. Indeed,
since $I_{l'k_{l'}(x)}\subset I_{lk_l(x)}$ we have $M_{l'k_{l'}(x)}\leq M_{lk_l(x)}$ and   thereby
$$
|\beta_{l'k_{l'}(x)}^0|\leq 2M_{l'k_{l'}(x)}2^{-l'(t(x)+1/2)}\leq 2 M_{lk_l(x)}2^{-l(t(x)+1/2)}\leq \bar\gamma \sqrt{\log n/n}.
$$
For a tree $\mT$, we denote with $\wt\mT_{int}$ a set of {\em pre-terminal nodes} such that both children are external nodes, i.e. 
\begin{equation}\label{eq:preterminal}
\wt\mT_{int}=\left\{(l,k)\in\mT_{int}\quad s.t.\quad \{(l+1,2k),(l+1,2k+1)\}\in\mT_{ext} \right\}.
\end{equation}
Note that for all $x\in[0,1]$ we have 
$$
\wt d_{l}(x)\geq  \wt d_{l+1}(x).
$$ 
In the sequel, $T$ denotes a set of trees $\mT\in\bT$ that (a) {\em capture signal} and (b) {\em that are suitably small locally}.
Formally, we define the set $T$  as 
\begin{equation}\label{eq:T}
T=\left\{\mT\in\bT:   l\leq \min_{x\in I_{lk}}\wt d_{l}(x)\quad\forall (l,k)\in\wt\mT_{int}\quad\text{and}\quad S(f_0,A;\upsilon)\subseteq \mT_{int} \right\}
\end{equation}
for some $A,\upsilon>0$ where 
\begin{equation}\label{eq:setS}
S(f_0;A;\upsilon)\equiv \{(l,k): |\beta_{lk}^0|>A\log^{1\vee \upsilon} n/\sqrt{n}\}.
\end{equation}
Going further, with $\mE(\mT)$ we denote the set of functions $f=\sum_{(l,k)\in\mT_{int}}\psi_{lk}\beta_{lk}$ that live on the tree skeleton $\mT$ and
\begin{equation}\label{eq:E}
\mE=\bigcup\limits_{\mT\in T}\, \mE(\mT)=\{f:\mT\in T\}.
\end{equation}

First, we show that $E_{f_0}\Pi(\mE^c\C Y)\rightarrow 0$. 
To begin, in Section \ref{sec:small_trees} below we show that the posterior concentrates on locally small trees.

\subsubsection{Posterior Concentrates on $\mE$}\label{sec:small_trees}
Our considerations will be conditional on the event
\begin{equation}\label{eq:event_A}
\mA_n=\left\{\max_{-1\leq l\leq L_n, 0\leq k<2^l} \epsilon_{lk}^2\leq 2\log(2^{L_n+1})\right\}
\end{equation}
which has a large probability in the sense that $P(\mA_n^c)\lesssim (\log n)^{-1}$. 
\begin{lemma}\label{lemma:locally_small}
Let $\wt d_{l}(x)$ be as in \eqref{eq:lowerb}. For the Bayesian CART prior from Section \ref{sec:bcart1} with a split probability $p_l=(1/\Gamma)^l$  we have,  on the event $\mA_n$ in \eqref{eq:event_A},  for $\Gamma>0$ large enough
\begin{equation}
\Pi\left[\mT: \exists (l,k)\in \wt\mT_{int} \,\,\,\, s.t.\,\,\,\, l>\min_{x\in I_{lk}}\wt d_{l}(x)\C Y\right]\rightarrow 0.\label{eq:locally_small}
\end{equation}
\end{lemma}

\begin{proof}
We can write
\begin{equation*}\label{eq:bound_prob1}
\Pi\left[\mT: \exists (l,k)\in \wt\mT_{int} \,\,\,\, s.t.\,\,\,\, l>\min_{x\in I_{lk}} \wt  d_{l}(x)\C Y\right]\leq \sum_{l\leq L_n}\sum_{k=0}^{2^l-1}\mathbb{I}[l>\min_{x\in I_{lk}} \wt d_{l}(x)]\Pi[(l,k)\in\wt\mT_{int}\C Y].
\end{equation*}
We denote with $\bT_{lk}$ the set of all trees $\mT$ such that $(l,k)\in\wt\mT_{int}$. Then
\begin{equation}\label{eq:pi_lk}
\Pi\left[ (l,k) \in \wt\mT_{int}  \C Y\right]  = \frac{ \sum_{\mathcal T\in \bT_{lk}} W_Y(\mT) }{ \sum_{\mathcal T} W_Y(\mT)}
\end{equation}
where, for $\b_\mT=(\beta_{lk}:(l,k)\in\mT_{int})'$ and  $\Y_\mT=(Y_{lk}:(l,k)\in\mT_{int})'$,
$$
W_Y(\mT)=\Pi(\mT)N_Y(\mT)\quad\text{with}\quad N_Y(\mT)=\int\e^{-\frac{n}{2}\|\b_\mT\|_2^2+n\Y_\mT'\b_\mT}\pi(\b_\mT)d\b_\mT.
$$
For a tree $\mT\in\bT_{lk}$, denote with $\mT^-$ the smallest subtree of $\mT$ that does {\em not} contain $(l,k)$ as a pre-terminal node, i.e. $\mT^-$ is obtained from $\mT$ by   turning $(l,k)$ into a terminal node.  We can then rewrite \eqref{eq:pi_lk} as
\begin{equation}\label{eq:prob_bound}
\Pi\left[ (l,k) \in \wt\mT_{int}  \C Y\right]  = \frac{ \sum_{\mathcal T\in \bT_{lk}} \frac{W_Y(\mT)}{W_Y(\mT^-)}W_Y(\mT^-)}{ \sum_{\mathcal T} W_Y(\mT)}.
\end{equation}
Assuming an independent product prior $\beta_{lk}\iid\mathcal{N}(0,1)$, we have
\begin{equation}\label{eq:ratio_W}
 \frac{W_Y(\mT)}{W_Y(\mT^-)}=\frac{\Pi(\mT)}{\Pi(\mT^-)} \frac{\e^{\frac{n^2}{2(n+1)}Y_{lk}^2}}{\sqrt{n+1}}.
\end{equation}
See Section 3.1 in \cite{castillo_rockova} for details on this derivation.
 Since $(l,k)$ is such that $l\geq \wt d_{l}(x)$ for some $x\in I_{lk}$,  we have $|\beta_{l\,k_l(x)}|\leq\bar\gamma\sqrt{\log n/n}$  from \eqref{eq:beta_bound} and thereby $Y_{lk}^2=(\beta_{lk}^0+\frac{1}{\sqrt{n}}\epsilon_{lk})^2\leq C_y\, \log n/n$ on the event $\mA_n$ for some $C_y>0$.
    Next, the prior ratio (under the Galton-Watson process prior) equals  
  $$
  \frac{\Pi(\mT)}{\Pi(\mT^-)}=\frac{p_l(1-p_{l+1})^2}{1-p_{l}}. 
   $$
 For $p_l=\Gamma^{-l}\leq 1/2$, we can bound this from above with $2\Gamma^{-l}$. 
Since for each $(\ell, k)$,  the mapping $\mT\rightarrow\mT^-$ is injective, we can bound \eqref{eq:prob_bound} with 
\begin{equation}\label{eq:aux33}
2\,\Gamma^{-l}\e^{C_y/2\log n} \frac{ \sum_{\mathcal T\in \bT_{lk}^-} W_Y(\mT)}{ \sum_{\mathcal T} W_Y(\mT)}\leq 2\,\Gamma^{-l}\e^{C_y/2\log n},
\end{equation}
where $\bT_{lk}^-$ corresponds to trimmed trees inside $\bT_{lk}$ whose pre-terminal node $(l,k)$ has been turned into a terminal node.
Writing  $\bar d_{lk}=\min\limits_{x\in I_{lk}}\wt d_l(x)$ and $\bar d=\min\limits_{0\leq l\leq L_n}\min\limits_{0\leq k<2^l} \bar  d_{lk}$, we can then bound
the probability in \eqref{eq:locally_small} with, for $\Gamma>2$ 
$$
2\,\e^{C_y/2\log n}\sum_{l=\bar d}^{L_n}\Gamma^{-l}\sum_{k=0}^{2^l-1}\mathbb{I}[l>\min_{x\in I_{lk}} \wt  d_{l}(x)]= \e^{C_y/2\log n}\sum_{l=\bar d}^{L_n}(\Gamma/2)^{-l}
\lesssim \e^{C_y/2\log n-\bar d\log(\Gamma/2)}
$$
Since $M(\cdot)$  and $\eta(\cdot)$ are bounded away from zero and $t(x)\geq t_1$ (see Assumption \ref{ass:local_holder}),  for a sufficiently large $n$ we  have 
$\wt d_l(x)=d_l(x)$ for all $x\in[0,1]$   and 
$$
\bar d\geq \underbar{C}+ \frac{1}{3}\log n-\frac{1}{3}\log\log n, 
$$
where  $0<\underbar{C}=\min\limits_{0\leq l\leq L_n}\min\limits_{0\leq k<2^l}\min\limits_{x\in I_{lk}}C_{l}(x)$. For a sufficiently large $\Gamma$, the right side goes to zero. \qedhere

\end{proof}

Next, with our Bayesian CART prior  we can deploy Lemma 2 of \cite{castillo_rockova} to find that, on the event $\mA_n$, 
\begin{equation}\label{eq:catch_signal}
\Pi\left[ \mT: S(f_0;A;1)\not\subseteq \mT_{int} \C Y\right]\rightarrow 0\quad\text{as $n\rightarrow\infty$}.
\end{equation}
where $S(f_0;A;1)$ was defined in \eqref{eq:setS}. We can thus conclude, together with Lemma \ref{lemma:locally_small} above, that $E_{f_0}\Pi(\mE^c\C Y)\rightarrow 0$.

\subsubsection{Controlling the Bias Term}\label{sec:bias}
The next step in the proof is to show that the class of trees $T$ in \eqref{eq:T} are good approximators of locally H\"{o}lder functions. 
\begin{lemma}\label{lemma:bias}
Let $f_0$ satisfy Assumption \ref{ass:local_holder}   and let  $\wt d_{l}(x)$ be as in \eqref{eq:lowerb}.
We define the {\em local bias}  as
\begin{equation}\label{eq:bias}
f_0^{\backslash  d}(x)=\sum_{l\leq L_n}\sum_{k=0}^{2^l-1}\mathbb{I}[l> \wt d_{l}(x)]\,\psi_{lk}(x)\beta_{lk}^0.
\end{equation}
With $\zeta_n(x)=(n/\log n)^{t(x)/[2t(x)+1]}$, the local bias  is uniformly small in the sense that 
\begin{equation}\label{eq:B_term}
B\equiv \sup_{x\in[0,1]}\left[\zeta_n(x)|f_0^{\backslash d}|\right]\leq  \bar C\quad\text{for some $\bar C>0$}.
\end{equation}
\end{lemma}
\proof
Using Lemma  \ref{lemma:coef2} and assuming $M(x)\leq \bar M$ we have for some $C_1>0$
\begin{align*}
B\equiv \sup_{x\in[0,1]}\left[\zeta_n(x)|f_0^{\backslash d}|\right]&\leq \sup_{x\in[0,1]}\left[\zeta_n(x)\sum_{l\leq L_n}\sum_{k=0}^{2^l-1}2^{l/2}\I[\l>\wt d_{l}(x)]|\beta_{lk}^0|\right]
\label{eq:B}\\
&\leq 2\bar M \sup_{x\in[0,1]}\left[\zeta_n(x)\sum_{l> \wt d_l(x)} 2^{-lt(x)}\right]\leq 2\bar M\, C_1\,\sup_{x\in[0,1]}\left[\zeta_n(x) 2^{-\wt d_l(x)t(x)}\right].
\end{align*}
From the definition of $\wt d_l(x)$ and $C_l(x)$ in \eqref{eq:lowerb}  and \eqref{eq:def_d},  we have under Assumption \ref{ass:local_holder} for some $\bar C>0$
$$
2^{-\wt d_l(x)t(x)}\leq \left(C_{l}(x)\right)^{-t(x)}\left(\frac{\log n}{n}\right)^{\frac{t(x)}{2 t(x)+1}}\leq \frac{\bar C}{2\bar M\,C_1} \left(\frac{n}{\log n}\right)^{\frac{t(x)}{2 t(x)+1}}.
$$

\subsubsection{The Main Proof}\label{sec:main_proof_wn}
With $\mE$ introduced in \eqref{eq:E}, we have shown in Section \ref{sec:small_trees} that $E_{f_0}\Pi(\mE^c\C Y)\rightarrow 0$.
We can then write, for $\mA_n$ introduced in \eqref{eq:event_A},
\begin{align*}
&E_{f_0}\Pi\left[f: \sup_{x\in[0,1]} \zeta_n(x) |f(x)-f_0(x)|>M_n\,\big|\,  Y\right]\leq P_{f_0}[\mA_n^c]+E_{f_0}\Pi[\mE^c\C Y]+ E_{f_0} \Pi_{\mE}\I_{\mA_n}
\end{align*}
where 
\begin{equation}\label{eq:pi_me}
\Pi_\mE\equiv \Pi\left[f\in\mE: \sup_{x\in[0,1]}\zeta_n(x)|f(x)-f_0(x)|>M_n\,\big|\,  Y\right].
\end{equation}
Using the Markov's inequality, one can bound the display above with 
\begin{align*}\
\Pi_\mE&\leq M_n^{-1}\int_\mE  \sup_{x\in[0,1]} \zeta_n(x)|f(x)-f_0(x)|d\Pi(f\C Y)\\
&\leq  M_n^{-1} \int_\mE \sup_{x\in[0,1]} \zeta_n(x)|f(x)-f_0^{d}(x)|d\Pi(f\C Y)+M_n^{-1}B \label{Pi_E},
\end{align*}
 where  $f_0^{d}=f_0-f_0^{\backslash d}$ with  $f_0^{\backslash d}$ introduced in \eqref{eq:bias} and where  $B$ was defined in \eqref{eq:B_term} and was shown to be $\mathcal{O}(1)$ in Lemma \ref{lemma:bias}.
We now focus on the integrand in the last display above.  For a function $f\in\mE(\mT)$ supported on $\mT\in T$ we have  
\begin{equation}\label{eq:bound1}
 \left|f(x)-f_0^{ d}(x)\right|\leq  \sum_{(l,k)\in\mT_{int}} \I_{x\in I_{lk}}2^{l/2}  |\beta_{lk}-\beta_{lk}^0|+\sum_{(l,k)\notin \mT_{int}; l\leq \wt d_l(x)} \I_{x\in I_{lk}}2^{l/2}|\beta_{lk}^0|.
\end{equation}
We now focus on the second term above. Since trees $\mT\in T$ catch large signals (definition of $T$ in \eqref{eq:T}), we have $|\beta_{lk}^0|<A\log n/\sqrt{n}$ for $(l,k)\notin\mT_{int}$ and thereby
\begin{equation}\label{eq:bound2}
 \sup_{x\in[0,1]}\left[\zeta_n(x)\sum_{(l,k)\notin \mT_{int}; l\leq \wt d_l(x)}2^{l/2}\I_{x\in I_{lk}}|\beta_{lk}^0|\right]\lesssim \frac{\log n}{\sqrt{n}}\sup_{x\in[0,1]}\zeta_n(x)2^{\frac{\wt d_l(x)}{2}}\lesssim \sqrt{\log n}.
\end{equation}
Above, we have used the fact that (for $M(x)\leq \bar M$ and $t_1\leq t(x)\leq 1$)
$$
\zeta_n(x)2^{\wt d_l(x)/2}\leq (2\bar M/\bar\gamma)^{1/[2t(x)+1]} \left(\frac{n}{\log n}\right)^{1/2}\lesssim \left(\frac{n}{\log n}\right)^{1/2}. 
$$
Regarding the first term in \eqref{eq:bound1}, we can write for a given tree $\mT\in T$
\begin{align}\label{eq:AT}
 A(\mT)\equiv&\int  \sup_{x\in[0,1]}\left[\zeta_n(x)\sum_{(l,k)\in\mT_{int}} \I_{x\in I_{lk}}2^{l/2}  |\beta_{lk}-\beta_{lk}^0|\right]d\,\Pi(\b\C \mT, Y)\\
\lesssim& \int  \max_{(l,k)\in\mT_{int}}|\beta_{lk}-\beta_{lk}^0|\sup_{x\in[0,1]}\left[\zeta_n(x) 2^{\wt d_l(x)/2}\right]d\,\Pi(\b\C \mT, Y)\nonumber\\
\lesssim& \sqrt{\frac{n}{\log n}} \int  \max_{(l,k)\in\mT_{int}}|\beta_{lk}-\beta_{lk}^0|d\,\Pi(\b\C\mT, Y).\nonumber
\end{align}
According to Lemma 3 of \cite{castillo_rockova}, the integral above is bounded by $C'\sqrt{\log n/n}$, which implies $A(\mT)\leq B_A$ uniformly for all $\mT\in T$ for some $B_A>0$.
We now put the pieces together. From the considerations above, we continue the calculations in \eqref{Pi_E} using \eqref{eq:bound1} and \eqref{eq:bound2} to obtain
\begin{align*}
\Pi_\mE&\leq M_n^{-1}\sum_{\mT\in T}\Pi[\mT\C Y]\int_{\mE(\mT)}   \sup_{x\in[0,1]}\zeta_n(x)\left|f(x)-f_0^{\backslash d}(x)\right|d\,\Pi(f\C Y,\mT) + o(1)\\
&\leq M_n^{-1}\left[\mathcal{O}(\sqrt{\log n})+B_A\right] +o(1).
\end{align*}
The upper bound goes to zero as long as $M_n$ is strictly faster than $\sqrt{\log n}$.

\subsection{Proof of Theorem \ref{thm:uq}}\label{proof:thm_uq}
We follow the strategy of the proof of Theorem 2 in \cite{castillo_rockova}.
We first show the set $\mathcal C_n$ has an optimal diameter, uniformly over the domain $[0,1]$. 

\subsubsection{Optimal Diameter}
We will use the following Lemma (a simple modification of Lemma S-11 in \cite{castillo_rockova}).
\begin{lemma}(Median Tree Estimator)\label{lemma:median_tree}
Consider the prior distribution as in Theorem \ref{thm:one} and let $\mT^*_Y$ be as in \eqref{bulktree}. Then there exists an event $\mA^*_n$ such that $P_{f_0}[\mA^*_n]=1+o(1)$ as $n\rightarrow\infty$ on which the tree $\mT^*_Y$ has the following two properties
\begin{itemize}
\item[(1)] With $S(f_0;A;\upsilon)$ defined in \eqref{eq:setS}, we have
$$
\mT^*_Y\supseteq S(f_0;A;1).
$$
\item[(2)] With $\wt d_l(x)$ as in \eqref{eq:def_d} and with $\wt \mT^*_Y$ denoting the pre-terminal nodes as defined in \eqref{eq:preterminal}
$$
l\leq \min_{x\in I_{lk}}\wt d_l(x)\quad\forall (l,k)\in\wt\mT^*_Y.
$$
\end{itemize}
\end{lemma}
\proof
Recall the notation  $L_n=\lfloor\log_2 n\rfloor$. We denote with 
$$
T_1=\{\mT: \mT_{int}\supseteq S(f_0;A;1)\}\quad\text{and}\quad
T_2=\{\mT: l\leq \min_{x\in I_{lk}}\wt d_l(x)\,\,\forall (l,k)\in\wt\mT_{int}\}
$$ and with $\mA_n$ the event \eqref{eq:event_A}. Then we have, using \eqref{eq:catch_signal}, 
\begin{align*}
P_{f_0}[\mT^*_Y\notin T_1\cap\mA_n]&\leq P_{f_0}[\{\exists (l,k)\in S(f_0;A;1)\cap\mT^{*c}_Y\}\cap \mA_n] \\
&\leq \sum_{l\leq L_n}\sum_{k=0}^{2^l-1}\mathbb{I}[(l,k)\in S(f_0;A;1)]P_{f_0}\left[\Pi\left((l,k)\notin \mT_{int}\C X\right)>0.5\cap \mA_n\right]\\
&\leq E_{f_0} \left(\sum_{l\leq L_n}\sum_{k=0}^{2^l-1}2\,\mathbb{I}[(l,k)\in S(f_0;A;1)] \Pi\left[\left((l,k)\notin \mT_{int}\C X\right)\right]\mathbb I(\mA_n)\right)=o(1).
\end{align*}
Next, we have from \eqref{eq:locally_small}
\begin{align*}
P_{f_0}[\mT^*_Y\notin T_2\cap\mA_n]&\leq P_{f_0}[\{\exists (l,k)\in\wt\mT^*_Y\quad\text{s.t.}\quad l>\min_{x\in I_{lk}}\wt d_l(x)\}\cap \mA_n] \\
&\leq \sum_{l\leq L_n}\sum_{k=0}^{2^l-1}\mathbb{I}[l>\min_{x\in I_{lk}} \wt d_{l}(x)]P_{f_0}\left[\Pi\left((l,k)\in\mT_{int}\C X\right)>0.5\cap \mA_n\right]\\
&\leq E_{f_0} \left(\sum_{l\leq L_n}\sum_{k=0}^{2^l-1}2\, \mathbb{I}[l>\min_{x\in I_{lk}} \wt d_{l}(x)]\Pi\left[\left((l,k)\in \mT_{int}\C X\right)\right]\mathbb I(\mA_n)\right)=o(1).
\end{align*}
Since $P_{f_0}(\mA_n)=1+o(1)$, one obtains $P_{f_0}[\{\mT^*_Y\notin T_1\}\cup\{\mT^*_Y\notin T_2\} ]=o(1)$.\qedhere

\medskip

From Lemma \ref{lemma:median_tree} it follows that, for some suitable increasing sequence $v_n$,
\begin{equation}\label{eq:diameter}
\sup_{f,g\in\mathcal C_n}\left[\sup_{x\in[0,1]} \frac{\zeta_n(x)}{v_n} |f(x)-g(x)|\right]=\mathcal{O}_{P_{f_0}}(1).
\end{equation}
Indeed, for any $f,g\in\cc_n$ we have
\begin{align*}
\sup_{x\in[0,1]} \left[\frac{\zeta_n(x)}{v_n} |f(x)-g(x)|\right]&\leq \sup_{x\in[0,1]} \left[\frac{\zeta_n(x)}{v_n} \left( |f(x)-\wh f_T(x)| +  |\wh f_T(x)-g(x)|\right)\right]\\
&\leq 2 \sup_{x\in[0,1]} \left[\frac{\zeta_n(x)\sigma_n(x)}{v_n}  \right].
\end{align*}
where $\sigma_n(x)$ was defined in \eqref{radiusprox}.
From the properties of the median tree in Lemma \ref{lemma:median_tree}, we know that there exists an event $\mA_n^*$ such that $P_{f_0}(\mA_n^*)=1+o(1)$ 
where the median tree satisfies $2^{l}\leq 2^{\wt d_l(x)}\lesssim (n/\log n)^{1/[2t(x)+1]}$ for all $(l,k_l(x))\in\wt\mT^*_Y$.
For any $x\in[0,1]$ we then have
$$
 \frac{\zeta_n(x)\sigma_n(x)}{v_n}\leq  \left(\frac{n}{\log n}\right)^{\frac{t(x)}{2t(x)+1}-\frac{1}{2}}\sum_{l=0}^{\Lmax} 2^{l/2}
 \mathbb{I}[(l,k_l(x))\in\mT^{*}_{Y\,int}]\lesssim 2^{\wt d_l(x)/2}\left(\frac{n}{\log n}\right)^{\frac{t(x)}{2t(x)+1}-\frac{1}{2}}.
$$
From the definition of $\wt d_l(x)$ we conclude that the right-hand-side is $\mathcal O(1)$ on the event $\mA^*_n$.
This concludes the statement \eqref{eq:diameter}.

\subsubsection{Confidence of the set $\cc_n$}
We first show that the median tree is a (nearly) rate-optimal estimator.
Denote with $\wh f^T_{lk}=\langle \wh f_T,\psi_{lk}\rangle$ and recall $\mathcal S(f_0;A;\upsilon)=\{ (l,k):\ |\beta_{lk}^0|\ge A\log^{1\vee\upsilon}{n}/\sqrt{n} \}$.  
Recall the definition of trees $T$ in \eqref{eq:T}. Let us consider the event
\begin{equation} \label{eventb}
B_n = \{ \mT^\star_Y\in T\} \cap \cA_n,
\end{equation}
where the noise-event $\cA_n$ is defined in \eqref{eq:event_A}. According to Lemma \ref{lemma:median_tree}, we have   $P_{f_0}(B_n)=1+o(1)$. Using similar arguments as around the inequality \eqref{eq:bound1}, on the event $B_n$, we have for some $M>0$
\begin{equation}\label{eq:rate_median_tree}
\sup_{x\in[0,1]}\zeta_n(x)|\wh f_T(x)-f_0(x)|\leq  M\,\sqrt{\log n}.
\end{equation}
Next, one needs to show that $\sigma_n(x)$ is appropriately large for each $x\in[0,1]$. 

Let $\La_n(x)$ be defined by, for $\mu(x)>0$ to be chosen below,
\begin{equation}\label{eq:lambda}
\frac{\mu(x)}{{\log n}}  \left(\frac{n}{\log{n}}\right)^{\frac{1}{2t(x)+1}} \le 2^{\La_n(x)} 
\le \frac{2\,\mu(x)}{{\log n}}\left(\frac{n}{\log n}\right)^{\frac{1}{2t(x)+1}}.
\end{equation}
We will use the following lemma which follows from the proof of Proposition 3 in \cite{hoffman_nickl}
\begin{lemma}\label{lemma:lb_beta}
Assume $f_0\in \mathcal{C}_{SS}(t(x),x,M(x),\eta(x))$. Then  for the sequence $\La_n(x)$ in \eqref{eq:lambda} there exists $C>0$ such that 
$l\geq \La_n(x)$ such that 
\begin{equation}\label{eq:lb_beta}
 |\beta_{lk_l(x)}^0|\geq    C\,2^{-\La_n(x)[t(x)+1/2]}.
\end{equation}
\end{lemma}
\proof
From the definition of local self-similarity, we have for some $c_1>0$
\begin{equation}\label{eq:aux_lb}
2^{-jt(x)}c_1\leq |K_j(f_0)(x)-f_0(x)|\leq \sum_{l\geq j}2^{l/2}|\beta_{lk_l(x)}^0|.
\end{equation}
Now, for all $N \geq 1$ there exists $j\geq \La_n(x)$ such that, using \eqref{eq:aux_lb}
\begin{align*}
|\beta_{jk_j(x)}^0|&\geq\frac{1}{N}\sum_{l=\La_n(x)}^{\La_n(x)+N-1}|\beta_{lk_l(x)}^0|\\
&\geq \frac{2^{-(\La_n(x)+N)/2}}{N(x)} \left(\sum_{l=\La_n(x)}^{\infty}2^{l/2}|\beta_{lk_l(x)}^0| - \sum_{l=\La_n(x)+N}^{\infty}2^{l/2}|\beta_{lk_l(x)}^0|\right)\\
&\geq \frac{2^{-(\La_n(x)+N)/2}}{N(x)} \left(2^{-\La_n(x)t(x)}c_1-c(t(x),N)2^{-(\La_n(x)+N)t(x)}\right)\\
&\geq \frac{2^{-(\La_n(x)+N)/2}}{2N} 2^{-\La_n(x)t(x)}c_1>\underline{c_1}2^{-\La_n(x)[t(x)+1/2]}.\qquad\qedhere
\end{align*}
where $\underline{c_1} = 2^{-N/2}c_1/(2N ) $ and $N$ is large enough. 

\medskip
Combining \eqref{eq:lambda} with \eqref{eq:lb_beta}, one can choose $\mu(x)$ such that for each $x\in[0,1]$ there exists $l\geq \La_n(x)$ such that
$$
 |\beta_{lk_l(x)}^0|> C  [2\mu(x)]^{-t(x)-1/2}\sqrt{\frac{\log{n}}{n}}(\log n)^{t(x)+1/2}\geq A\log n/\sqrt{n}.
$$
Since this is a signal node (i.e. $\beta_{lk_l(x)}^0\in S(f_0;A;1)$), it will be captured by the median tree, i.e. $\wh f^T_{lk_l(x)}\neq 0$. One deduces that the term $(l,k_l(x))$ in the sum defining $\sigma_n(x)$ is nonzero on the event $B_n$, so that
\begin{equation}\label{eq:lb_sigma}
\sigma_n(x) \ge v_n\sqrt{\frac{\log{n}}{n}} |\psi_{lk_l(x)}|
\ge v_n\sqrt{\frac{\log{n}}{n}} 2^{\La_n(x)/2}\geq\frac{\sqrt{ \mu(x)}v_n}{\sqrt{\log n}}\left(\frac{\log n}{n}\right)^{\frac{t(x)}{2t(x)+1}}.
 \end{equation}
For $v_n$ faster than $\log n$ for all $x\in[0,1]$ one has $\sigma_n(x)\geq \sqrt{\log n}/\zeta_n(x)$ and from \eqref{eq:rate_median_tree} one obtains 
 for suitable $v_n$
 \begin{equation}\label{eq:B_n_subset}
 B_n\subset\left\{\sup_{x\in[0,1]}\left[ \frac{1}{\sigma_n(x)}|\wh f_T(x)-f_0(x)|\right] \leq 1/2 \right\}
 \end{equation}
 the desired coverage property, since 
 $$
 P_{f_0}( f_0\in C_n)=P_{f_0}\left(\sup_{x\in[0,1]}\left[ \frac{1}{\sigma_n(x)}|\wh f_T(x)-f_0(x)|\right] \leq 1\right)
 \geq P_{f_0}(B_n)=1+o_{P_{f_0}}(1).
 $$

\subsubsection{Credibility of the set $\cc_n$}
We want to show that
$$
\Pi[\cc_n\C Y]=1+o_{P_{f_0}}(1).
$$
We note that the posterior distribution 
 and the median estimator $\wh f_T$  converge at a rate at $x\in[0,1]$ strictly faster than $\sigma_n(x)$  on the event $B_n$, using again the lower bound on $\sigma_n(x)$ in \eqref{eq:lb_sigma}. In particular, because of
 \eqref{eq:B_n_subset} we can write
\begin{align}
&E_{f_0}\left(\Pi\left[\sup\limits_{x\in[0,1]}\frac{1}{\sigma_n(x)}|f(x)-\wh f_T(x)| \leq 1\,\big |\,  Y\right] \right)\geq \\
&\quad\quad\quad\quad E_{f_0}\left(\Pi\left[\sup\limits_{x\in[0,1]}\frac{1}{\sigma_n(x)}|f(x)- f_0(x)| \leq 1/2\,\big |\,  Y\right] \mathbb{I}_{B_n}\right)+o_{P_{f_0}}(1).
\end{align}
The right side converges to $1$ in  $P_{f_0}$-probability, which concludes the proof of the theorem.

\subsection{Proof of Theorem \ref{th:conc:spikeslab}}\label{sec:proof:th:conc:spikeslab}
The proof follows the lines of \cite{hoffmann}, with several refinements to allow for weaker constraints on the inclusion probabilities $\omega_l$'s. 
For some suitable $B>0$, we define an event (for $L_{max}=\lfloor \log_2 n\rfloor$)
\begin{equation}\label{An}
\mathcal A_{n,B} = \{ | \varepsilon_{lk} | \leq \sqrt{ 2 [\log 2^l+  B\log n]} \quad  \forall (l , k) \quad\text{such that}\quad  l \leq L_{max} \}
\end{equation}
which satisfies 
$P_{f_0}(\mathcal A_{n,B}^c ) \leq \frac{2\log n}{n^B}. $
First, we show an auxiliary Lemma which is reminiscent of Lemma 1 in \cite{hoffmann}. We define 
$S(f_0,A)=\{(l,k): l\leq L_{max}\,\,\text{and}\,\, |\beta_{lk}^0|>A\sqrt{\log n/n}\}.$

\begin{lemma} \label{lem:smallcoef}
Under the assumptions of Theorem \ref{th:conc:spikeslab}  there exists  $a>0$ determined by $\delta >1/4$ defined in Assumption  \ref{cond:omegal}  and $A>a$ such that, uniformly over $\mathcal C(t,M, \eta)$, we have
\begin{equation} \label{eq:smallcoef}
E_{f_0} \Pi \left[ \mathcal T \cap S(f_0, a)^c \neq \emptyset \C Y  \right] = o(1) \quad\text{as } n\rightarrow\infty
\end{equation}
and 
\begin{equation} \label{eq:largecoef}
E_{f_0}\Pi\left[ S(f_0;A)\not\subseteq \mT \C Y\right]=o(1) \quad\text{as } n\rightarrow\infty.
\end{equation}

\end{lemma}

\begin{proof}

We first prove \eqref{eq:smallcoef} folowing  \cite{hoffmann}, except that we have a weaker condition on the prior on $\mT$. Note that 
\begin{equation*}
\begin{split}
\Pi\left[ \mathcal T \cap S(f_0, a)^c \neq \emptyset \C Y\right] &= \Pi\left[ \exists (l,k) \in \mathcal T \cap S(f_0, a)^c  \C Y\right] \\
&\leq 
\sum_{l \leq L_{max}} \sum_{k=0}^{2^l-1}\I\left[(l,k)\notin S(f_0, a)\right] \Pi\left[ (l,k) \in \mathcal T  \C Y\right]. 
\end{split}
\end{equation*}
We denote by $\bT$ the set of all subsets of wavelet coefficients $\beta_{lk}$ up to the maximal depth $L_{max}=\lfloor \log_2 n\rfloor$. 
Then
$$ 
\Pi\left[ (l,k) \in \mathcal T  \C Y\right]  = \frac{ \sum_{\mathcal T \in\bT} \I\left[(l,k) \in \mathcal T \right] \Pi( \mathcal T )  m_n(\mathcal T) }{ \sum_{\mathcal T \in\bT} \Pi( \mathcal T ) m_n(\mathcal T)}
$$
where 
$$
m_n(\mathcal T) =\prod_{(l,k)\notin\mT}\e^{-\frac{n}{2}Y_{lk}^2}\times
\prod_{(l,k)\in\mT}\int \e^{-\frac{n}{2}|Y_{lk}-\beta_{lk}|^2}\pi_{lk}(\beta_{lk})d\beta_{lk}.
$$
For a set $\mT$ such that $(l,k)\in\mT$ we denote with $\mT^{-}=\mT\backslash\{(l,k)\}$.
Due to the fact that the marginal likelihood factorizes, we obtain
$$
\Pi\left[ (l,k) \in \mathcal T  \C Y\right]  = 
 \frac{\sum_{\mathcal T \in\bT} \I\left[(l,k) \in \mathcal T \right] R(\mT,\mT^-){
 \Pi( \mathcal T^- )  m_n(\mathcal T^-)}}{\sum_{\mathcal T \in\bT} \Pi( \mathcal T ) m_n(\mathcal T)},
$$
where, invoking \eqref{cond:priorMT:SS}, we obtain
$$
R(\mT,\mT^-)\coloneqq \frac{\Pi( \mathcal T )  m_n(\mathcal T)}{
 \Pi( \mathcal T^- )  m_n(\mathcal T^-)}\leq \sqrt{2\pi}n^{-1/2}\, C\times C_T\times w_l\times \e^{\frac{n}{2}Y_{lk}^2}.
 $$
This yields
$$
\Pi\left[ (l,k) \in \mathcal T  \C Y\right]\leq  \sqrt{2\pi}n^{-1/2}\, C\times C_T\times w_l\times \e^{\frac{n}{2}Y_{lk}^2}.
$$
On the event $A_{n,B}$ in \eqref{An} we have $|\varepsilon_{lk}|\leq \sqrt{2(1+B)\log n}$ and for $(l,k)\notin S(f_0,a)$ we have $|\beta_{lk}|<a\sqrt{\log n/n}$. We use the fact that for any $b\in(0,1)$ 
$$
\frac{n}{2}Y_{lk}^2\leq\frac{1-b}{2}\varepsilon_{lk}^2+ \frac{b}{2}\,\varepsilon_{lk}^2+|\varepsilon_{lk}|a\sqrt{\log n}+\frac{a^2}{2}\log n
$$ 
to find that for $\wt C\equiv  \sqrt{2\pi}\times C\times C_T$
$$
\Pi\left[ \mathcal T \cap S(f_0, a)^c \neq \emptyset \C Y\right] \leq  \wt C
\, n^{\frac{a^2-1}{2}+b(1+B)} \sum_{l \leq L_{max}}2^l\,w_l \sum_{k: (l,k)\notin S(f_0;a)} \e^{\frac{1-b}{2}\varepsilon_{lk}^2
+\frac{a}{2}|\varepsilon_{lk}|\sqrt{\log n}}.
$$
Since, using \eqref{cond:omegal}, 
\begin{align*}
E_{f_0} \sum_{l \leq L_{max}}2^l\,w_l \sum_{k: (l,k)\notin S(f_0;a)} \e^{\frac{1-b}{2}\varepsilon_{lk}^2
+\frac{a}{2}|\varepsilon_{lk}|\sqrt{\log n}}&\leq 2 \sum_{l \leq L_{max}} 2^l\omega_l \int_{0}^{\infty}
\e^{-\frac{b}{2}x^2+\frac{a}{2}x\sqrt{\log n}}\d x\\
&\leq \frac{2\sqrt{2\pi} n^{a^2/(2b)}}{\sqrt b}    \sum_{l \leq L_{max}} 2^l\omega_l.
\end{align*}
This yields
$$
E_{f_0} \Pi\left[ \mathcal T \cap S(f_0, a)^c \neq \emptyset \C Y\right] \leq 
\frac{2\sqrt{2\pi}   \wt C\, n^{\frac{a^2-1}{2}+b(1+B)+a^2/(2b)}}{\sqrt b}
 \sum_{l \leq L_{max}} 2^l\omega_l.
$$
We can find $b$ and $a$ such that $a^2+2b(1+B)+a^2/b<\delta$ and thereby, using the assumption \eqref{cond:omegal},
$$
E_{f_0} \Pi\left[ \mathcal T \cap S(f_0, a)^c \neq \emptyset \C Y\right] \lesssim L_{max}n^{-c/2}
$$
for $c=\delta-a^2+2b(1+B)+a^2/b>0$. Assuming $\delta>1/4$, we can choose  $a^2<(\delta-1/4)\frac{1}{1+8(1+B)}$.
This proves the first statement \eqref{eq:smallcoef}.
{ We now prove that there exists $A>0$ such that on the event $\mA_{n,B}$ we have \eqref{eq:largecoef}.
We have
$$
\Pi\left[ S(f_0,A)\not\subseteq \mT \C Y\right] \leq\sum_{(l,k)\in S(f_0,A)} \Pi[(l,k)\notin \mT\C Y]
$$
For $\mT$ such that $(l,k)\notin \mT$,  denote with $\mT^+=\mT\cup \{(l,k)\}$. Then
$$
\Pi\left[ (l,k) \notin \mathcal T  \C Y\right]  = 
 \frac{\sum_{\mathcal T \in\bT} \I\left[(l,k) \notin \mathcal T \right] R(\mT,\mT^+){
 \Pi( \mathcal T^+ )  m_n(\mathcal T^+)}}{\sum_{\mathcal T \in\bT} \Pi( \mathcal T ) m_n(\mathcal T)},
$$
where (choosing $R>C_\beta+\sqrt{2(1+B)\log n/n}$ for a suitably large $C_\beta$)
$$
R(\mT,\mT^+)\coloneqq \frac{\Pi( \mathcal T )  m_n(\mathcal T)}{
 \Pi( \mathcal T^+ )  m_n(\mathcal T^+)}\leq \frac{n^{1/2}}{\sqrt{2\pi}\,c_T\,w_l\,c_R\, c} \e^{-\frac{n}{2}Y_{lk}^2}.
 $$
Above, we have used the fact that on the event $\mA_{n,B}$ we have $|Y_{lk}|\leq U\equiv C_\beta+\sqrt{2(1+B)\log n/n}\leq R$ and
for some $c>0$ 
\begin{align*}
 \int \e^{-\frac{n}{2}|Y_{lk}-\beta_{lk}|^2}\pi_{lk}(\beta_{lk})&\geq c_R \int_{-R}^R \e^{-\frac{n}{2}|Y_{lk}-\beta_{lk}|^2}=
 c_R\frac{\sqrt{2\pi}}{\sqrt{n}}[\Phi(R;Y_{lk},1/n)-\Phi(-R;Y_{lk},1/n)]\\
 &\geq c_R\frac{\sqrt{2\pi}}{\sqrt{n}}[\Phi(U;Y_{lk},1/n)-\Phi(-U;Y_{lk},1/n)]\geq  c_R\frac{\sqrt{2\pi}}{\sqrt{n}}(\Phi(2U\sqrt{n};0,1)-1/2)\\
 &\geq c_R\frac{\sqrt{2\pi}}{\sqrt{n}}c
 \end{align*}
 where $\Phi(x;\mu,\sigma)$ is a cdf of a normal distribution with mean $\mu$ and variance $\sigma$.
On the event $\mA_{n,B}$ (in \eqref{An}) we have  $|\varepsilon_{lk}|\leq \sqrt{2(1+B)\log n}$ and 
$$
Y_{lk}^2=[(\beta^0_{lk})^2+\varepsilon_{lk}/\sqrt{n}]^2\geq (\beta^0_{lk})^2/2-4(1+B)\log n/n>[A^2-4(1+B)]\log n/n.
$$
This yields (using the prior assumption \eqref{cond:omegal}
$$
\Pi\left[ \mT: S(f_0,A)\not\subseteq \mT \C Y\right]\leq \frac{ |S(f_0,A)|}{\sqrt{2\pi}\,c_T\,c_R}n^{B_\omega+1/2-[A^2-4(1+B)]/2}<n^{-A^2/4}
$$
for $A^2/4>2(1+B)+B_\omega+1/2$.}\qedhere
\end{proof}

\medskip

Now, we complete the proof of Theorem \ref{th:conc:spikeslab}.
We deploy Lemma \ref{lem:smallcoef} to find $A>a>0$ such that 
for 
$
S(f_0, b)  = \{(l,k):  |\beta^0_{lk}| \geq b \sqrt{\log n/n} \}
$
we obtain, on the event $\mathcal A_{n,B}$, 
$$
\Pi\left[ \mT: \mathcal T \subset S(f_0, a) \C Y\right]\rightarrow 1\quad \text{and} \quad\Pi\left[ \mT: S(f_0,A)\not\subseteq \mT \C Y\right]\rightarrow 0\quad\text{as } n\rightarrow\infty.
$$
Similarly as in the proof of Theorem \ref{thm:one}, we then define
 $$
 T = \{ \mT: S(f_0,A) \subset \mT \subset S(f_0,a)\},
 $$ 
 and we denote with $\mE(\mT)$ the set of of functions $f(x)=\sum_{(l,k)\in\mT}\psi_{lk}(x)\beta_{lk}$. 
We also assume that the bound $A>0$ can be chosen large enough such that 
$$
\sup_{f_0\in \mathcal C(t,M,\eta)} 
E_{f_0}\Pi\left[ \max_{(l,k)\in S(f_0,A)}|\beta_{lk}-\beta^0_{lk}|>A\sqrt{\log n/n}\,\big|\,  Y\right]\lesssim \frac{\log n}{n^B}
$$
This is indeed the case, as shown in the proof of Theorem 3.1 in \cite{hoffmann}.
 {From the definition of local H\"{o}lder balls we have
 $$
\{(l,k): l\leq \min_{x\in I_{lk}} \wt d_l(x,A)\} \subset S(f_0, a)\subset\{(l,k): l\leq \min_{x\in I_{lk}} \wt d_l(x,a)\},
 $$
 where $\wt d_l(x,a)$ is defined as in \eqref{eq:def_d} using $a$ instead of $\bar\gamma$.
We note that for each $f\in\mE(\mT)$ with a coefficient sequence  $\{\beta_{lk}\}$ that satisfies
\begin{equation}\label{eq:concentration_signal}
\max_{(l,k)\in S(f_0,A)}|\beta_{lk}-\beta^0_{lk}|\leq A\sqrt{\log n/n}
\end{equation}
we have for $\zeta_n(x)= \left(\frac{n}{\log n}\right)^{\frac{t(x)}{2t(x)+1}}$
\begin{align}
\zeta_n(x)|f(x)-f_0(x)|\leq& \,\zeta_n(x)|f_0(x)^{\backslash d}|\nonumber \\
&  +\zeta_n(x)\left[\sum_{(l,k)\in\mT} \I({x\in I_{lk}})2^{l/2}  |\beta_{lk}-\beta_{lk}^0|+\sum_{(l,k)\notin \mT} \I({x\in I_{lk}})2^{l/2}|\beta_{lk}^0|\right]\label{eq:second_third_term}
\end{align}
 where   $f_0^d(x)=f_0-f_0^{\backslash d}$ and where  $f_0^{\backslash d}$ is the bias term defined in \eqref{eq:bias} which satisfies $\sup_{x\in[0,1]}\zeta_n(x)|f_0(x)^{\backslash d}|=\mathcal O(1)$ according to Lemma \ref{lemma:bias}.
Focusing on the last term in \eqref{eq:second_third_term}, we know that for each $\mT\in T$, we have 
$
|\beta_{lk}^0|<A\sqrt{\log n/n}$
for $(l,k)\notin\mT$
and $l\leq \min_{x\in I_{lk}}\wt d_{l}(x,A)$. Thereby, we obtain
$$
\sup_{x\in[0,1]}\zeta_n(x)\sum_{(l,k)\notin\mT} \I({x\in I_{lk}})2^{l/2}|\beta_{lk}^0|\lesssim \sqrt{\log n/n}\sup_{x\in[0,1]}\zeta_n(x)2^{\wt d_l(x,A)/2}=\mathcal O(1). 
$$
Regarding the middle term in \eqref{eq:second_third_term}, we use the property \eqref{eq:concentration_signal} and the fact that $(l,k)\in\mT$ and $x\in I_{lk}$ implies  $l\leq \wt d_l(x,a)$. Then
$$
\sup_{x\in[0,1]}\zeta_n(x)\sum_{(l,k)\in\mT} \I({x\in I_{lk}})2^{l/2}  |\beta_{lk}-\beta_{lk}^0|\lesssim \sqrt{\log n/n}\sup_{x\in[0,1]}\zeta_n(x)2^{\wt d_l(x,a)/2}=\mathcal O(1).
$$
This completes the proof of Theorem \ref{th:conc:spikeslab}. $\qedhere$
}

\section{Proof of Theorem \ref{th:GaussianP} }\label{sec:proof_lb}

The proof of Theorem \ref{th:GaussianP} is based on Corollary 2.1 and Theorem 2.2 of \cite{rousseau:szabo:17} for the regression case with wavelet priors and the proof is very similar to the proof of Propositions 3.1 and 3.2 of \cite{rousseau:szabo:17} . 
We first determine $\epsilon_n(\lambda)$ defined by 
\begin{equation}\label{eq:define_epsilon} 
\Pi\left[\|\b-\b_0\|_2 \leq K \epsilon_n(\lambda)\C \lambda \right] = \e^{-n \epsilon_n(\lambda)^2}
\end{equation}
for some $K>0$ and where $\lambda$ is the unknown hyper-parameter $L, \tau$ and  $\alpha$ in cases T1, T2 and T3, respectively. 
We assume that $f_0$ satisfies \eqref{spatialBeta} and determine $\epsilon_n(\lambda)$ and $\epsilon_{n,0}  = \inf_\lambda \epsilon_n(\lambda) $ in all 3 cases (T1)-(T3). 
\smallskip

\noindent \textbf{Case T1.} 
The main difference with Lemma 3.1 of \cite{rousseau:szabo:17} (further referred to as RS17) is that the parameter space is different. 
We denote with $\b_{L} = (\beta_{lk}: l \leq L, k \in I_{l})'$ for $I_l=\{0,1,\dots,2^l-1\}$ where $2^{l} = |I_l|$. Since $\beta_{lk}=0$ for $l>L$, we can write $\|\b-\b_0\|_2^2=\|\b_L-\b_L^0\|_2^2+\sum_{l>L}\sum_{k=0}^{2^l-1}(\beta_{lk}^{0})^2$. For $s_n^2=K^2 \epsilon_n(L)^2-\sum_{l>L}\sum_{k=0}^{2^l-1}(\beta_{lk}^{0})^2$ we use the same arguments as in Lemma 3.1 of RS17 to conclude that
\begin{equation*}
\begin{split}
\Pi\left( \|\b_L - \b_{L}^0 \|_2  \leq s_n \C L \right) \asymp \e^{ 2^L  \log (s_n\,2^{-L/2})(1 + o(1))}
\end{split}
\end{equation*}
 as in the proof of Lemma 3.1 of \cite{rousseau:szabo:17}. We then obtain a variant of equation (A1) in RS17
 $$
 s_n^2+\sum_{l>L}\sum_{k=0}^{2^l-1}(\beta_{lk}^0)^2=\frac{K^22^L}{n}\log\left(\frac{2^{L/2}}{s_n}\right)(1+o(1)).
 $$ 
 In addition, 
$$ 
\sum_{l>L}\sum_{k=0}^{2^l-1} (\beta_{lk}^0)^2\leq  M_1 2^{-\alpha_1 L } (1/2 + o(1))
$$
so that
\begin{equation}\label{T1:epsilon}
\epsilon_{n,0} \lesssim (n/\log n)^{-\alpha_1/(2\alpha_1+1)}\quad\text{and}\quad \epsilon_n(L) \asymp   \sum_{l>L}\sum_k \beta_{lk}^{02}+\frac{  2^L \log n}{n } . 
\end{equation}
For all $\beta_{0}$ satisfying also \eqref{spatialBeta2}, we obtain that 
$$
\epsilon_{n,0} \asymp (n/\log n)^{-\alpha_1/(2\alpha_1+1)}\quad\text{and}\quad
 \epsilon_n(L) \gtrsim  M_1 2^{-\alpha_1 L } +\frac{  2^L \log n}{n } .  
 $$

\noindent \textbf{Case T2 and T3.} 
Similarly as in the proof of Lemma 3.2 in \cite{rousseau:szabo:17}, it can be seen that (see equation (3.4) in RS17 or Theorem 4 of \cite{kuelbs})
$$
-\log  \Pi ( \| \b\|_2 \leq K \epsilon \C \alpha, \tau ) \asymp (K\epsilon/\tau)^{-1/\alpha}.
$$
Note that this equivalence is valid for the non-truncated prior and  remains valid under the priors defined in T2 and T3 for any positive $\alpha$.\footnote{The case $\alpha$ close to 0 is of no importance here since the associated $\epsilon_{n}(\lambda)$ is much bigger than $\epsilon_{n,0}$}
Similarly as in the proof of Lemma 3.2 \cite{rousseau:szabo:17}, we bound from above
$$ \inf_{h \in \mathbb H^{\alpha,\tau}; \| h-\b_0\|_2 \leq \epsilon} \|h\|_{\mathbb H^{\alpha,\tau}}^2 $$
by $\|\be_0\|_{H^{\alpha,\tau}}^2$ if $\alpha_1 >\alpha+1/2$ where, under \eqref{spatialBeta},
$$ \|\be_0\|_{H^{\alpha,\tau}}^2 = \tau^{-2}\sum_{l,k}2^{(2\alpha+1)l}\beta_{lk}^{02} \leq \frac{M_1 }{ \tau^2} \sum_{l,k}2^{(2\alpha+1-2\alpha_1)l} \lesssim   \frac{M_1 }{ \tau^2[2\alpha_1-1-2\alpha]} $$
and by
$$ \|\be_0\|_{H^{\alpha,\tau}}^2 = \tau^{-2}\sum_{l}^{L_\epsilon} \sum_k2^{(2\alpha+1)l}\beta_{lk}^{02} \lesssim   \frac{M_1 2^{(2\alpha+1-2\alpha_1)L_\epsilon} }{ \tau^2[2\alpha+1-2\alpha_1]} , \quad M_1 2^{-2\alpha_1L_\epsilon} =\epsilon^2$$
 if $\alpha_1 <\alpha+1/2$ or by 
 $$ \|\be_0\|_{H^{\alpha,\tau}}^2  \lesssim  \frac{M_1 L_\epsilon }{ \tau^2}, \quad \text{if } \quad  \alpha_1 =\alpha+1/2.$$
 We thus obtain that  if $\alpha_1 \neq \alpha +1/2$ then 
$$\epsilon_n(\alpha, \tau)\lesssim  n^{-\alpha/(2\alpha+1)} \tau^{1/(2\alpha+1) } + \left( \frac{ 1 }{ n\tau^2 (\alpha_1 - \alpha -1/2)}\right)^{\frac{ \alpha_1}{2\alpha+1} \wedge \frac{1}{2}} $$
while if $\alpha_1 = \alpha +1/2$ 
$$\epsilon_n(\alpha, \tau)\lesssim  n^{-\alpha/(2\alpha+1)} \tau^{1/(2\alpha+1) } + \left( \frac{ \log(n\tau^2) }{ n\tau^2}\right)^{ \frac{1}{2}}. $$
Note   that if $f_0$ also  follows \eqref{spatialBeta2}, we can bound from below  (similarly to \cite{rousseau:szabo:17}) for all $h \in \mathbb H^{\alpha,\tau}; \| h-\b_0\|_2 \leq \epsilon_n(\alpha,\tau)$) when $\alpha_1 < \alpha + 1/2$ , $L_n>1$, 
\begin{align*}
  \|h\|_{\mathbb H^{\alpha,\tau}}^2 &\geq \tau^{-2}\sum_{l\leq L_n} \sum_{k \in I_{l1}}2^{(2\alpha+1)l} [\beta_{lk}^{02} -2 |\beta_{lk}^0||\beta_{lk}^0-h_{lk} | ]\\
  &\geq \frac{ m_1 }{ 2\tau^2(2\alpha+1-2\alpha_1) }  2^{(2\alpha+1-2\alpha_1)L_n} -
  \frac{ 2 M_1 }{\tau^2} 2^{(2\alpha-\alpha_1+1/2)L_n}\sum_{l\leq L_n} \sum_{k \in I_{l1}}|\beta_{lk}^0-h_{lk} |\\
  &\geq \frac{ m_1 }{ 2\tau^2(2\alpha+1-2\alpha_1) }  2^{(2\alpha+1-2\alpha_1)L_n} - C2^{(2\alpha-\alpha_1+1)L_n}\epsilon_n(\alpha, \tau) \\
  &= \frac{ m_1 }{ 2\tau^2(2\alpha+1-2\alpha_1) }  2^{(2\alpha+1-2\alpha_1)L_n} \left( 1 - \frac{ C 2^{\alpha_1L_n}}{m_1}\epsilon_n(\alpha, \tau)\right)
\end{align*}
for some $C>0$ and choosing $L_n$ equal to $L_n = \left\lfloor \frac{ \log \left(\frac{m_1}{2C\epsilon_n(\alpha, \tau)} \right) }{ \alpha_1 \log 2 } \right\rfloor $, 
we bound 
$$
  \|h\|_{\mathbb H^{\alpha,\tau}}^2 \gtrsim \epsilon_n(\alpha,\tau)^{ -  (2\alpha+1-2\alpha_1)/\alpha_1} $$
  which leads to 
$$\epsilon_n(\alpha, \tau)\gtrsim  n^{-\alpha/(2\alpha+1)} \tau^{1/(2\alpha+1) } + \left( \frac{ 1 }{ n\tau^2 (\alpha_1 - \alpha -1/2)}\right)^{\frac{ \alpha_1}{2\alpha+1} } $$
while if $\alpha_1 = \alpha +1/2$ 
$$\epsilon_n(\alpha, \tau)\gtrsim  n^{-\alpha/(2\alpha+1)} \tau^{1/(2\alpha+1) } + \left( \frac{ \log(n\tau^2) }{ n\tau^2}\right)^{\frac{ \alpha_1}{2\alpha+1} \wedge \frac{1}{2}} $$
and if $ \alpha_1 <\alpha +1/2$, since $\| \b_0\|_2 \geq c>0$ for some fixed $c$, 
$$\epsilon_n(\alpha, \tau)\gtrsim  n^{-\alpha/(2\alpha+1)} \tau^{1/(2\alpha+1) } + \left( \frac{ \|\b_0\|_2 }{ n\tau^2 (\alpha_1 - \alpha -1/2)}\right)^{1/2}. $$
The lower bounds thus match the previous upper bounds.  Minimizing in $\alpha$ in the case T3 these upper and lower bounds lead to  choosing $\alpha = \alpha_1$ and 
$\epsilon_{n,0} \asymp n^{-\alpha_1/(2\alpha_1+1)}$ while minimizing in $\tau$ (Case T2) 
 with $\alpha_1 < \alpha+1/2$, the minimum is obtained by considering $\tau \asymp n^{-(\alpha_1-\alpha)/(2\alpha_1+1)}$ and 
 $\epsilon_{n,0} \asymp n^{-\alpha_1/(2\alpha_1+1)}$ while if $\alpha_1 \geq \alpha+1/2$ it is obtained with  $\tau \asymp n^{-1/(4\alpha+4)}$ leading to $\epsilon_{n,0} \asymp n^{-(2\alpha+1)/(4\alpha+4)}$ up to a $\log n $ term. 
\smallskip

Having quantified $\epsilon_{n,0}$ for the three cases,   the upper bound \eqref{eq:statement1} follows directly from Theorem 2.3 of RS17.
\medskip

Regarding the lower bound, we let $\Lambda_0 = \{ \lambda : \epsilon_n(\lambda) \leq M_n \epsilon_{n,0}\}$ with $M_n $ going to infinity and $\lambda$ either  $L,\tau$ or $\alpha$ in cases T1, T2 and T3, respectively. 
Then following from the proofs of {Propositions 3.1 and 3.2 of  \cite{rousseau:szabo:17} and since the priors on $\lambda$ satisfy condition H1 (using Lemma 3.5 and 3.6 of \cite{rousseau:szabo:17})} one obtains
\begin{equation*}
\Pi( \lambda \in \Lambda_0 \C Y) = 1 + o_{P_{f_0}}(1).
\end{equation*} 
From this and the remark that for all $\b=\{\beta_{lk}\}$ such that $\beta_{lk}=0$ for $l>{L_{max}}=\lfloor\log_2(n/2)\rfloor$ we have for some $C_1>0$
\begin{equation*} 
\begin{split}
\|f_{\beta_0} - f_\beta\|_n &\leq \|f_{\beta_0,L_{max}} - f_{\beta}\|_n+ C_1\,M_1n^{-\alpha_1} = \|f_{\beta_0,L_{max}} - f_{\beta}\|_2+C_1\, M_1n^{-\alpha_1} \\
&\geq \|f_{\beta_0,L_{max}} - f_{\beta}\|_n- C_1\,M_1n^{-\alpha_1} = \|f_{\beta_0,L_{max}} - f_{\beta}\|_2-C_1\, M_1n^{-\alpha_1}. 
\end{split}
\end{equation*}
Together with the fact $\|f_{\beta_0,L_{max}} - f_{\beta}\|_2 = \|\b_{L_{max}}^0 -\b\|_2$ we obtain for  some $\delta>0$
$$
B_n = \left\{ f: \|f-f_0\|_{1/2,1} \leq n^{-\delta} \epsilon_n(\alpha_1) \right\}
$$  
the following (with $l_n(\b)=\log \Pi(Y\C\b)$ and $m_n(\lambda)=\int_{\b}\e^{ l_n(\beta) - l_n(\beta_0)}d\,\Pi(\b)$)
\begin{equation*}
\begin{split}
 \Pi( B_n \C Y) 
 & = \Pi\left(B_n \cap \{ \lambda \in \Lambda_0 \} \C Y\right) +o_{P_{f_0}}(1)\\
 & \leq \frac{ \int_{\lambda \in \Lambda_0 } \int_{B_n} \e^{ l_n(\beta) - l_n(\beta_0)}d\,\Pi(\b\C\lambda) \,d\Pi(\lambda)}{  \int_{\lambda \in \Lambda_0 } m_n(\lambda) d\,\Pi(\lambda)}.
 \end{split}
\end{equation*} 

\noindent 
\textbf{The Case T1.} 
We have $\lambda = L$ and set $L_{n,1} $ such that $L_{n,1} = \lfloor \log (L_0 (n/\log n)^{1/(2\alpha_1+1)})/ \log 2\rfloor$ for some suitable $L_0>0$. Then  
$ \Pi ( L \in \Lambda_0 \C Y) = 1 + o_{P_{f_0}}(1)$ and 
and for all $L \in \Lambda_0$ under \eqref{spatialBeta}  and \eqref{spatialBeta2} we have
$$
2^{L_{n,1}}M_n^{-2/\beta_1} \lesssim  2^{L} \lesssim 2^{L_{n,1}} M_n^2.
$$  
Moreover 
\begin{equation*}
\begin{split}
\Pi(B_n \C Y) = \sum_{L\in \Lambda_0 } \Pi( B_n \C Y, L ) \Pi(L\C Y) + o_{P_{f_0}}(1).
 \end{split}
\end{equation*}

 It will be useful to rewrite \eqref{eq:f0} in a vector notation. For any $L\geq 1$, 
 we denote with  $\b_L$  a vector with coordinates $\beta_j = \beta_{lk}$ for $j = 2^{l-1} + k-1$. 
{If $L \leq L_{max}$ and $\be $ is according to the  \textit{model $L$}, i.e. $\beta_{lk} = 0$ for all $l >L$, the log-likelihood at $\be$ (conditionally on the model $L$) can be written as 
$\ell_n(\be) = \ell_n(\hat \be_L) - \frac{ (\be - \hat \be_L)^t \Psi_L^t \Psi_L (\be - \hat \be_L) }{ 2 \sigma^2} + C$, where 
$$
\Psi_L(i,j)  =\psi_{\ell k} (x_i)  \quad \text{for}\quad j = 2^{\ell-1}+k-1, \quad  i\leq n, \quad \text{and for}\quad \hat \be_L   =\frac{  \Psi_L^t Y}{n},
$$ 
since  $\Psi_L^t \Psi_L =n I$ so that $\ell_n(\be) =  - \frac{n \|\be - \hat \be_L\|^2}{ 2\sigma^2 } + C'$.}
This implies, together with the Gaussian prior on the $\beta_{lk}$ when $l \leq L$, that for all $L\leq \Lmax$ the conditional posterior given $L$ is Gaussian with a mean 
$ \tilde \b_L = \frac{ n \hat \b_L}{ 1 + n}$ and a variance $ \sigma^2 I/(n+1)$. In the following, we  write $\tilde \b_L$  as the subvector of $\tilde \b_L$  whose coordinates correspond to $j = 2^{l-1}+k-1$ with $k \in I_{l2}$ and $\mathcal N(\tilde \theta_{L,2}, \sigma^2I_2/(n+1))$ as a Gaussian vector with a mean $\tilde \theta_{L,2}$ and a covariance matrix $\sigma^2I_2/(n+1)$, where $I_2$ is the identity matrix of dimension $|I_2| = \sum_{l\leq L} I_{l2} \asymp c 2^L $ for some $c>0$  and  $L\in \Lambda_0$. We have
\begin{equation*}
\begin{split}
 \Pi ( B_n \C Y, L ) &\leq P( \|\mathcal N(\tilde \theta_{L,2}, \sigma^2I_2/(n+1) )\|\leq  \delta \epsilon_{n,0} ) \\
 &\leq  P( \|\mathcal N(0,I_2)\|^2 \leq (n+1)\frac{ n^{-2\delta} \epsilon_{n,0}^2}{ \sigma^2 }  )\\
 &=P\left( \mathcal X^2(|I_2|) \leq  (n+1)\frac{ n^{-2\delta} \epsilon_{n,0}^2}{ \sigma^2 }\right) \lesssim \e^{ - c' |I_2|}
\end{split}
\end{equation*}
 for some $c'>0$, since 
$$  (n+1)\frac{ n^{-2\delta} \epsilon_{n,0}^2}{ \sigma^2 } \lesssim n^{-2\delta} (\log n)^q n^{1/(2\alpha_1+1)} \lesssim n^{-\delta} |I_2| \quad \text{ for some } q>0$$
and
$$ \Pi ( B_n \C Y) = o_p(1) .$$ 
Note that the same holds true if the prior is not Gaussian and if $\alpha_1>1/2$.

\section{Intermediate Results  for Theorem \ref{thm:np}}\label{sec:proof_thm_np}
We first   describes some properties of the Gram matrix induced by irregular designs.
Note that Lemma \ref{lemma:cor} implies that, under the balancing  Assumption \ref{ass:design}, we have for the $j^{th}$ column $X_j$ of $X$ with $j=2^l+k$
\begin{equation}\label{eq:design_norm}
\|X_j\|_2^2=2^{l}n_{lk}\leq 2{n\, (C+l)}\quad\text{and}\quad \|X_j\|_1=2^{l/2}n_{lk}\leq \frac{2n\, (C+l)}{2^{l/2}}.
\end{equation}
and for $i=2^{l_2}+k_2$
\begin{equation}\label{eq:covariance}
|X_j'X_i| \leq  C_d {\sqrt{n}\log^\upsilon n}\,2^{\frac{l}{2}}\mathbb I\{(l_2,k_2)\,\text{is a descendant of}\, (l,k)\}.
\end{equation}

\medskip

Recall the notation of pre-terminal nodes $\wt\mT_{int}$ in \eqref{eq:preterminal} and let $\mX=\{x_i:1\leq i\leq n\}$. We will also be denoting with $\lambda_{min}(A)$ and $\lambda_{max}(A)$ the minimal and maximal eigenvalues of a matrix $A$. 
The idea behind the proof is similar to the one of Theorem \ref{thm:one}. We will be using the same definition of $\wt d_{l}(x)$ in \eqref{eq:lowerb},   $T$ in \eqref{eq:T},
$S(f_0;A;\upsilon)$ in \eqref{eq:setS} and $\mE$ in \eqref{eq:E}.
First, we show that $E_{f_0}\Pi(\mE^c\C Y)\rightarrow 0$. 

\smallskip

To this end, in Section \ref{sec:small_trees2} we show that the posterior concentrates on locally small trees  and in Section \ref{sec:signal_catch} we show that the posterior trees catch signal nodes. These results will be conditional on the set $\mA$ in \eqref{eq:setA_np}.
The complement of this set has a vanishing probability $P_{f_0}(\mA^c)\leq 2/ p\rightarrow 0$.

\subsubsection{Posterior Concentrates on Locally Small Trees}\label{sec:small_trees2}
 We now show that 
\begin{equation}
\Pi\left[\mT: \exists (l,k)\in \wt\mT_{int} \,\,\,\, s.t.\,\,\,\, l>\min_{x\in I_{lk}\cap \mX}\wt d_{l}(x)\C Y\right]\rightarrow 0.\label{eq:locally_small_np}
\end{equation}
on the set
\begin{equation}\label{eq:setA_np}
\mA=\{\bm\varepsilon: \|X'\bm\varepsilon\|_\infty\leq 2\|X\|\sqrt{\log p}\},
\end{equation}
where $\|X\|=\max\limits_{1\leq j\leq p}\|X_j\|_2$. 

\smallskip
To prove this statement, we follow the route of Lemma \ref{lemma:locally_small} for the white noise model. 
The irregular design requires non-trivial modifications of the proof due to the induced correlation among predictors. Similarly as in the proof of Lemma \ref{lemma:locally_small}, we
denote with $\mT^-$ the  sub-tree of $\mT$  obtained by deleting a deep node $(l_1,k_1)$ which corresponds to the column $X_j$ where $j=2^{l_1}+k_1$ and which satisfies  $l_1\geq \wt d_l(x)$ (as in \eqref{eq:lowerb}) for some $x\in I_{l_1k_1}$ and thereby $|\beta_{l_1k_1}|\lesssim \sqrt{\log n/n}$.
Then we have
\begin{equation}\label{eq:ratio_minus}
\frac{N_\mT(Y)}{N_{\mT^-}(Y)}=\frac{1}{\sqrt{1+g_n}}\exp\left\{ \frac{1}{2} Y'[X_\mT\Sigma_\mT X_\mT' -X_{\mT^-}\Sigma_{\mT^-} X_{\mT^-}' ]Y\right\}.
\end{equation}
Using Lemma \ref{lemma:dif}, we  simplify the exponent in \eqref{eq:ratio_minus} to find for $c_n=g_n/(g_n+1)$
$$
\frac{N_\mT(Y)}{N_{\mT^-}(Y)}=\frac{1}{\sqrt{1+g_n}} \exp\left\{\frac{c_n|X_j'(I-P_{\mT^-})Y|^2}{2Z}\right\}.
$$
First, we bound the term 
 \begin{align}
|X_j'(I-P_{\mT^-})Y|^2&= |X_j'(I-P_{\mT^-})(X_j\beta_j^0+X_{\bmT}\b_{\bmT}^0+\bm \nu)|^2\nonumber\\
&\leq 2  |X_j'(I-P_{\mT^-})X_j|^2|\beta_j^0|^2+
2|X_j'(1-P_{\mT^-})(X_{\bmT}\b_{\bmT}^0+\bm \nu)|^2.\label{eq:two_terms}
 \end{align}
Using the design assumption \eqref{eq:design_norm},  the first term satisfies (since $\lambda_{max}(I-P_{\mT^-})=1$)  
 $$
 \frac { |X_j'(I-P_{\mT^-})X_j|^2|\beta_j^0|^2}{Z}=Z|\beta_j^0|^2\leq \|X_j\|_2^2\log n/n\lesssim \log^2 n,
 $$
 where we used the fact that $(l_1,k_1)$ is deep and thereby $|\beta_j^0|\lesssim \sqrt{\log n/n}$.
Using   \eqref{eq:covariance}, the H\"{o}lder condition  and  the  assumption $t_1> 1/2$,
we obtain   (since $\lambda_{max}(I-P_{\mT^-})=1$) 
 \begin{align*}
 |X_j'(I-P_{\mT^-})X_{\bmT}\b_{\bmT}^0|&\leq C_d\sqrt{n}\log^\upsilon n
\sum_{(l_2,k_2)}\mathbb I[(l_2,k_2)\,\,\text{is a descendant of $(l_1,k_1)$}] 2^{l_1/2}2^{-l_2(t_1+1/2)}\\
&\leq   C_d\sqrt{n}\log^\upsilon n\sum_{l_2=l_1+1}^{L_{max}} 2^{l_2-l_1}\, 2^{l_1/2}\,2^{-l_2(t_1+1/2)} \\
&\lesssim   \frac{C_d\sqrt{n}\log^\upsilon n}{2^{l_1/2}} \lesssim \sqrt{n}\log^{\upsilon} n.
 \end{align*}
Regarding the second term in \eqref{eq:two_terms}, on the event $\mA$, we have  from the  Lemma \ref{lemma:bias2}
 \begin{equation}\label{eq:X_bias}
 |X_j'(I-P_{\mT^-})\bm\nu|\leq |X_j'(F_0-X\b_0+\bm \varepsilon)|\lesssim \sqrt{n}\log^{1\vee\upsilon}n.
 \end{equation}

Now, we find a lower bound for $Z=X_j'(I-P_{\mT^-})X_j$. 
From the proof of Lemma \ref{lemma:dif}, we can see that $1/Z$ is a `submatrix' of $(X_\mT'X_\mT)^{-1}$.
The eigenvalue of this `submatrix' will be smaller than the maximal eigenvalue of the entire matrix $(X_\mT'X_\mT)^{-1}$ (from the interlacing eigenvalue theorem \citep{matrix_book})
and thereby
$$
1/Z\leq \lambda_{\max} (X_\mT'X_\mT)^{-1}=1/\lambda_{min}(X_\mT'X_\mT).
$$
From Lemma \ref{lemma:eigen} we have  
\begin{equation}\label{eq:min_eigen}
\lambda_{min}(X_\mT'X_\mT)\geq \munderbar{\lambda}\, n\,\quad \text{for some}\quad \munderbar{\lambda}>0.
\end{equation}
From $Z\geq \munderbar{\lambda}\,n$  we then obtain for some suitable $C>0$
 $$
 \frac{N_\mT(Y)}{N_{\mT^-}(Y)}\leq \exp\left(C \log^{2(1\vee \upsilon)} n\right).
 $$
 We can now continue as in the proof of Theorem \ref{thm:one} by plugging-in the likelihood ratio above in the expression \eqref{eq:aux33}. Earlier in the proof of Lemma \ref{lemma:locally_small}, the likelihood ratio was of the order $\e^{C\log n}$. Here, we have a larger logarithmic factor which can be taken care off by choosing $p_l=(\Gamma)^{-l^{2(1\vee\upsilon)}}$ as the split probability. We then conclude  \eqref{eq:locally_small_np} using the same strategy as in the proof of Lemma 
 \ref{lemma:locally_small} for the white noise.

\subsubsection{Catching Signal}\label{sec:signal_catch}
We now show  that, on the event $\mA$ in \eqref{eq:setA_np},  
\begin{equation}\label{eq:catch_signal_np}
\Pi\left[ \mT: S(f_0;A;\upsilon)\not\subseteq \mT \C Y\right]\rightarrow 0\quad\text{as $n\rightarrow\infty$},
\end{equation}
where $\upsilon$ is the balancing constant in the design Assumption \ref{ass:design}.

\smallskip

The proof of \eqref{eq:catch_signal_np} follows the route of Lemma 3 in \cite{castillo_rockova} with nontrivial alterations due to the fact that we now have the regression model where the regression matrix is not orthogonal. Suppose that 
$ (l_1,k_1)\in S(f_0;A;\upsilon) $
is a signal node for some $A>0$ and let $\mT$ be such that $(l_1,k_1)\notin\mT$. 
We grow a branch from $\mT$ that extends towards $(l_1,k_1)$ to obtain an enlarged tree $\mT^+\supset\mT$. In other words $\mT^+$ is the smallest tree that contains $\mT$ and $(l_1,k_1)$ as an internal node. For details, we refer to Lemma 3 in \cite{castillo_rockova}.
 We define $K=|\mT^+_{int}\backslash \mT_{int}|$ and write
\begin{equation}\label{eq:ratio2}
\frac{N_\mT(Y)}{N_{\mT^+}(Y)}=(1+g_n)^{K/{2}}\exp\left\{ \frac{1}{2} Y'[X_\mT\Sigma_\mT X_\mT' -X_{\mT^+}\Sigma_{\mT^+} X_{\mT^+}' ]Y\right\}.
\end{equation}
We denote with $\mT=\mT^-\rightarrow\mT^1\rightarrow\dots\rightarrow\mT^K=\mT^+$ the sequence of nested trees obtained by adding one additional internal node towards $(l_1,k_1)$.
Then using Lemma \ref{lemma:dif} we find
\begin{align}\label{eq:ratio2}
\frac{N_\mT(Y)}{N_{\mT^+}(Y)}&=(1+g_n)^{K/2}\prod_{j=1}^K\exp\left\{ \frac{c_nY'(P_{\mT^{j-1}}-P_{\mT^j})Y}{Z_j} \right\}\\
&=(1+g_n)^{K/2}\prod_{j=1}^K\exp\left\{- \frac{c_n |X_{[j]}'(I-P_{j-1})Y|^2}{Z_j} \right\},
\end{align}
where 
$$
P_j=X_{\mT^j}(X_{\mT^j}'X_{\mT^j})^{-1}X_{\mT^j}'\quad\text{and}\quad Z_j=X_{[j]}'(I-P_{j-1})X_{[j]}
$$
and where $X_{[j]}$ is the column added at the $j^{th}$ step of  branch growing.
Let $X_{[K]}$ be the {\em last} column to be added to $X_{\mT^+}$, i.e. the {\em signal} column associated with $(l_1,k_1)$. We will be denoting simply $\beta_{[K]}^0\equiv \beta^0_{l_1k_1}$ the coefficient associated with $X_{[K]}$.
Then (using the fact that $P_{K-1}$ is a projection matrix onto the columns of $X_{\mT^{K-1}}$)
$$
|X_{[K]}'(I-P_{K-1})Y|^2=|X_{[K]}'(I-P_{K-1})X_{[K]}\beta_{[K]}^0+ X_{[K]}'(I-P_{K-1})X_{\bmT^K}\b_{\bmT^K}^0+X_{[K]}'(I-P_{K-1})\bm \nu|^2
$$
Using the inequality $(a+b)^2\geq a^2/2-b^2$, we find that
$$
\frac{|X_{[K]}'(I-P_{K-1})Y|^2}{Z_K}\geq \frac{Z_K|\beta_{[K]}^0|^2}{2}- \frac{1}{Z_K}|X_{[K]}'(I-P_{K-1})X_{\bmT^K}\b_{\bmT^K}^0+X_{[K]}'(I-P_{K-1})\bm \nu|^2.
$$
Next, since all entries in $X_{\bmT^K}$ are either descendants of $(l_1,k_1)$ or are orthogonal to $X_{[K]}$ we have  (using similar arguments as before in Section \ref{sec:small_trees2})
\begin{align*}
|X_{[K]}'(I-P_{K-1})X_{\bmT^K}\b_{\bmT^K}^0|&\leq |X_{[K]}'X_{\bmT^K}\b_{\bmT^K}^0|\lesssim \sqrt{n}\log^{\upsilon}n.
\end{align*}
Using Lemma \ref{lemma:bias2}, we find that $|X_{[K]}'(I-P_{K-1})\bm \nu|\leq \sqrt{n}\log^{1\vee\upsilon} n$  which yields 
$$
\frac{|X_{[K]}'(I-P_{K-1})Y|^2}{Z_K}\geq \frac{Z_K|\beta_{[K]}^0|^2}{2}- \frac{C_1n\log^{2(1\vee \upsilon)} n}{Z_K}\quad\text{for some $C_1>0$}.
$$
The term $Z_K=X_{[K]}'(I-P_{K-1})X_{[K]}$ is a submatrix of the matrix $(X_{\mT^K}'X_{\mT^K})^{-1}$ and by our assumption \eqref{eq:min_eigen}
we have $Z_K\geq \munderbar\lambda n$ which yields (for $g_n=n$) and from the assumption $|\beta_{[K]}^0|>A\log^{1\vee\upsilon} n/\sqrt{n}$ for some sufficiently large $A>0$ and some $C_2,C_3>0$
$$
\frac{N_\mT(Y)}{N_{\mT^+}(Y)}\leq \e^{K\log (1+g_n)- C_2\log^{2(1\vee\upsilon)} n}=\exp\left\{-C_3\log^{2(1\vee \upsilon)} n\right\}.
$$
As was shown in the proof of Lemma 3 in \cite{castillo_rockova}, we have $\Pi(\mT)/\Pi(\mT^+)\leq 4\Gamma^{2l_1^2}$, and thereby for some $C_4>0$
$$
\Pi[(l_1,k_1)\notin \mT_{int}\C Y]\leq \exp\left\{C_4(\log\Gamma)\log^2n-C_3\log^{2(1\vee \upsilon)}n\right\}.
$$
 Thereby
$$
\sum_{(l_1,k_1)\in S(f_0;A;\upsilon)}\Pi[(l_1,k_1)\notin \mT_{int}\C Y]\leq \e^{-C_3\log^2n}2^{L_{max}+1}\lesssim  \e^{-C_3/2\log^2n} \rightarrow 0.
$$
This concludes the proof of \eqref{eq:catch_signal_np}.

\section{Proof of Lemma  \ref{lem:small:nx} in Section \ref{sec:proof_repulsive}} \label{pr:small:nx}

To prove Lemma \ref{lem:small:nx}, we split $ \left\{  n_{I_x} \leq s_n(\delta)\right\} $ into $B_{n,1} = \left\{  n_{I_x} \leq s_n(\delta)\right\}  \cap \{ | \bar y_{I_{x,1}} -f_0(x) | \leq M_0 \varepsilon_n\} $,  $B_{n,2} = \left\{  n_{I_x} \leq s_n(\delta)\right\} \cap \{ | \bar y_{I_{x,1}} -f_0(x) | > M_0 \varepsilon_n\} \cap \{ n_{I_{x,1} } \leq s_n(\delta_1)\} $ and $B_{n,3} = \left\{  n_{I_x} \leq s_n(\delta)\right\}  \cap \{ | \bar y_{I_{x,1}} -f_0(x) | > M_0 \varepsilon_n\} \cap \{ n_{I_{x,1} } > s_n(\delta_1)\} $ where $\sqrt{\delta_1 } M_0 > 2 u_0$. 

We first consider $B_{n,1}$. 
We have for $\bar I = I_x \cup I_{x,1}$   and writing $S=S'\cup I_x\cup I_{x,1}$
\begin{align*}
\Pi( B_{n,1} | D_n ) &= \frac{ \sum\limits_{S=S'\cup I_x\cup I_{x,1} } m(S') m( \bar I) \pi_S( S'\cup \bar I )  \1_{B_{n,1}}\frac{ m(I_x)m(I_{x,1} )\pi_S( S'\cup I_x \cup I_{x,1} ) }{ m( \bar I)\pi_S( S'\cup \bar  I) } }{ \sum_{S } m(S)\pi_S(S)}.
\end{align*}
Moreover on  $\Omega_{n,x}(u_1) \cap \Omega_{n,y}(u_0)$, 
 \begin{align*}
 \frac{ m(I_x)m(I_{x,1} ) }{ m( \bar I) } &\leq \frac{ c_1^2\sqrt{2\pi} \sqrt{ n_{\bar I} } }{c_0 \sqrt{ n_{ I_x} }\sqrt{ n_{ I_x,1}} }\exp\left\{ \frac{ n_x}{ 2 }(\bar y_{I_x}-\bar y_{\bar I})^2+ \frac{ n_{I_{x,1}}}{ 2 }(\bar y_{I_{x,1}}-\bar y_{\bar I})^2\right\}\\
 &= \frac{ c_1^2\sqrt{2\pi} \sqrt{ n_{\bar I} } }{c_0 \sqrt{ n_{ I_x} }\sqrt{ n_{ I_x,1}} }\exp\left\{ \frac{ n_{I_{x,1}} n_x}{ 2n_{\bar I}  }(\bar y_{I_x}-\bar y_{I_{x,1}})^2\right\}
 \end{align*}
 Note that on  $\Omega_{n,x}(u_1)$ we have
 $$
 p_0|I_x|\leq n_{I_x}/n+u_1\sqrt{p_1\log n/n}  \lesssim s_n(\delta)/n = \delta \varepsilon_n^{1/t(x)}
 $$
and that  $  n_{I_{x,1}} n_{I_x}/ n_{\bar I} \leq n_{I_x} \leq s_n(\delta)$. Moreover, 
 \begin{equation}\label{eq:decompose} 
 \bar y_{I_x}-\bar y_{I_{x,1}} = \bar \epsilon_{I_x} + \bar \beta_{0, I_x}-f_0(x)+ f_0(x) -\bar y_{I_{x,1}} 
 \end{equation}
 and  $ |\bar \beta_{0, I_x}-f_0(x)| \leq M |I_x|^{t(x)}  \leq \delta^{t(x) } C \varepsilon_n$ for some $C$ independent on $\delta$ and $n$. 
Therefore when  $  | \bar y_{I_{x,1}} -f_0(x) | \leq M_0 \varepsilon_n $
 \begin{align*}
 \frac{ n_{I_{x,1}} n_{I_x}}{ 2n_{\bar I}  }(\bar y_{I_x}-\bar y_{I_{x,1}})^2 &\leq \frac{  n_{I_x} \bar \epsilon_{I_x}^2 }{2 } + C \delta^{2t(x)+1} \log n  + \frac{\delta M_0^2 \log n  }{ 2} +  \sqrt{n_x}|\bar \epsilon_{I_x}| \sqrt{\delta } [ M_0 +C \delta^{t(x)+1/2} ] \sqrt{\log n}\\
 & \leq  \frac{  n_{I_x} \bar \epsilon_{I_x}^2 }{2 } + \delta \log n [ C \delta^{2t(x) } + M_0^2/2 + (M_0 +C \delta^{t(x)}) u_0] \leq \frac{  n_{I_x} \bar \epsilon_{I_x}^2 }{2 }  + \delta M_0^2 \log n
 \end{align*}
 on $\Omega_n$, as soon as $M_0$ is large enough (independently of $\delta$) .
  
Moreover, on $\Omega_n$ we can also bound $n_{I_x} \bar \epsilon_{I_x}^2 $ by $u_0^2 \log n$ so that for all $b\in (0,1)$, 
so that 
\begin{align*}
 \frac{ m(I_x)m(I_{x,1} ) }{ m( \bar I) } 
 &\leq \frac{ c_1^2\sqrt{2\pi} \sqrt{ n_{\bar I}}}{c_0 \sqrt{ n_{ I_x} }\sqrt{ n_{ I_{x,1}}}}n^{\delta M_0^2} e^{\frac{  n_{I_x} \bar \epsilon_{I_x}^2 }{2 } } \\
 & \leq \frac{ c_1^2\sqrt{2\pi} \sqrt{ n_{\bar I}}}{c_0 \sqrt{ n_{ I_x} }\sqrt{ n_{ I_{x,1}}}}n^{[\delta M_0^2+ bu_0^2/2]} e^{\frac{ (1-b) n_{I_x} \bar \epsilon_{I_x}^2 }{2 } } 
 \end{align*}
and denoting $Z_b^{I_x}\equiv \exp\left\{ \frac{  (1-b) n_{I_x} \bar \epsilon_{I_x}^2 }{2 } \right\}$
and using the fact that 
$$ E( Z_b^{I_x} |X) = \int \frac{ e^{(1-b) u^2/2 - u^2/2} }{\sqrt{2\pi}} du = 1/\sqrt{b} <\infty $$
we obtain on $\Omega_n$,  
\begin{align*}
\Pi( B_{n,1} | D_n ) &\leq \frac{ \sqrt{2\pi}n^{ bu_0^2/2 + \delta M_0^2 } c_1^2}{ c_0} \\
   & \times \quad \frac{ \sum\limits_{S=S'\cup\bar I } m(S') m( \bar I) \pi_S( S'\cup \bar I ) \sum\limits_{\bar I = I_x\cup I_{x,1}} \1_{B_{n,1}}\frac{ |I_{x,1}|^B| I_x|^B\sqrt{ n_{\bar I} } }{ |\bar I|^B\sqrt{ n_{ I_x} }\sqrt{ n_{ I_x,1}}}Z_b^{I_x} }{
\sum\limits_{S=S'\cup \bar I } m(S') m( \bar I) \pi_S( S'\cup \bar I )}.
\end{align*}
Note that for any $\bar I$ containing $x$, there are many possible choices for $(I_x, I_{x,1}) $ such that $\bar I = I_{x,1} \cup I_x$. 
Also  $n_{\bar I} /[n_{I_{x,1}}n_{ I_x} ] \leq 2/(n_{I_{x,1}} \wedge n_{ I_x})$ so that, choosing without loss of generality $n_{ I_x} \leq n_{I_{x,1}}$, 
$$
\frac{ |I_{x,1}|^B| I_x|^B\sqrt{ n_{\bar I} } }{ |\bar I|^B\sqrt{ n_{ I_x} }\sqrt{ n_{ I_x,1}}  } \leq \frac{ \sqrt{2} |I_{x}|^B }{ \sqrt{ n_{ I_x} }  } \lesssim \frac{ \ell_x^{B-1/2} }{ \sqrt{n} L_n^{B-1/2}} 
$$
where $\ell_x $ is the number of units (i.e. subintervals at the finest level) in $I_x$, i.e. $\ell_x \asymp L_n |I_x| $. 

Hence, there exists $\gamma<0$ such that for any  $u_n=o(1)$,  writing $I_{x,1} = \bar I\setminus  I_{x}$ and using Markov inequality, 
\begin{align*}
 P \left( \Pi( B_{n,1 }| D_n )> u_n \right)  
 &\lesssim  o(1/n)+ \sum_{ \bar I: x \in \bar I} 
 P\left[ \sum\limits_{I_x\subset \bar I} \1_{{n_{I_x} \leq s_n(\delta)}}  \ell_x^{ B-1/2}
Z_b^{I_x}> \frac{\gamma u_n L_n^{B-1/2} \sqrt{n} }{  n^{ bA+ \delta M_0^2}}  \right] \\
&\lesssim  o(1/n)+ \frac{ n^{ bA+\delta M_0^2 } }{ \sqrt{b}u_n L_n^{B-1/2} \sqrt{n}}  \sum_{ \bar I: x \in \bar I} 
 \sum_{l_x, I_{x,1}} \1_{n_{I_x} \leq  s_n(\delta)}  \ell_x^{ B-1/2} \\
&\lesssim o(1/n)+ \frac{ n^{ bA+\delta M_0^2 } L_n }{ u_n L_n^{B-1/2} \sqrt{n}} \sum_{\ell_x=1}^{L_n s_n(\delta)/n} \ell_x^{B+1/2} \lesssim  o(1/n)+\frac{n^{ bA+\delta M_0^2 }  [L_n s_n(\delta)/n]^{B+3/2} }{ u_n L_n^{B-3/2} \sqrt{n}} \\
& \lesssim  o(1/n)+ (\log n)^q n^{ bA+\delta M_0^2]} n^{\frac{ 5 t(x) + 1 - B}{ 2t(x)+1} }  = o(1/n)
\end{align*}
for some $q>0$  as soon as $ B > 9$,  by choosing $b, \delta$ small enough.

\smallskip
We now study $B_{n,2}$. When $n_{I_{x,1}} < s_n(\delta_1)$ with $\delta_1 \geq \delta$ we have $|I_x\cup I_{x,1}| \leq p_1  s_n(\delta+\delta_1)/n$ and by the H\"older condition on $f_0$ at $x$ we obtain for some $M>0$
$$ 
|\bar \beta_{0, I_x} - \bar \beta_{0, I_{x,1}}| \leq 2M [p_1  s_n(\delta+\delta_1)/n]^{t(x)}
$$
so that
$$ 
\bar y_{I_x}-\bar y_{I_{x,1}} = \bar \epsilon_{I_x} - \bar \epsilon_{I_{x,1}} +  \bar \beta_{0, I_x} - \bar \beta_{0, I_{x,1}} = \bar \epsilon_{I_x} - \bar \epsilon_{I_{x,1}} +O( \delta_1^{ t(x) } \varepsilon_n). 
$$
Consider the event 
$$ 
\bar \Omega_{n,2} = \left\{    \forall \bar I\,\, s.t.\,\, n_{\bar I} \leq s_n(\delta+\delta_1)\,\,\text{and}\,\, x \in \bar I \, : \frac{ \sqrt{n_{I_x}}\sqrt{n_{I}} |\bar \epsilon_{I_x} - \bar \epsilon_{I}| }{ \sqrt{n_{I_x}+n_I} } \leq u_1' \sqrt{ \log n} \right\} . 
$$
Then since for each $(\bar I, I_x)$ , $\sqrt{n_{I_x}}\sqrt{n_{I}} (\bar \epsilon_{I_x} - \bar \epsilon_{I })/ \sqrt{n_{I_x}+n_I} \sim \mathcal N(0,1)$ and since the number of $(\bar I, I_x)$ satisfying $\bar I =  I_x\cup I_{x,1}$ and
$ n_{I_x}, n_{I_{x,1}} \leq s_n(\delta)$ is bounded by 
$$ \sum_{\ell_x}^{ c L_n s_n(\delta)/n}	\ell_x \underbrace{  L_n s_n(\delta_1)/n }_{\text{bound on number of  }  I_{x,1}} \lesssim  [L_n s_n(\delta_1)/n]^3$$
 as soon as $(u_1')^2 > 6t(x)/(2t(x) +1)+2$ we have $P(\bar \Omega_{n,2} ) =1 +o(1/n) $.
On $\bar \Omega_{n,2}$, 
\begin{align*}
\frac{ n_{I_{x,1}} n_{I_x}}{ 2n_{\bar I}  } (\bar y_{I_x}-\bar y_{I_{x,1}})^2 &\leq \frac{ n_{I_{x,1}} n_x}{ 2n_{\bar I}  }(\bar \epsilon_{I_x}-\bar \epsilon_{I_{x,1}})^2 + a(\delta_1) \log n \\
\end{align*}
for some  $a(\delta_1)>0$ which goes to 0 when $\delta_1$ goes to 0 and similarly to before,  for all $1>b>0$
\begin{align*}
\mathbb P (\Pi( B_{n2}  | D_n )> u_n )
 &\leq 
 \mathbb P \left(\frac{ n^{b u_1'+a(\delta_1)} }{ \sqrt{n} L_n^{B-1/2}} \max_{\bar I: x\in \bar I, n_ {\bar I} \leq s_n(\delta_1+\delta)}   \sum_{I_x \subset \bar I} \1_{n_ {I_x} \leq s_n(\delta)}
  \ell_x^{ B-1/2} e^{\frac{(1-b) n_{I_{x,1}} n_x}{ 2n_{\bar I}  }(\bar \epsilon_{I_x}-\bar \epsilon_{I_{x,1}})^2}  > u_n \right)\\
&\leq \frac{  n^{b u_1'+a(\delta_1)} }{ \sqrt{b} u_n\sqrt{n} L_n^{B-1/2} } 
   \sum_{\bar I: x \in \bar I} \1_{ n_ {\bar I} \leq 2s_n(\delta_1) } \sum_{I_x\subset I} \1_{n_ {I_x} \leq s_n(\delta_1)}   \ell_x^{ B-1/2}
     \\
     &\leq  \frac{ n^{b u_1'+a(\delta_1)} }{ \sqrt{b}u_n\sqrt{n} L_n^{B-1/2}}  
[L_ns_n(\delta_1)/n]^{ B+5/2} = O( n^{\frac{ -B + 5t(x)}{ 2t(x)+2} +{b u_1'+a(\delta_1)} })=o(1/n)
 \end{align*}
 since $B>9$ by choosing $b, \delta_1$ small enough.

Finally, we study $B_{n,3}$. 
Since $|I_x| \leq p_1 s_n(\delta)/n$ we can choose  a point in the grid $x_1 \in I_{x,1}$, such that $|x-x_1| \leq 2 p_1 s_n(\delta)/n$, so that the H\"older condition of $f_0$ at $x$ implies that 
$$ |f_0(x)-f_0(x_1)| \leq M (2p_1)^{t(x)} \delta^{t(x) } \varepsilon_n(x) .$$ 
Moreover, since $t$ is H\"older $\alpha$ for some $\alpha>0$ on $[x,x_1]$ (note that for $n$ large enough $|x-x_1|$ is arbitrarily small)  we have
$$
|t(x_1)- t(x) | \leq L_0 \delta^{\alpha} (n/\log n)^{- \alpha/(2t(x)+1)}
$$ 
and $\varepsilon_n(x_1):=(n/\log n)^{-2t(x_1)/(2t(x_1)+1)} = (n/\log n)^{-2t(x)/(2t(x)+1)} ( 1 +o(1) )$. Hence, choosing $M_0 > 2(3p_1)^{t(x)}$,  for $n$ large enough, 
\begin{align*}
\Pi & \left(  B_{n,3} \C D_n \right)  \\
& \leq \Pi  \left(  \{ | \bar y_{I_{x,1}} -f_0(x_1) | >M_0 \varepsilon_n(x_1)/2 \} \cap \{ n_{I_{x,1}} > s_n(\delta_1) \} |D_n \right)= o_P(1/n) 
\end{align*}
from Lemma \ref{lem:mIx} and Theorem \ref{th:repulsive} is proved by choosing $M_0> 4/\sqrt{\delta_1}$.  

\section{Auxiliary Results} \label{sec:aux}

\subsubsection{Auxiliary Lemmata in the Proof of Theorem \ref{sec:proof_thm_np}}

\begin{lemma}\label{lemma:dif}
We denote with $P_{\mT}= X_{\mT}(X_{\mT}'X_{\mT})^{-1}X_{\mT}'$ the projection matrix and with
  $Z=\|X_j\|^2_2-X_j' P_{\mT^-}X_j$. Then 
\begin{equation}
Y'[P_\mT-P_{\mT^-}]Y=\frac{Y'(I-P_{\mT^-})X_jX_j'(I-P_{\mT^-})Y}{Z}.\label{eq:dif}
\end{equation}
\end{lemma}
\proof 
We can write
$$
(X_\mT'X_\mT)^{-1}=\left(
\begin{matrix}
X_{\mT^-}'X_{\mT^-} & X_{\mT^-}'X_j\\
X_j'X_{\mT^-}&\|X_j\|^2_2
\end{matrix}
\right)^{-1}=
\left(
\begin{matrix}
\Sigma_{\mT^-}+{\Sigma_{\mT^-}X_{\mT^-}'X_jX_j'X_{\mT^-}\Sigma_{\mT^-}}/Z &\quad\,\,- {\Sigma_{\mT^-}X_{\mT^-}'X_j}/Z \\
-{X_j'X_{\mT^-}\Sigma_{\mT^-}}/Z &\quad\,\, 1/{Z}\\
\end{matrix}
\right)
$$
Next, noting that $X_\mT=(X_{\mT^-}, X_j)$
$$
X_\mT\Sigma_\mT=
\left(
\begin{matrix}
X_{\mT^-}\left[\Sigma_{\mT^-}+\frac{\Sigma_{\mT^-}X_{\mT^-}'X_jX_j'X_{\mT^-}\Sigma_{\mT^-}}{Z}\right]-\frac{X_j{X_j'X_{\mT^-}'\Sigma_{\mT^-}}}{Z}, &\quad -\frac{P_{\mT^-}X_j}{Z}+\frac{X_j}{Z}\\
\end{matrix}
\right)
$$
which yields 
\begin{align*}
P_\mT=&P_{\mT^-}+\frac{1}{Z}\left[P_{\mT^-}X_jX_j'P_{\mT^-} - X_j{X_j'P_{\mT^-}}-P_{\mT^-}'X_jX_j'+X_jX_j'\right].
\end{align*}
We then obtain \eqref{eq:dif}.\qed

 \begin{lemma}\label{lemma:bias2}
 Let $X_j$ be the $j^{th}$ column in the matrix $X$ and let $\b_0=(\beta_1^0,\dots,\beta^0_{p})'$ be the vector of multiscale coefficients $\langle \psi_{lk},f_0\rangle$ for $f_0\in\mathcal C(t,M,\eta)$ where $t,M,\eta$ satisfy Assumption \ref{ass:local_holder} with $t_1> 1/2$.
 Then, on the event $\mathcal A$, we have 
 $$ 
 |X_j'\bm\nu|=|X_j'(F_0-X\b_0+\bm \varepsilon)|\lesssim\sqrt{n}\log^{1\vee\upsilon}n.
 $$
 \end{lemma}
 \proof
 From the definition of the set $\mA$ in \eqref{eq:setA_np} we know that $|X_j'\bm \varepsilon|\lesssim \sqrt{n}\log n$.
 Next, we decompose the bias term $|X_j'(F_0-X\b_0)|$ into resolutions $\Lmax<l\leq \wtLmax$ that are within the spam of the matrix $X$ and  higher resolutions  $l>\wtLmax$ for which the balancing Assumption \ref{ass:design} is no longer required. Then, using \eqref{eq:covariance}, we obtain
\begin{align*}
 |X_j'(F_0-X\b_0)|\leq&\,\, C_d\sqrt{n}\log^\upsilon n\sum_{l=\Lmax+1}^{\wtLmax}2^{l-l_1}2^{l_1/2}2^{-l(t_1+1/2)} \\
 &+\|X_j\|_1\left\|\sum_{l>\wt L_{max}}\sum_k\psi_{lk}(x)\beta_{lk}^0\right\|_\infty.
 \end{align*}
 The first term above can be bounded by a constant multiple of $2^{-l_1/2}\sqrt{n}\,  \log^\upsilon n$ when $t_1> 1/2$. Regarding the second term,
under the assumption $t_1> 1/2$ and using the fact that $\wtLmax=\mathcal O[\log_2(n/\log n)]$, we obtain for each $x\in[0,1]$ 
\begin{equation*}\label{eq:bias3}
\left|\sum_{l>\wt L_{max}}\sum_k\psi_{lk}(x)\beta_{lk}^0\right|\leq \sum_{l> \wtLmax}2^{l/2}|\beta_{lk_l(x)}^0|\leq \sum_{l> \wtLmax}2^{-lt(x)}\lesssim 2^{-\wtLmax/2}\lesssim \sqrt{\frac{\log n}{n}}.
\end{equation*}
Using \eqref{eq:design_norm} we find that $\|X_j\|_1\lesssim n$. Going back to \eqref{eq:X_bias} to we conclude that   
\begin{equation*}\label{eq:cross_term}
 |X_j'\bm\nu|\lesssim \sqrt{n}\log^\upsilon n  +\sqrt{n}\log n\lesssim \sqrt{n}\log^{1\vee\upsilon}n.\qedhere
\end{equation*}

\subsubsection{Proof of Lemma \ref{lemma:eigen}}\label{sec:proof_lemma_eigen}
The diagonal elements of $X'X$, denoted with  $a(i)$, satisfy
$$
2\,c\,n\leq a(i)\equiv\|X_i\|_2^2\leq 2n[C+\log_2\lfloor C_x\sqrt{n/\log n}\rfloor].
$$ 
For a given node $(l_1,k_1)\in\mT$ with  $i=2^{l_1}+k_1$ we denote with $a(\backslash  i)$  the sum of absolute off-diagonal terms in the $i^{th}$ row of $X'X$.
Under the Assumption \ref{ass:design}, we show below (for some $C_m>0$)
\begin{equation}\label{eq:off_diagonal}
a(\backslash i)=\sum_{(l,k)\neq (l_1,k_1)}|X_i'X_{2^l+k}|\leq C_m\,n\log^{\upsilon-1/2}n.
\end{equation}
In order to show \eqref{eq:off_diagonal}, we split the sum into nodes $(l,k)\in P(l,k)$ that are predecessors of $(l_1,k_1)$ and nodes $(l,k)\in D(l,k)$ that are descendants of $(l_1,k_1)$.
Using  \eqref{eq:covariance} and the fact that there are $2^{l-l_1}$ descendants at each layer $l>l_1$ we have 
\begin{align*}
\sum_{(l,k)\in D(l_1,k_1)}|X_i'X_{2^l+k}|&\leq 
C_d {\sqrt{n}\log^\upsilon n}\,\sum_{l=l_1+1}^{L_{max}}2^{l-l_1}\,2^{\frac{l_1}{2}}\lesssim  {\sqrt{n}\log^\upsilon n}2^{\Lmax}\lesssim
n\log^{\upsilon-1/2}n\end{align*}
and
$$
\sum_{(l,k)\in P(l_1,k_1)}|X_i'X_{2^l+k}|\leq C_d {\sqrt{n}\log^\upsilon n}\, \sum_{l=0}^{l_1-1} 2^{\frac{l}{2}}=
 {\sqrt{n}\log^\upsilon n}\,2^{l_1/2}\lesssim n.
$$
From \eqref{eq:off_diagonal}  one obtains $a(i)-a(\backslash i)>2n[c - C_m\,\log^{\upsilon-1/2}n]>0$ for some suitable $c>0$ and for $\upsilon\leq 1/2$.  The Gershgorin circle theorem \citep{matrix_book} then yields
\begin{align*}
\min [a(i)-a(\backslash i)]&\leq \lambda_{min}(X'X)\leq \lambda_{max}(X'X)\leq \max [a(i)+a(\backslash i)].\quad\quad\qedhere
\end{align*}

\smallskip

\begin{lemma}\label{lemma:design}
Assume that $x_i\iid U[0,1]$. Then for $t_l= \frac{C_d}{2} \frac{\log^\upsilon n}{\sqrt{n}\,2^{l/2}}$ we have
$$
P\left(|\bar n_{lk}-\munderbar n_{lk}|\leq  2nt_l \quad \forall (l,k)\,\,s.t.\,\,l\leq \wt{L}_{max}\right)=1+o(1).
$$
\end{lemma}
\proof
Under  the uniform random design, both $n_{lk}^R$ and $n_{lk}^L$ are distributed according to  $Bin(n,2^{-(l+1)})$.
We can write
  \begin{align*}
  &P(|\bar n_{lk}-\munderbar n_{lk}|\leq  2nt_l \quad \forall (l,k)\,\,s.t.\,\,l\leq \wt L_{max})\geq \\
  &P(| n_{lk}^R-n2^{-(l+1)}|\leq  nt_l\quad \text{and}\quad |n_{lk}^L-n2^{-(l+1)}|\leq  nt_l \quad\forall (l,k)\,\,s.t.\,\,l\leq \wt L_{max})=\\
  &1-P\left(\cup_{l,k}\{| n_{lk}^R-n2^{-(l+1)}|>  nt_l\}\cup \{ |n_{lk}^L-n2^{-(l+1)}|>  nt_l\} \right).
   \end{align*}
   We show that the probability on the right-hand side above is $o(1)$. For $V=n2^{-(l+1)}(1-2^{-(l+1)})$ we note that $2^{\wtLmax}=\mathcal O(n/\log n)$ and thereby, using the Bernstein inequality with $b=3/(C_d\log n)$,
   \begin{align*}
   \sum_{lk} P\left(| n_{lk}^R-n2^{-(l+1)}|>  nt_l\right)&\leq \sum_{l=0}^{\wtLmax} 2^l\exp\left(-  \frac{n^2t_l^2}{2(V+b/3 nt_l)}\right)\\
   &\leq \sum_{l=0}^{\wtLmax} 2^l\exp\left(-  \frac{C_d^2\log^{2\upsilon} n}{4\times2^l}\times\frac{1}{2^{-l}+b\,C_d/3\log n/({\sqrt{n}2^{l/2}})}\right)\\
   &\leq  \sum_{l=0}^{\wtLmax} 2^l\exp\left( -C_d^2/8\log^{2\upsilon}n\right)=\mathcal O\left(\frac{n}{\log n}\right) \times n^{-C_d^2/8\log^{2\upsilon-1}n}.
   \end{align*}
    The last sum  is $o(1)$ when $\upsilon>1/2$ and when $C_d^2/8\geq1$ when $\upsilon=1/2$. \qedhere

   \begin{lemma}(Random Design)\label{lemma:random_design}
 Assume that the design points $x_i$'s are iid with density $p$ bounded from above and below on $[0,1]$ by $p_1 $ and $p_0$, respectively. Then
 \begin{equation}\label{eq:bern1}
 \mathbb P\left( |n_I - n\times p(I)| \leq u_1\sqrt{np(I)\log n } \right) \leq     \e^{- \gamma u_1\log n}, \quad \text{where}\quad \gamma =  3/4 \sqrt{ p_0 C_1},
 \end{equation}
 when $u_1$ is chosen large enough.  
  \end{lemma}
  \proof
 The Bernstein inequality (Theorem 2.8.4 in \cite{bernstein}) implies that for all possible  intervals $I$ in the partition,  we have 
 $$
 \mathbb P\left( |n_I - n\times p(I)| \leq u_1\sqrt{n p(I)\log n } \right) \leq    2\exp\left\{-\frac{ n p(I)u_1^2 \log n }{2np(I)(1-p(I)) + 2/3 u_1\sqrt{n p(I)\log n} } \right\}
 $$
For  $p(I)\geq p_0 C_1 \sqrt{\log n /n}$ we obtain
 $$ 
 np(I)(1-p(I)) + u_1\sqrt{n p(I)\log n }/3 \leq  np(I) \left[ 1 + \frac{ u_1 }{3\sqrt{  p_0 C_1 } }\right].
 $$
With $u_1\geq 3\sqrt{p_0C_1}$ we obtain the desired statement \eqref{eq:bern1}.
\qedhere

\bibliographystyle{imsart-number}

\clearpage


\end{document}